\documentclass[a4paper,12pt]{amsart}
\usepackage{amssymb}
\usepackage{amsmath}
\usepackage{ifthen}
\usepackage{hyperref}
\usepackage{setspace}
\nonstopmode \numberwithin{equation}{section}
\newtheorem{thm}[equation]{Theorem}
\newtheorem{cor}[equation]{Corollary}
\newtheorem{lem}[equation]{Lemma}
\newtheorem{prop}[equation]{Proposition}

\newenvironment{pf}[1][]{%
 \vskip 3mm
 \noindent
 \ifthenelse{\equal{#1}{}}%
  {{\slshape Proof. }}%
  {{\slshape #1.} }%
 }%
{\qed\bigskip}

\begin{document}
\bibliographystyle{amsplain}
\newcommand{\A}{{\mathcal A}}
\newcommand{\B}{{\mathcal B}}
\newcommand{\es}{{\mathcal S}}
\newcommand{\IR}{{\mathbb R}}
\newcommand{\RR}{{\mathcal R}}
\newcommand{\IC}{{\mathbb C}}
\newcommand{\IN}{{\mathbb N}}
\newcommand{\K}{{\mathcal K}}
\newcommand{\uhp}{{\mathbb H}}
\newcommand{\Z}{{\mathbb Z}}
\newcommand{\N}{{\mathcal N}}
\newcommand{\M}{{\mathcal M}}
\newcommand{\SCC}{{\mathcal{SCC}}}
\newcommand{\CC}{{\mathcal C}}
\newcommand{\st}{{\mathcal{SS}}}
\newcommand{\D}{{\mathbb D}}
\newcommand{\remark}{\vskip .3cm \noindent {\sl Remark.} \@}
\newcommand{\remarks}{\vskip .3cm \noindent {\sl Remarks.} \@}
\newcommand{\UCV}{{\mathcal UCV}}
\newcommand{\UST}{{\mathcal UST}}
%%%%%%%%%%%%%%%%%%%%%%%%%%%%%%%%%%%%%%%%%%%%%%%%%%%%%%%%%%%%%%%%%%%%%%%%%%%%%%%%%5
%%%%%%%%%%%%%%%%%%%%%%%%%%%%%%%%%%%%%%%%%%%%%%%%%%%%%%%%%%%%%%%%%%%%%%
\def\be{\begin{equation}}
\def\ee{\end{equation}}
\newcommand{\bee}{\begin{enumerate}}
\newcommand{\eee}{\end{enumerate}}
\newcommand{\pays}{\!\!\!\!}
\newcommand{\pay}{\!\!\!}
\newcommand{\blem}{\begin{lem}}
\newcommand{\elem}{\end{lem}}
\newcommand{\bthm}{\begin{thm}}
\newcommand{\ethm}{\end{thm}}
\newcommand{\bcor}{\begin{cor}}
\newcommand{\ecor}{\end{cor}}
\newcommand{\beg}{\begin{example}}
\newcommand{\eeg}{\end{example}}
\newcommand{\begs}{\begin{examples}}
\newcommand{\eegs}{\end{examples}}
\newcommand{\bdefe}{\begin{defn}}
\newcommand{\edefe}{\end{defn}}
\newcommand{\bprob}{\begin{prob}}
\newcommand{\eprob}{\end{prob}}
\newcommand{\bcon}{\begin{conj}}
\newcommand{\econ}{\end{conj}}
`\newcommand{\bprop}{\begin{prop}}
\newcommand{\eprop}{\end{prop}}
\newcommand{\bpf}{\begin{pf}}
\newcommand{\epf}{\end{pf}}
\newcommand{\ba}{\begin{array}}
\newcommand{\ea}{\end{array}}
\newcommand{\beq}{\begin{eqnarray}}
\newcommand{\beqq}{\begin{eqnarray*}}
\newcommand{\eeq}{\end{eqnarray}}
\newcommand{\eeqq}{\end{eqnarray*}}

\title[Geometric Properties of Generalized Hypergeometric Functions]{Geometric Properties of Generalized Hypergeometric Functions}

\author[K.Chandrasekran]{K. Chandrasekran}
\address{K. Chandrasekran \\ Department of Mathematics \\ MIT Campus, Anna University \\ Chennai 600 044, India}
\email{kchandru2014@gmail.com}

\author[D. J. Prabhakaran ]{D.J. Prabhakaran}
\address{D.J. Prabhakaran \\ Department of Mathematics \\ MIT Campus, Anna University \\ Chennai 600 044, India}
\email{asirprabha@gmail.com}

\subjclass[2000]{30C45}
\keywords{Generalized Hypergeometric Function, Univalent Functions, Starlike Functions, Convex Functions, Uniformly Starlike Functions and Uniformly Convex Functions.}
\maketitle
\begin{abstract}
In this article, Using Hadamard product for $_4F_3\left(^{a_1,\, a_2,\, a_3,\, a_4}_{b_1,\, b_2,\, b_3};z\right)$ hypergeometric function with normalized analytic functions in the open unit disc, an operator $\mathcal{I}^{a_1,a_2,a_3,a_4}_{b_1,b_2,b_3}(f)(z)$ is introduced.  Geometric properties of $_4F_3\left(^{a_1,\, a_2,\, a_3,\, a_4}_{b_1,\, b_2,\, b_3};z\right)$  hypergeometric functions are discussed for various subclasses of univalent functions. Also, we consider an operator $\mathcal{I}^{  a,\frac{b}{4},\frac{b+1}{4},\frac{b+2}{4},\frac{b+3}{4} }_{ \frac{c}{4}, \frac{c+1}{4}, \frac{c+2}{4},\frac{c+3}{4} }(f)(z)$$= z\, _5F_4\left(^{a,\frac{b}{4},\frac{b+1}{4},\frac{b+2}{4},\frac{b+3}{4}}_{\frac{c}{4}, \frac{c+1}{4}, \frac{c+2}{4},\frac{c+3}{4}}; z\right)*f(z)$, where, $_5F_4(z)$ hypergeometric function and the $*$ is usual Hadamard  product. In the main results, conditions are determined on $ a,b,$ and $c$ such that the function $z\, _5F_4\left(^{a,\frac{b}{4},\frac{b+1}{4},\frac{b+2}{4},\frac{b+3}{4}}_{\frac{c}{4}, \frac{c+1}{4}, \frac{c+2}{4},\frac{c+3}{4}}; z\right)$ is in the each of the classes  $ \es^{*}_{\lambda} $, $ \mathcal{C}_{\lambda}$,  ${\UCV} $ and $\es_p$. Subsequently, conditions on $a,\,b,\,c,\, \lambda,$ and $\beta$ are determined using the integral operator such that functions belonging to $\mathcal{R}(\beta)$ and $\es$ are mapped onto each of the classes $ \es^{*}_{\lambda}  $, $ \mathcal{C}_{\lambda}$,  ${\UCV},$ and $\es_p$.
\end{abstract}
\begin{sloppypar}
\section{Introduction and Preliminaries}
\begin{sloppypar}
The classical Complex Analysis is one of the important branches of Mathematics. Geometric function theory is a branch of Complex Analysis, in which there are several results based on the geometric properties of analytic functions. In particular, the theory of Univalent functions remains an active field of current research. Researchers from all over the world have contributed to the rapid progress of this field.%

Moving forward, the generalized hypergeometric function is discussed. The importance of special functions in the field of Geometric function theory was felt after \cite{ch1de-Branges-1985} proved the famous Bieberbach conjecture using hypergeometric functions. He used  \cite{ch1Askey-Gasper-1976} inequality for Jacobi polynomial which involves the generalized hypergeometric functions
$$ _3F_2\left(\begin{array}{c}
                       k+\frac{1}{2},\, k-n,\, n+k+2 \\
                       2k+1, \, k+\frac{3}{2}
                     \end{array}; e^{-t}\right)\,
{and}\,
 _4F_3\left(\begin{array}{c}
                       k+\frac{1}{2},\, n+k+2,\,k,\, k-n,\,  \\
                       k+1,\, 2k+1, \, k+\frac{3}{2}
                     \end{array}; e^{-t}\right).$$
From that time onwards, many eminent scholars in the field of Geometric function theory have established many results in connection with special functions and univalent functions. Many research articles were established in this field  using Gaussian hypergeometric functions, whereas very few research articles are available on generalized hypergeometric functions. For more details \cite{ch1Koepf-2007,ch1Koepf-Schmersau-1996}. In this thesis, many results are obtained relating to generalized hypergeometric functions and univalent functions. Some basic definitions and theorems on generalized hypergeometric functions which are relevant to present study are listed. \cite{ch1Andrews-Askey-Roy-1999-book}
For any non-zero complex variable $a$, the Pochhammer symbol (or shifted factorial) is defined as
\beqq
(a)_n&=&\left\{ \begin{array}{lll}
 a(a+1)(a+2)\cdots(a+n-1),\, n>1 \\
 1,  \,     n=0&
\end{array} \right.
\eeqq
In terms of  Euler gamma function,  the Pochhammer symbol can be written as $$(a)_n = \frac{\Gamma(n+a)}{\Gamma(a)},\quad n=0,1,2,\cdots$$ where $a$ is neither zero nor a negative integer. Further some of the identities related to above definition are stated.%
\blem \label{ch1Lemma1} \cite{ch1Rain-1960-Mac}
For any non-zero complex variable $a$, $k$ is a positive integer and $n$ is a non-negative integer,\cite{ch1Rain-1960-Mac}  is defined as
\beqq
 (a)_{kn}\,=\, k^{kn}  \left(\frac{a}{k}  \right)_{n}   \left(\frac{a+1}{k}  \right)_{n} \cdots  \left(\frac{a+k-1}{k}  \right)_{n}.
\eeqq
\elem
In particular, for the case $ k=2 , 3 $ we have the following%
\begin{itemize}
  \item[i.]  $ (a)_{2n}\,=\, 2^{2n}  \left(\frac{a}{2}  \right)_{n}   \left(\frac{a+1}{2}  \right)_{n}$.
  \item[ii.] $ (a)_{3n}\,=\, 3^{3n}  \left(\frac{a}{3}  \right)_{n}   \left(\frac{a+1}{3}  \right)_{n}  \left(\frac{a+2}{3}  \right)_{n}$.
\end{itemize}%

\cite{ch1Andrews-Askey-Roy-1999-book} The generalized hypergeometric series %\cite{Andrews-Askey-Roy-1999-book,Slater-1966-book}
is defined by
\beqq
_pF_q\left(
\begin{array}{ccc}
  a_{1}, & a_{2},  \cdots  & a_{p} \\
  b_{1}, & b_{2},  \cdots & b_{q}
\end{array} ; \displaystyle z
\right) =\displaystyle\sum_{n=0}^{\infty}\frac{(a_{1})_n \cdots (a_{p})_n}{(b_{1})_n \cdots (b_{q})_n  (1)_{n}} z^n.
\eeqq

This series converges absolutely for all $z$ if $p\leq q$ and for $|z|<1$ if $p=q+1$, and it diverges for all $z\neq 0$ if $p > q+1$.
For $|z|=1$ and $p = q+1$, the series $_{p}F_{q}(a_{1}\ldots,a_{p};b_{1}\ldots b_{q};z)$ converges absolutely if $Re(\sum b_i - \sum a_i) >0.$
The series converges conditionally if $z=e^{i\theta}\neq 1$ and $-1 < Re(\sum b_i - \sum a_i) \leq 0$ and diverges if $Re(\sum b_i - \sum a_i)\leq -1.$%
\bthm\label{thm1eq0} \cite{ch1Rain-1960-Mac}  For  $ p \leq q+1 $, if   $Re(b_{1})\, >\, Re(a_{1})\, > 0$, if none of $  b_{1},  b_{2},  \cdots  b_{q}$ is zero or a negative integer and if\, $ | z | < 1  $, then
\beqq\label{inteq6}
\displaystyle  _pF_q\left(
\begin{array}{ccc}
  a_{1}, & a_{2},  \cdots  & a_{p} \\
  b_{1}, & b_{2},  \cdots & b_{q}
\end{array} ; \displaystyle z
\right) &=& \frac{\Gamma(b_{1})}{\Gamma(a_{1})\,\Gamma(b_{1} - a_{1} )}\int_{0}^{1} t^{a_{1}-1} (1-t)^{b_{1}-a_{1} -1} \,\nonumber\\
 && \times
_{p-1}F_{q-1}\left(
\begin{array}{ccc}
   a_{2},  \cdots  & a_{p} \\
   b_{2},  \cdots & b_{q}
\end{array} ; \displaystyle zt
\right) \, dt.
\eeqq
If $ p \leq q$, the condition $ | z | < 1  $ may be omitted.
\ethm

If $p=2$ and $q=1$, the generalized hypergeometric function becomes the Gaussian hypergeometric function. Further, the formal definition of Gaussian hypergeometric function and related results are stated. \cite{ch1Andrews-Askey-Roy-1999-book}  The Gaussian hypergeometric function is defined as
\beqq\label{inteq003}
 _2F_1(a,b;c;z) = \sum^{\infty}_{n=0} \frac{(a)_n(b)_n}{(c)_n(1)_n} z^n, \quad |z| < 1
\eeqq
where $a,b,c \in {\mathbb C}$ and  $c$ is neither zero nor a negative integer.
\bthm\label{thm1eq0} \cite{ch1Andrews-Askey-Roy-1999-book}  For  $Re(c)\, >\, Re(b)\, > 0$,
\beqq\label{inteq6}
_2F_1(a,b;c;z)&=& \frac{\Gamma(c)}{\Gamma(b)\,\Gamma(c-b)}\, \int_{0}^{1} t^{b-1} (1-t)^{c-b-1} ( 1 - z t )^{-a} \, dt.
\eeqq
\ethm
It is interesting to see that at $z=1$ Gauss proved following summation formula:
\bthm\label{thm1eq0} \cite{ch1Andrews-Askey-Roy-1999-book}  If $Re(c)\, >\, Re(b)\, > 0$ and $Re(c-a-b)\, >\, 0$, then
\beqq
  _2F_1(a,b;c;1) = \frac{\Gamma (c)\,\Gamma (c-a-b) }{\Gamma (c-a) \Gamma (c-b)}.
\eeqq
\ethm%
The Gaussian hypergeometric function $_2F_1(a,b;c;z)$  has been studied extensively by a number of authors in the context of
deriving geometric properties such as convexity, starlikeness, close-to-convexity and univalency. These were obtained by determining conditions on the parameters $a,\,b,\,c.$\\

If $p=3$ and $q=2$, the generalised hypergeometric function  is known as Clausen's hypergeometric function. Originally, Clausen obtained the series by squaring the Gaussian hypergeometric series
$$\displaystyle _2F_1(a,b;a+b+1/2; z)^2 \,=\,  _3F_2(2a,2b,a+b;2a+2b,a+b+1/2;z).$$
\cite{ch1Slater-1966-book}
The Clausen's hypergeometric function $_3F_2(a,b,c;d,e;z)$ is defined as
\beqq
_3F_2(a,b,c;d,e;z)=\sum_{n=0}^{\infty}\frac{(a)_n(b)_n(c)_n}{(d)_n(e)_n(1)_n}z^n, \quad |z| < 1,
\eeqq
where $ a,b,c,d,e\in \IC$  with  $d,\, e\, \neq 0,-1,-2,-3\cdots$.%
\bthm\label{thm1eq0} \cite{ch1Slater-1966-book} If $Re(d)\, >\, Re(b)\, > 0$, if none of $d,\, e$ is zero or a negative integer and if $ | z | < 1  $,
\beqq
 \displaystyle  _3F_2\left(
\begin{array}{ccc}
  a, & b, & c \\
  d, & e
\end{array} ; \displaystyle z
\right) &=& \frac{\Gamma(d)}{\Gamma(b)\,\Gamma(d - b)}\,\int_{0}^{1} t^{b-1} (1-t)^{d-b-1}{}
_{2}F_{1}\left(
\begin{array}{ccc}
   a, & c \\
   e
\end{array}; \displaystyle zt
\right) dt.
\eeqq
\ethm
Recently, \cite{ch1Mill-Paris-2012-ITSF} derived a summation formula involving $n$. Also, \cite{ch1Shpot-Srivas-2015-AMC} found a new summation formula involving $m$ and $n$. The restricted form of their summation formula is as follows:
\bthm\label{thm1eq0} \cite{ch1Shpot-Srivas-2015-AMC}  For  $a,\, b,\, c\, >0$, \, $c \neq b$\, and $a < \min\{1,\, b+1,\, c+1\}$,
\beqq\label{ch1inteq6}
\displaystyle  _3F_2\left(
\begin{array}{ccc}
  a, & b, & c \\
  b+1, & c+1
\end{array} ; \displaystyle 1
\right)
&=& \frac{bc}{c-b}\Gamma(1-a)\left[\frac{\Gamma(b)}{\Gamma(1-a+b)}-\frac{\Gamma(c)}{\Gamma(1-a+c)}\right].
\eeqq
\ethm%
In \cite{ch1Driver-Johnston-2006} derived a summation formula for Clausen's hypergeometric function using the Lemma \ref{ch1Lemma1} and Gauss summation formula, which is given below:
\bthm\label{thm1eq0} \cite{ch1Driver-Johnston-2006}  For  $Re(c)\, >\, Re(b)\, > 0$,
\beqq\label{inteq6}
\displaystyle  _3F_2\left(
\begin{array}{ccc}
  a, & \frac{b}{2}, & \frac{b+1}{2} \\
  \frac{c}{2}, &  \frac{c+1}{2}
\end{array} ; \displaystyle z
\right)&=& \frac{\Gamma(c)}{\Gamma(b)\,\Gamma(c-b)}\,\int_{0}^{1} t^{b-1} (1-t)^{c-b-1} ( 1 - z t^{2} )^{-a} \, dt.  \nonumber
\eeqq
\ethm
\bthm\label{thm1eq0} \cite{ch1Driver-Johnston-2006}  If $Re(c)\, >\, Re(b)\, > 0$ and $Re(c-a-b)\, >\, 0$, then
\beqq\label{ch1inteq06}
\displaystyle  _3F_2\left(
\begin{array}{ccc}
  a, & \frac{b}{2}, & \frac{b+1}{2} \\
  \frac{c}{2}, &  \frac{c+1}{2}
\end{array} ; \displaystyle 1
\right)
&=& \frac{\Gamma(c)\,\Gamma(c-a-b)}{\Gamma(c-a)\,\Gamma(c-b)}\, _{2}F_1(a,b;c-a;-1).
\eeqq
\ethm
Clausen's hypergeometric function $_3F_2(a,b,c;d,e;z)$ has been studied by only few authors in connection with Geometric function theory. In particular, \cite{ch1Ponnu-saba-1997} considered the generalized hypergeometric functions and tried to find conditions
on the parameters so that $z\,_3F_2(a,b,c;d,e;z)$ has some geometric properties.\\

If $p=4$ and $q=3$, the generalised hypergeometric function is reduced to the following: \cite{ch1Slater-1966-book}
The hypergeometric function is defined as
\beqq\label{inteq5}
\displaystyle  _4F_3\left(
\begin{array}{cccc}
  a_{1}, & a_{2},& a_{3}, & a_{4} \\
  b_{1}, & b_{2}, & b_{3}
\end{array} ; \displaystyle z
\right)=\sum_{n=0}^{\infty}\frac{(a_1)_n(a_2)_n(a_3)_n(a_4)_n}{(b_1)_n(b_2)_n(b_3)_n(1)_n}z^n, \quad |z| < 1,
\eeqq
where, $ a_1,a_2,a_3,a_4,b_1,b_2,b_3\in \IC$  with $b_1,\, b_2,\, b_3\, \neq 0,-1,-2,-3\cdots$.
\bthm\label{thm1eq0} \cite{ch1Driver-Johnston-2006}  If   $Re(b_{1})\, >\, Re(a_{1})\, > 0$, if none of $  b_{1},  b_{2}, b_{3}$ is zero or a negative integer and if $ | z | < 1  $,
\beqq\label{inteq6}
\displaystyle  _4F_3\left(
\begin{array}{cccc}
  a_{1}, & a_{2},& a_{3}, & a_{4} \\
  b_{1}, & b_{2}, & b_{3}
\end{array} ; \displaystyle z
\right)&=& \frac{\Gamma(b_{1})}{\Gamma(a_{1})\,\Gamma(b_{1} - a_{1} )}\, \int_{0}^{1} t^{a_{1}-1} (1-t)^{b_{1}-a_{1} -1} \,\nonumber\\
&& \qquad \times
_{3}F_{2}\left(
\begin{array}{ccc}
   a_{2}, & a_{3}, & a_{4} \\
   b_{2}, & b_{3}
\end{array} ; \displaystyle zt
\right) \, dt.
\eeqq
\ethm
\end{sloppypar}
\begin{sloppypar}
The action of  some important integral operators on various classes have been studied by many authors. These operators are listed below:

\cite{ch1Alexander-1915}
For $f\in {\A},$ $$F_{0}(z) = \int_{0}^{1}{\frac {f(tz)}{t}}dt = \int_{0}^{z}{\frac {f(t)}{t}}dt,$$ is known as Alexander transform.

Subsequently,   \cite{ch1Libera-1965} introduced an operator known as Libera's operator and showed that the classes $ {\K, \es^*}$ and $\CC$ are closed under this operator.

For $f\in {\A},$  $$F_{1}(z) =2{\int_{0}^{1}{f(tz)}dt} ={\frac{2}{z}}\int_{0}^{z}{f(t)}dt,$$ is defined as Libera's integral operator.%

\cite{ch1Bernardi-1969} introduced an operator in such a way that Alexander and Libera transforms become special cases of his operator.%

Let $f$ be in ${\A}.$ Then $$F_{c}(z) = \frac{c + 1}{z^{c}} \int_{0}^{z} t^{c - 1} f(t)dt,$$
with $c\geq 0,$ is called Bernardi's integral operator.

It is noted that, $c=0$ gives the Alexander transform and $c=1$ gives the Libera's operator. Bernard proved that the classes ${\K, \es^*, \CC}$ are closed under his operator.

\cite{ch1Hohlov-1889}
For $f\in {\A},$  the Hohlov  operator is defined as   $I_{a,b,c} (f)$ by $$[I_{a,b,c} (f)] (z) = z {}_2 F_1 (a,b;c;z) \ast f(z),$$ where $\ast$ stands for the usual Hadamard product.%

Using the integral representation of the Gaussian hypergeometric function, one can write $$[I_{a,b,c} (f)] (z) = \frac{\Gamma (c)}{\Gamma (b) \Gamma(c-b)} \int^{1}_{0} t^{b-1} (1-t)^{c-b-1} \frac{f(tz)}{t} dt \ast \frac{z}{(1-z)^a}.$$
It is observed that the Alexander, Libera and Bernardi transforms can be obtained as special cases of Hohlov convolution operator.%

Now, the work carried out by us will be briefed here. The integral operators using generalized hypergeometric function are introduced and conditions are imposed on the parameters such that the generalized hypergeometric function lies in the each of class $ \es^{*}_{\lambda}$, $ \mathcal{C}_{\lambda}$,  ${\UCV} $ and $\es_p$. Subsequently, conditions on the parameters are determined using integral operators, so that the functions belonging to $\mathcal{R}(\beta)$ and $\es$ are mapped onto each of the classes $ \es^{*}_{\lambda}  $, $ \mathcal{C}_{\lambda}$,  ${\UCV} $ and $\es_p$.

The integral operator $  \mathcal{I}^{a_1,a_2,a_3,a_4}_{b_1,b_2,b_3}(f) $ is introduced  in second section using the  hypergeometric function
$ _4F_3\left(a_1,a_2,a_3,a_4;b_1,b_2,b_3; z\right)$. Condition is determined on $ a,b, c, $ such that the function $ z\, _4F_3\left(a,\frac{b}{3},\frac{b+1}{3},\frac{b+2}{3};\frac{c}{3}, \frac{c+1}{3}, \frac{c+2}{3}; z \right)  $ is in the each of the classes  $ \es^{*}_{\lambda}  $, $ \mathcal{C}_{\lambda}$,  ${\UCV} $ and $\es_p$. And also conditions on $ a,b, c$ are determined using the integral operator such that functions belonging to  $\mathcal{R}(\beta)$ and $\es$  are mapped onto each of the classes $ \es^{*}_{\lambda} $, $ \mathcal{C}_{\lambda}$,  ${\UCV} $ and $\es_p$.%

The integral operator $\mathcal{I}^{  a,\frac{b}{4},\frac{b+1}{4},\frac{b+2}{4},\frac{b+3}{4} }_{ \frac{c}{4}, \frac{c+1}{4}, \frac{c+2}{4},\frac{c+3}{4} }(f)(z)$ is defined in the third section with usual Hadamard  product. Condition is determined on $ a,b, c$ such that the function $z\, _5F_4\left(^{a,\frac{b}{4},\frac{b+1}{4},\frac{b+2}{4},\frac{b+3}{4}}_{\frac{c}{4}, \frac{c+1}{4}, \frac{c+2}{4},\frac{c+3}{4}}; z\right)$ is in the each of the classes  $ \es^{*}_{\lambda}  $, $ \mathcal{C}_{\lambda}$,  ${\UCV} $ and $\es_p$. Subsequently, conditions on $a,\,b,\,c,\, \lambda$ and $\beta$ are determined using the integral operator such that functions belonging to $\mathcal{R}(\beta)$ and $\es$ are mapped onto each of the classes $ \es^{*}_{\lambda}  $, $ \mathcal{C}_{\lambda}$,  ${\UCV} $ and $\es_p$.%

Finally, we have enlisted only the books and the research papers consisting of the results that are used directly to prove our results in the present study are given at the end of each section.
\end{sloppypar}

\end{sloppypar}
\newpage
\begin{sloppypar}
\section[Univalence, Starlikeness and Convexity Properties of $_4F_3\left(^{a_1,\, a_2,\, a_3,\, a_4}_{\, b_1,\, b_2,\, b_3};z\right)$  Hypergeometric Functions using convolution technique]{Univalence, Starlikeness and Convexity properties of $_4F_3\left(^{a_1,\, a_2,\, a_3,\, a_4}_{\, b_1,\, b_2,\, b_3};z\right)$ Hypergeometric Functions using convolution technique}
\begin{abstract}
Using Hadamard product for $_4F_3\left(^{a_1,\, a_2,\, a_3,\, a_4}_{b_1,\, b_2,\, b_3};z\right)$ hypergeometric function with normalized analytic functions in the open unit disc, an operator $\mathcal{I}^{a_1,a_2,a_3,a_4}_{b_1,b_2,b_3}(f)(z)$ is introduced.  Geometric properties of $_4F_3\left(^{a_1,\, a_2,\, a_3,\, a_4}_{b_1,\, b_2,\, b_3};z\right)$  hypergeometric functions are discussed for various subclasses of univalent functions.
\end{abstract}
\maketitle
\subsection{Introduction and preliminaries}
The integral operator plays an important role to characterize, various subclasses of univalent functions in geometric function theory. The important integral operators are Alexander, Libera, and Bernadi and they are the particular case of the Hohlov convolution operator, which is none other than convolution of the normalized analytic univalent function with Gaussian hypergeometric function. Recently, Chandrasekran and Prabhakaran introduced an integral operator involving the Clausen’s hypergeometric function and using this operator, they derived geometric properties of various subclasses of univalent function. In sequel we are interested in operators that involves more generalized hypergeometric function with more parameters. In 2006, Driver and Johnston paper \cite{ch2Driver-Johnston-2006} introduced an integral representation of $_4F_3\left(^{a_1,\, a_2,\, a_3,\, a_4}_{b_1,\, b_2,\, b_3};z\right)$ hypergeometric function. So, in this article we introduced an integral operator that involving $_4F_3\left(^{a_1,\, a_2,\, a_3,\, a_4}_{b_1,\, b_2,\, b_3};z\right)$ hypergeometric function.\\
Now, we will state the basic concepts of geometric function theory, which will helpful prove our main results.%

Let, $\A$ denote the family of functions $f$ of the form
\beq\label{ch2inteq0}
f(z)= z+\sum_{n=2}^{\infty}\, a_n\,z^n
\eeq that are normalized analytic in the open unit disc $\D =\{z:\, |z|<1\}$.%

The class of univalent, starlike and convex functions denoted by $\es$, \, $\es^{*}$ and $\CC$, respectively. (For more details  \cite{ch2Peter-L-Duren-book-1983,ch2A-W-Goodman-1983-book}). If $\displaystyle f \in \A$ is belongs to $\es$, then,
\beq\label{ch2inteq00}
|a_n|\leq n,\text{ for}\, n\geq 2.
\eeq

We consider the classes $\es^{*}_{\lambda},\,  \lambda >0,$ and $\mathcal{C}_{\lambda}$ are defined as $\mathcal{S}^{*}_{\lambda}\, =\, \left\{f\in \mathcal{A}\, |\, \left|\frac{zf'(z)}{f(z)}-1\right| \, < \, \lambda,\, z\in \D \right\}$ and ${\IC_{\lambda}}=\left\{f\in \mathcal{A}\, |\, zf'(z)\in \mathcal{S}^{*}_{\lambda}\right\},$ respectively.%

A sufficient condition for which the function $f$ to be in $\es^{*}_{\lambda}$  and $\mathcal{C}_{\lambda}$, respectively, are as follows:
\beq\label{ch2inteq2}
\displaystyle \sum_{n=2}^{\infty}(n+\lambda-1)|a_n| \leq \lambda.
\eeq
and
\beq\label{ch2inteq}
\displaystyle \sum_{n=2}^{\infty}\, n\, (n+\lambda-1)|a_n| \leq \lambda.
\eeq

Goodman\cite{ch2Good-1991-Ann-PM,ch2Good-1991-JMAA} introduced the concept of uniformly convex and uniformly starlike functions and denoted as $\UCV$ and $\UST$, respectively. Subsequently, R{\o}nning \cite{ch2Ronn-1993-Proc-ams} and Ma and Minda \cite{ch2Ma-Minda-1992-Ann-PM} independently gave the one variable analytic characterization of the class $\UCV$. %

In 1995, K. G. Subramanian et al. \cite{ch2Subram-Murugu-1995} was proved that if
\beq\label{ch2lem4eq1}
\sum_{n=2}^{\infty}\,n\, (2n-1)|a_n|\leq 1,
\eeq then the function of the form (\ref{ch2inteq0}) is in $\UCV$.%

The  subclass $\es_p$ of starlike functions introduced by R{\o}nning \cite{ch2Ronn-1993-Proc-ams} is defined as
\beqq\label{ch2inteq3}
\displaystyle \es_p = \{ F \in \es^{*}|F(z)=zf'(z),\, f(z) \in \UCV \}.
\eeqq

 A sufficient condition for a function $f$ of form (\ref{ch2inteq0}) to belong to $\es_p$ is given by
 \beq \label{ch2lem2eq1}
\sum_{n=2}^{\infty}(2n-1)|a_n|\leq 1.
\eeq

The more general case $\es_p(\alpha)$ was proved in \cite{ch2Subra-Sudharsan-1998}. Let, $  f(z)= z+\sum_{n=2}^{\infty}\, a_n\,z^n $ and  $g(z)= z+\sum_{n=2}^{\infty}\, b_n\,z^n $ be analytic in $\D$. Then, the Hadamard product or convolution of $f(z)$  and $g(z)$ is defined by $f(z)*g(z)= z+\sum_{n=2}^{\infty} a_nb_n z^n.$%

For $\beta <1$, Let, ${\mathcal R}(\beta)=\{f\in {\mathcal A}:  \exists \  \eta \in \left(-\frac{\pi}{2}, \frac{\pi}{2}\right)\, |\, {\rm Re}\, \ [ e^{i\eta}(f'(z)-\beta)] > 0, \quad z\in\D \}.$ It is evident that ${\mathcal R}(\beta) \subset {\mathcal S}$ when $\beta \ge 0$ and for each $\beta <0 ,\ \  {\mathcal R}(\beta) $ also contains non-univalent functions. For more details, refer \cite{ch2Anbu-Parva-2000,ch2Parva-prabha-2001-Far-East} Suppose that $\displaystyle f \in \A$ is in the class $ {\mathcal R}(\beta)$. Then, by \cite{ch2MacGregor-1962-Trans-ams}, we have
\beq\label{ch2inteq30}
|a_n|\leq \frac{2(1-\beta)}{n},\, n \geq 2.
\eeq

For any non-zero complex variable $a$, the ascending factorial (or Pochhammer symbol) is defined as $(a)_0\,=\,1,\, {\rm and }\, (a)_n\,=\,a(a+1)\cdots (a+n-1),\, {\rm for}\, n\,=\,1,2,3,\cdots.$%

The hypergeometric function $_4F_3(z)=\,  _4F_3\left(^{a_1,a_2,a_3,a_4}_{b_1,b_2,b_3};z\right)$ is defined as
\beqq\label{ch2inteq5}
_4F_3(z)=\sum_{n=0}^{\infty}\frac{(a_1)_n(a_2)_n(a_3)_n(a_4)_n}{(b_1)_n(b_2)_n(b_3)_n(1)_n}z^n; \, \, \, a_1,a_2,a_3,a_4,b_1,b_2,b_3 \in \IC,
\eeqq provided $b_1,\, b_2,\, b_3\, \neq 0,-1,-2,-3\cdots,$ which is an analytic functions in unit disc $\D$.%

The Gaussian Hypergeometric function $_2F_1(a,b;c;z)$ have been studied extensively by a number of authors in the context of deriving geometric properties such as convexity, starlikeness, close-to-convexity, and univalency. The Clausen's hypergeometric function $_3F_2(a,b,c;d,e;z)$ have been studied by only few authors. In \cite{ch2Chandru-prabha-2019}, Chandrasekran and Prabhakaran have introduced an integral operator and derived the geometric properties for the Clausen's hypergeometric Series $z\, _3F_2(a,b,c;b+1,c+1;z)$, in which, the numerator and denominator parameters differs by arbitrary negative integers. Recently, Univalence, Starlikeness and Convexity properties of an integral $\mathcal{I}^{a,\frac{b}{2},\frac{b+1}{2}}_{\frac{c}{2}, \frac{c+1}{2}}(f)(z)$ were obtained based on its Taylor's coefficient of various subclasses of univalent functions by the authors Chandrasekran and Prabhakaran \cite{ch2Chandru-prabha-2020}.%

In this article, the authors introduced an integral operator $\mathcal{I}^{a_1,a_2,a_3,a_4}_{b_1,b_2,b_3}(f)(z)$ involving $_4F_3(z)$ hypergeometric function by means of convolution (or Hadamard  product). i.e.,%

For $f\in \mathcal{A}$, we define the operator $\mathcal{I}^{a_1,a_2,a_3,a_4}_{b_1,b_2,b_3}(f)(z)$
\beq\label{ch2inteq7}
\mathcal{I}^{a_1,a_2,a_3,a_4}_{b_1,b_2,b_3}(f)(z) &=& z\, _4F_3\left(^{a_1,a_2,a_3,a_4}_{b_1,b_2,b_3}; z\right)*f(z)\\& =& z+\sum_{n=2}^{\infty} A_n\, z^n,\nonumber
\eeq
with $A_1=1$ and for $n > 1,$
\beq\label{ch2inteq007}
A_n&=&\frac{(a_1)_{n-1}\left(a_2\right)_{n-1}\left(a_3\right)_{n-1}\left(a_4\right)_{n-1}}{\left(b_1\right)_{n-1}\left(b_2\right)_{n-1}\left(b_3\right)_{n-1}(1)_{n-1}}\, a_n.
\eeq

In 2006, Driver and Johnston  \cite{ch2Driver-Johnston-2006} derived a summation formula for $_4F_3\left(1\right)$ hypergeometric function in terms of Gaussian hypergeometric function. We recall their summation formula as follows:
\begin{sloppypar}
\beq\label{ch2inteq6}
_4F_3\left(1\right)&=& \frac{\Gamma(c)\,\Gamma(c-a-b)}{\Gamma(c-a)\,\Gamma(c-b)}\,\left(\sum_{n=0}^{\infty}\frac{(a)_k\,(-)^k\,(b)_k}{(c-a)_k}\right)\\
&&\qquad\qquad\qquad\qquad\qquad\qquad\qquad\times\,_{2}F_1(-k,b+k;c-a+k;-1)\nonumber
\eeq
where, $_4F_3\left(1\right)=\, _4F_3\left(^{a,\frac{b}{3},\frac{b+1}{3},\frac{b+2}{3}}_{\frac{c}{3}, \frac{c+1}{3}, \frac{c+2}{3}}; 1\right)$ provide $Re(c-a-b)\, >\, 0.$
\end{sloppypar}

In particular, we are focus on the following particular values of the parameter for our investigation in this article.%

\beq\label{ch2inteq7}
\mathcal{I}^{a,\frac{b}{3},\frac{b+1}{3},\frac{b+2}{3}}_{\frac{c}{3}, \frac{c+1}{3}, \frac{c+2}{3}}(f)(z) &=& z\, _4F_3\left(^{a,\frac{b}{3},\frac{b+1}{3},\frac{b+2}{3}}_{\frac{c}{3}, \frac{c+1}{3}, \frac{c+2}{3}}; z\right)*f(z)\\
 &=& z+\sum_{n=2}^{\infty} A_n\, z^n,\nonumber
\eeq
with $A_1=1$ and for $n > 1,$
\beq\label{ch2inteq007}
A_n&=&\frac{(a)_{n-1}\left(\frac{b}{3}\right)_{n-1}\left(\frac{b+1}{3}\right)_{n-1}\left(\frac{b+2}{3}\right)_{n-1}}{\left(\frac{c}{3}\right)_{n-1}\left(\frac{c+1}{3}\right)_{n-1}\left(\frac{c+2}{3}\right)_{n-1}(1)_{n-1}}\, a_n.
\eeq

\subsection{Main Results and Proofs.}

\begin{sloppypar}
The following Lemma is useful to prove our main results.
\blem \label{ch2lem1eq1}
Let, $a,b,c > 0$.
\begin{enumerate}
\item For $ c > a+b+1$.\\
\begin{flushleft}
$\displaystyle\sum_{n=0}^{\infty} \frac{(n+1)\,(a)_n\, \left(\frac{b}{3}\right)_n\, \left(\frac{b+1}{3}\right)_n\, \left(\frac{b+2}{3}\right)_n }
{\left(\frac{c}{3}\right)_n\, \left(\frac{c+1}{3}\right)_n\, \left(\frac{c+2}{3}\right)_n\, (1)_n}$
\end{flushleft}
\begin{eqnarray*}
&=& \frac{\Gamma(c)\, \Gamma(c-a-b)}{\Gamma(c-a)\,\Gamma(c-b)} \, \bigg( \sum_{n=0}^{\infty}\,\bigg(\frac{(a)_{n+1}\,(-1)^n\,(b)_{n+3}}{n!\,(c-a)_{n+2}\,(c-a-b-1)}\bigg)\cr \cr
&& \qquad \times \, _{2}F_1(-n,b+3+n;c-a+2+n;-1)\,\cr \cr
&& +\,\sum_{n=0}^{\infty}\,\bigg( \frac{(a)_n\, (-1)^n\,(b)_{n}\,}{n!\,(c-a)_{n}\,}\bigg)\,_{2}F_1(-n,b+n;c-a+n;-1)\bigg).
\end{eqnarray*}
\item For $ c > a+b+2$,\\
\begin{flushleft}
$\displaystyle\sum_{n=0}^{\infty} \frac{(n+1)^2\,(a)_n\, \left(\frac{b}{3}\right)_n\, \left(\frac{b+1}{3}\right)_n\, \left(\frac{b+2}{3}\right)_n }
{\left(\frac{c}{3}\right)_n\, \left(\frac{c+1}{3}\right)_n\, \left(\frac{c+2}{3}\right)_n\, (1)_n}$
\end{flushleft}
\begin{eqnarray*}
&=&   \frac{\Gamma(c)\, \Gamma(c-a-b)}{\Gamma(c-a)\,\Gamma(c-b)} \,\bigg( \sum_{n=0}^{\infty}\,\bigg(\frac{(a)_{n+2}\,(-1)^n\,(b)_{n+6}}{n!\,(c-a)_{n+4}\,(c-a-b-2)_2}\bigg)\cr \cr
&& \qquad \qquad\times \, _{2}F_1(-n,b+6+n;c-a+4+n;-1)\,\cr \cr
&&+ 3\,\sum_{n=0}^{\infty}\,\bigg(\frac{(a)_{n+1}\,(-1)^n\,(b)_{n+3}}{n!\,(c-a)_{n+2}\,(c-a-b-1)}\bigg)\,\cr \cr
&&\qquad\qquad\times\, _{2}F_1(-n,b+3+n;c-a+2+n;-1)\, \cr \cr
&& +\,\sum_{n=0}^{\infty}\,\bigg( \frac{(a)_n\, (-1)^n\,(b)_{n}\,}{n!\,(c-a)_{n}\,}\bigg)\,_{2}F_1(-n,b+n;c-a+n;-1)\bigg).
\end{eqnarray*}
\item For $ c > a+b+3$,\\
\begin{flushleft}
$\displaystyle\sum_{n=0}^{\infty} \frac{(n+1)^3\,(a)_n\, \left(\frac{b}{3}\right)_n\, \left(\frac{b+1}{3}\right)_n\, \left(\frac{b+2}{3}\right)_n }
{\left(\frac{c}{3}\right)_n\, \left(\frac{c+1}{3}\right)_n\, \left(\frac{c+2}{3}\right)_n\, (1)_n}$
\end{flushleft}
\begin{eqnarray*}
&=&   \frac{\Gamma(c)\, \Gamma(c-a-b)}{\Gamma(c-a)\,\Gamma(c-b)} \,\bigg(  \sum_{n=0}^{\infty}\,\bigg(\frac{(a)_{n+3}\,(-1)^n\,(b)_{n+9}}{n!\,(c-a)_{n+6}\,(c-a-b-3)_3}\bigg)\, \cr \cr
&& \qquad\qquad \times \, _{2}F_1(-n,b+9+n;c-a+6+n;-1)\,\cr \cr
&& + 6 \sum_{n=0}^{\infty}\,\bigg(\frac{(a)_{n+2}\,(-1)^n\,(b)_{n+6}}{n!\,(c-a)_{n+4}\,(c-a-b-2)_2}\bigg)\cr \cr
&&\qquad\qquad \times \, _{2}F_1(-n,b+6+n;c-a+4+n;-1)\,\cr \cr
&& + 7\,\sum_{n=0}^{\infty}\,\bigg(\frac{(a)_{n+1}\,(-1)^n\,(b)_{n+3}}{n!\,(c-a)_{n+2}\,(c-a-b-1)}\bigg)\cr \cr
&&\qquad\qquad \times \, _{2}F_1(-n,b+3+n;c-a+2+n;-1)\,\cr \cr
&& +\,\sum_{n=0}^{\infty}\,\bigg( \frac{(a)_n\, (-1)^n\,(b)_{n}\,}{n!\,(c-a)_{n}\,}\bigg)\,_{2}F_1(-n,b+n;c-a+n;-1)\bigg).
\end{eqnarray*}
\item For $a\neq 1,\, b\neq 1,\, 2,\, 3$ and $c > \max\{a+2,\, a+b-1\}$,\\
\begin{flushleft}
$\displaystyle\sum_{n=0}^{\infty} \frac{(a)_n\, \left(\frac{b}{3}\right)_n\, \left(\frac{b+1}{3}\right)_n\, \left(\frac{b+2}{3}\right)_n }
{\left(\frac{c}{3}\right)_n\, \left(\frac{c+1}{3}\right)_n\, \left(\frac{c+2}{3}\right)_n\, (1)_{n+1}}$
\end{flushleft}
\begin{eqnarray*}
&=& \frac{\Gamma(c)\, \Gamma(c-a-b)}{\Gamma(c-a)\,\Gamma(c-b)} \,\\  \\
&&\qquad\times\, \bigg(\sum_{n=0}^{\infty} \, \frac{(c-a-1)\,(c-a-b)\,(a-1)_n\, (-1)^n\,(b-3)_{n}\,}{n!\, (a-1)\, (b-3)_3\,(c-a-2)_{n}}\cr \cr &&\times \, _{2}F_1(-n,b-3+n;c-a-2+n;-1)-\left(\frac{(c-3)_3}{(a-1)\, (b-3)_3}\right)\bigg).
\end{eqnarray*}
\end{enumerate}
\elem
\end{sloppypar}
\bpf(1) Using Pochhammer symbol, we can formulate
\begin{flushleft}
$\displaystyle \sum_{n=0}^{\infty} \frac{(n+1)\, (a)_n\, \left(\frac{b}{3}\right)_n\, \left(\frac{b+1}{3}\right)_n\, \left(\frac{b+2}{3}\right)_n }{\left(\frac{c}{3}\right)_n\, \left(\frac{c+1}{3}\right)_n\, \left(\frac{c+2}{3}\right)_n\, (1)_n}$
\end{flushleft}
\begin{eqnarray*}
&=& \displaystyle \sum_{n=0}^{\infty} \frac{(a)_{n+1}\, \left(\frac{b}{3}\right)_{n+1}\, \left(\frac{b+1}{3}\right)_{n+1}\, \left(\frac{b+2}{3}\right)_{n+1}}{\left(\frac{c}{3}\right)_{n+1}\, \left(\frac{c+1}{3}\right)_{n+1}\, \left(\frac{c+2}{3}\right)_{n+1}\, (1)_{n}} + \displaystyle \sum_{n=0}^{\infty} \frac{(a)_n\, \left(\frac{b}{3}\right)_n\, \left(\frac{b+1}{3}\right)_n\, \left(\frac{b+2}{3}\right)_n}{\left(\frac{c}{3}\right)_n\, \left(\frac{c+1}{3}\right)_n\, \left(\frac{c+2}{3}\right)_n\, (1)_n}
\end{eqnarray*}
Using the formula (\ref{ch2inteq6}) and using the fact that $\Gamma(a+1)= a\Gamma(a)$, the aforementioned equation reduces to
\begin{flushleft}
$\displaystyle\sum_{n=0}^{\infty} \frac{(n+1)\,(a)_n\, \left(\frac{b}{3}\right)_n\, \left(\frac{b+1}{3}\right)_n\, \left(\frac{b+2}{3}\right)_n }
{\left(\frac{c}{3}\right)_n\, \left(\frac{c+1}{3}\right)_n\, \left(\frac{c+2}{3}\right)_n\, (1)_n}$
\end{flushleft}
\begin{eqnarray*}
&=&  \frac{\Gamma(c)\, \Gamma(c-a-b)}{\Gamma(c-a)\,\Gamma(c-b)} \,\bigg( \sum_{n=0}^{\infty}\,\bigg(\frac{(a)_{n+1}\,(-1)^n\,(b)_{n+3}}{n!\,(c-a)_{n+2}\,(c-a-b-1)}\bigg)\cr
&&   \qquad \times \, _{2}F_1(-n,b+3+n;c-a+2+n;-1)\,\cr
&&   +\,\sum_{n=0}^{\infty}\,\bigg( \frac{(a)_n\, (-1)^n\,(b)_{n}\,}{n!\,(c-a)_{n}\,}\bigg)\,_{2}F_1(-n,b+n;c-a+n;-1)\bigg).
\end{eqnarray*}
Hence, (1) is proved.%

(2) Using ascending factorial notation and by adjusting coefficients suitably, we can write
\begin{flushleft}
  $\displaystyle \sum_{n=0}^{\infty} \frac{(n+1)^2\, (a)_n\, \left(\frac{b}{3}\right)_n\, \left(\frac{b+1}{3}\right)_n\, \left(\frac{b+2}{3}\right)_n}{\left(\frac{c}{3}\right)_n\, \left(\frac{c+1}{3}\right)_n\, \left(\frac{c+2}{3}\right)_n\, (1)_n}$
\end{flushleft}
\begin{eqnarray*}
&=& \displaystyle\sum_{n=0}^{\infty} \frac{(a)_{n+2}\, \left(\frac{b}{3}\right)_{n+2}\, \left(\frac{b+1}{3}\right)_{n+2}\, \left(\frac{b+2}{3}\right)_{n+2}}{\left(\frac{c}{3}\right)_{n+2}\, \left(\frac{c+1}{3}\right)_{n+2}\, \left(\frac{c+2}{3}\right)_{n+2}\, (1)_{n}} \\
&&+ 3 \displaystyle\sum_{n=0}^{\infty} \frac{(a)_{n+1}\, \left(\frac{b}{3}\right)_{n+1}\, \left(\frac{b+1}{3}\right)_{n+1}\, \left(\frac{b+2}{3}\right)_{n+1}}{\left(\frac{c}{3}\right)_{n+1}\, \left(\frac{c+1}{3}\right)_{n+1}\, \left(\frac{c+2}{3}\right)_{n+1}\, (1)_{n}}+ \displaystyle\sum_{n=0}^{\infty} \frac{(a)_n\, \left(\frac{b}{3}\right)_n\, \left(\frac{b+1}{3}\right)_n\, \left(\frac{b+2}{3}\right)_n}{\left(\frac{c}{3}\right)_n\, \left(\frac{c+1}{3}\right)_n\, \left(\frac{c+2}{3}\right)_n\, (1)_n}
\end{eqnarray*}
Using the formula (\ref{ch2inteq6}) and using the fact that $\Gamma(a+1)= a\Gamma(a)$, the aforementioned equation reduces to
\begin{flushleft}
$\displaystyle\sum_{n=0}^{\infty} \frac{(n+1)^2\,(a)_n\, \left(\frac{b}{3}\right)_n\, \left(\frac{b+1}{3}\right)_n\, \left(\frac{b+2}{3}\right)_n }
{\left(\frac{c}{3}\right)_n\, \left(\frac{c+1}{3}\right)_n\, \left(\frac{c+2}{3}\right)_n\, (1)_n}$
\end{flushleft}
\begin{eqnarray*}
&=&   \frac{\Gamma(c)\, \Gamma(c-a-b)}{\Gamma(c-a)\,\Gamma(c-b)} \,\bigg( \sum_{n=0}^{\infty}\,\bigg(\frac{(a)_{n+2}\,(-1)^n\,(b)_{n+6}}{n!\,(c-a)_{n+4}\,(c-a-b-2)_2}\bigg)\cr
&& \qquad \times \, _{2}F_1(-n,b+6+n;c-a+4+n;-1)\,\cr
&& + 3\,\sum_{n=0}^{\infty}\,\bigg(\frac{(a)_{n+1}\,(-1)^n\,(b)_{n+3}}{n!\,(c-a)_{n+2}\,(c-a-b-1)}\bigg)\cr
&& \qquad \times\, _{2}F_1(-n,b+3+n;c-a+2+n;-1)\,\cr
&& +\,\sum_{n=0}^{\infty}\,\bigg( \frac{(a)_n\, (-1)^n\,(b)_{n}\,}{n!\,(c-a)_{n}\,}\bigg)\,_{2}F_1(-n,b+n;c-a+n;-1)\bigg),
\end{eqnarray*}
Which completes the proof of (2).%

(3)  Using $(n+1)^3 = n(n-1)(n-2)+6(n)(n-1)+7(n)+1$, we can easily obtain that
\begin{flushleft}
$\displaystyle \sum_{n=0}^{\infty} \frac{(n+1)^3 \,(a)_n\, \left(\frac{b}{3}\right)_n\, \left(\frac{b+1}{3}\right)_n\, \left(\frac{b+2}{3}\right)_n }
{\left(\frac{c}{3}\right)_n\, \left(\frac{c+1}{3}\right)_n\, \left(\frac{c+2}{3}\right)_n\, (1)_n}$
\end{flushleft}
\begin{eqnarray*}
&=& \sum_{n=0}^{\infty} \frac{(a)_{n+3}\, \left(\frac{b}{3}\right)_{n+3}\, \left(\frac{b+1}{3}\right)_{n+3}\, \left(\frac{b+2}{3}\right)_{n+3} }{\left(\frac{c}{3}\right)_{n+3}\, \left(\frac{c+1}{3}\right)_{n+3}\, \left(\frac{c+2}{3}\right)_{n+3}\, (1)_{n}}\\
&& + 6\displaystyle \sum_{n=0}^{\infty} \frac{ (a)_{n+2}\, \left(\frac{b}{3}\right)_{n+2}\, \left(\frac{b+1}{3}\right)_{n+2}\, \left(\frac{b+2}{3}\right)_{n+2}}{\left(\frac{c}{3}\right)_{n+2}\, \left(\frac{c+1}{3}\right)_{n+2}\, \left(\frac{c+2}{3}\right)_{n+2}\, (1)_{n}} \cr
&&+7 \displaystyle \sum_{n=0}^{\infty} \frac{(a)_{n+1}\, \left(\frac{b}{3}\right)_{n+1}\, \left(\frac{b+1}{3}\right)_{n+1}\, \left(\frac{b+2}{3}\right)_{n+1}}{\left(\frac{c}{3}\right)_{n+1}\, \left(\frac{c+1}{3}\right)_{n+1}\, \left(\frac{c+2}{3}\right)_{n+1}\, (1)_{n}} + \displaystyle \sum_{n=0}^{\infty} \frac{(a)_n\, \left(\frac{b}{3}\right)_n\, \left(\frac{b+1}{3}\right)_n\, \left(\frac{b+2}{3}\right)_n}{\left(\frac{c}{3}\right)_n\, \left(\frac{c+1}{3}\right)_n\, \left(\frac{c+2}{3}\right)_n\, (1)_n}
\end{eqnarray*}
Using the formula (\ref{ch2inteq6}) and using the fact that $\Gamma(a+1)= a\Gamma(a)$, the aforementioned equation reduces to
\begin{flushleft}
$\displaystyle\sum_{n=0}^{\infty} \frac{(n+1)^3\,(a)_n\, \left(\frac{b}{3}\right)_n\, \left(\frac{b+1}{3}\right)_n\, \left(\frac{b+2}{3}\right)_n }
{\left(\frac{c}{3}\right)_n\, \left(\frac{c+1}{3}\right)_n\, \left(\frac{c+2}{3}\right)_n\, (1)_n}$
\end{flushleft}
\begin{eqnarray*}
&=&   \frac{\Gamma(c)\, \Gamma(c-a-b)}{\Gamma(c-a)\,\Gamma(c-b)} \,\bigg(  \sum_{n=0}^{\infty}\,\bigg(\frac{(a)_{n+3}\,(-1)^n\,(b)_{n+9}}{n!\,(c-a)_{n+6}\,(c-a-b-3)_3}\bigg)\cr
&& \qquad \times \, _{2}F_1(-n,b+9+n;c-a+6+n;-1)\,\cr
&&+ 6 \sum_{n=0}^{\infty}\,\bigg(\frac{(a)_{n+2}\,(-1)^n\,(b)_{n+6}}{n!\,(c-a)_{n+4}\,(c-a-b-2)_2}\bigg)\, _{2}F_1(-n,b+6+n;c-a+4+n;-1)\,\cr
&&+ 7\,\sum_{n=0}^{\infty}\,\bigg(\frac{(a)_{n+1}\,(-1)^n\,(b)_{n+3}}{n!\,(c-a)_{n+2}\,(c-a-b-1)}\bigg)\, _{2}F_1(-n,b+3+n;c-a+2+n;-1)\,\cr
&& +\,\sum_{n=0}^{\infty}\,\bigg( \frac{(a)_n\, (-1)^n\,(b)_{n}\,}{n!\,(c-a)_{n}\,}\bigg)\,_{2}F_1(-n,b+n;c-a+n;-1)\bigg),
\end{eqnarray*}
which completes the proof.%

(4)  Let, $a\neq 1$, $b\neq 1,\, 2,\, 3$, $c > \max \{a+2, a+b-1\}$, it is found that
\begin{flushleft}
$\displaystyle\sum_{n=0}^{\infty} \frac{(a)_n\, \left(\frac{b}{3}\right)_n\, \left(\frac{b+1}{3}\right)_n\, \left(\frac{b+2}{3}\right)_n }
{\left(\frac{c}{3}\right)_n\, \left(\frac{c+1}{3}\right)_n\, \left(\frac{c+2}{3}\right)_n\, (1)_{n+1}} $
\end{flushleft}
\begin{eqnarray*}
&=& \left(\frac{(c-3)\, (c-2)\, (c-1)}{(a-1)\, (b-3)\, (b-2)\, (b-1)}\right)\\
&&\qquad\times \displaystyle\left(\displaystyle\sum_{n=0}^{\infty}\frac{(a-1)_n\, (b-3)_n\, (b-2)_n\, (b-1)_n}{(c-3)_n\, (c-2)_n\, (c-1)_n\, (1)_n}-1 \right) \\
&=& \frac{\Gamma(c)\, \Gamma(c-a-b)}{\Gamma(c-b)\, \Gamma(c-a)} \left(\frac{(c-a-1)\, (c-a-b)}{(a-1)\, (b-3)_3}\right)\cr
&&\qquad\times \left(\displaystyle\sum_{n=0}^{\infty}\frac{(a-1)_n\, (-1)^n\, (b-3)_n}{n!\,(c-a-2)_n}\right)\, _2F_1 (\begin{array}{cccc}
            -n, & b-3+n; & c-a-2+n; & -1
            \end{array})\\
&& - \left(\frac{(c-3)(c-2)(c-1)}{(a-1)(b-1)(b-2)(b-3)}\right).
\end{eqnarray*}
Hence the desired result follows.
\epf
\subsection{$z\, _4F_3\left(^{a,\,\frac{b}{3},\, \frac{b+1}{3},\, \frac{b+2}{3}}_{\frac{c}{3},\, \frac{c+1}{3},\, \frac{c+1}{3}};z\right)$ to belong to the classes $ \es^{*}_{\lambda},\, \CC_{\lambda},\, \UCV,\,$ and $\, \es_p,\,  0 < \lambda  \leq 1$}
\bthm\label{ch2thm1eq0}
 Let, $a,\, b \in {\Bbb C} \backslash \{ 0 \} $,\, $c > 0$\, and $c > |a|+|b|+1.$ A sufficient condition for the function $z\, _4F_3\left(^{a,\,\frac{b}{3},\, \frac{b+1}{3},\, \frac{b+2}{3}}_{\frac{c}{3},\, \frac{c+1}{3},\, \frac{c+2}{3}};z\right) $ to belong to the class $ \es^{*}_{\lambda}, \,  0 < \lambda  \leq 1 $ is that
\beq\label{ch2thm1eq1}
\frac{\Gamma(c)\, \Gamma(c-|a|-|b|)}{\Gamma(c-|a|)\, \Gamma(c-|b|)}\bigg[\sum_{n=0}^{\infty}\left(\frac{(|a|)_{n+1}\,(-1)^n\, (|b|)_{n+3}}{n!\,(c-|a|)_{n+2}\, (c -|a|-|b|-1)}\right)\qquad\qquad\qquad&&\cr
\qquad\times\, _{2}F_1(-n,|b|+3+n;c-|a|+2+n;-1)&& \cr
+\, \lambda\, \left(\,\sum_{n=0}^{\infty}\, \frac{(|a|)_n\, (-1)^n\, (|b|)_n}{n!\,(c-|a|)_n}\right)\, _{2}F_1(-n,|b|+n;c-|a|+n;-1)\bigg] &\leq& 2\lambda.
\eeq
\ethm
\bthm\label{ch2thm10eq1}
Let,  $a, b \in {\Bbb C} \backslash \{ 0 \} ,\, c > 0,\,$ $c > |a|+|b|+2$ and $0 < \lambda \leq 1$. A sufficient condition for the function $z\, _4F_3\left(^{a,\,\frac{b}{3},\, \frac{b+1}{3},\, \frac{b+2}{3}}_{\frac{c}{3},\, \frac{c+1}{3},\, \frac{c+2}{3}};z\right)$ to belong to the class $ \mathcal{C}_{\lambda}$ is that
 \beq\label{ch2thm10eq10}
 \left(\frac{\Gamma(c)\,\Gamma(c-|a|-|b|)}{\Gamma(c-|a|)\, \Gamma(c-|b|)}\right)\,\bigg[\sum_{n=0}^{\infty}\,\left(\frac{\,(|a|)_{n+2}\,(-1)^n\,(|b|)_{n+6}}{n!\,(c-|a|)_{n+4}\, (c-|a|-|b|-2)_{2}}\,\right) \qquad\qquad\quad\,\nonumber&&\cr
 \times\, _{2}F_1(-n,|b|+6+n;c-|a|+4+n;-1)\,&& \cr
  + \,\sum_{n=0}^{\infty}\, \left(\frac{(\lambda+2)\, (|a|)_{n+1}\,(-1)^n\,(|b|)_{n+3}}{n!\,(c-|a|)_{n+2}\, (c-|a|-|b|-1)}\right)\, \qquad\qquad\qquad\qquad\qquad\,\,\nonumber&&\cr
 \times\, _{2}F_1(-n,|b|+3+n;c-|a|+2+n;-1)&&\cr
 +\, \lambda\, \sum_{n=0}^{\infty}\frac{(|a|)_n\,(-1)^n\, (|b|)_{n}}{n!\,(c-|a|)_{n}}\, _{2}F_1(-n,|b|+n;c-|a|+n;-1)\bigg] &\leq& 2\lambda.
   \eeq
\ethm
\bthm\label{ch2thm7eq1}
 Let, $a,\, b \in {\Bbb C} \backslash \{ 0 \} ,\,c > 0$  and $c > |a|+|b|+2 .$  A sufficient condition for the function $z\, _4F_3\left(^{a,\,\frac{b}{3},\, \frac{b+1}{3},\, \frac{b+2}{3}}_{\frac{c}{3},\, \frac{c+1}{3},\, \frac{c+2}{3}};z\right) $ to belong to the class ${\UCV}$ is that
 \beq\label{ch2thm7eq10}
 \frac{\Gamma(c)\, \Gamma(c-|a|-|b|)}{\Gamma(c-|a|)\,\Gamma(c-|b|)}\,\bigg( \sum_{n=0}^{\infty}\,\bigg(\frac{2\, (|a|)_{n+2}\,(-1)^n\,(|b|)_{n+6}}{n!\,(c-|a|)_{n+4}\,(c-|a|-|b|-2)_2}\bigg)\qquad\qquad\qquad\qquad\qquad &&\cr
  \qquad\qquad \times\, _{2}F_1(-n,|b|+6+n;c-|a|+4+n;-1)\,&&\cr
+5\,\sum_{n=0}^{\infty}\,\bigg(\frac{(|a|)_{n+1}\,(-1)^n\,(|b|)_{n+3}}{n!\,(c-|a|)_{n+2}\,(c-|a|-|b|-1)}\bigg)\,_{2}F_1(-n,|b|+3+n;c-|a|+2+n;-1)\,&&\cr
 +\, \sum_{n=0}^{\infty}\,\bigg( \frac{(|a|)_n\, (-1)^n\,(|b|)_{n}\,}{n!\,(c-|a|)_{n}\,}\bigg)\,_{2}F_1(-n,|b|+n;c-|a|+n;-1)\bigg)\leq&2.&
\eeq
\ethm
\bthm\label{ch2thm4eq0}
 Let, $a,\, b \in {\Bbb C} \backslash \{ 0 \} $,\, $c > 0$   and $c > |a|+|b|+1.$  A sufficient condition for the function $z\, _4F_3\left(^{a,\,\frac{b}{3},\, \frac{b+1}{3},\, \frac{b+2}{3}}_{\frac{c}{3},\, \frac{c+1}{3},\, \frac{c+2}{3}};z\right) $ to belong to the class $ \es_p$ is that
 \beq\label{ch2thm4eq1}
\frac{\Gamma(c)\, \Gamma(c-|a|-|b|)}{\Gamma(c-|a|)\,\Gamma(c-|b|)}\,\bigg(2\, \sum_{n=0}^{\infty}\,\bigg(\frac{(|a|)_{n+1}\,(-1)^n\,(|b|)_{n+3}}{n!\,(c-|a|)_{n+2}\,(c-|a|-|b|-1)}\bigg)\quad\qquad\qquad\qquad&&\cr
 \qquad\qquad\qquad \times \, _{2}F_1(-n,|b|+3+n;c-|a|+2+n;-1)\,&&\cr
+\,\,\sum_{n=0}^{\infty}\,\bigg( \frac{(|a|)_n\, (-1)^n\,(|b|)_{n}\,}{n!\,(c-|a|)_{n}\,}\bigg)\,_{2}F_1(-n,|b|+n;c-|a|+n;-1)\bigg) &\leq&2.
\eeq
\ethm
\subsection{Results of $\mathcal{I}^{a,\,\frac{b}{3},\, \frac{b+1}{3},\, \frac{b+2}{3}}_{\frac{c}{3},\, \frac{c+1}{3},\, \frac{c+2}{3}}(f)$ maps $ \mathcal{R}(\beta)$ into classless $ \es^{*}_{\lambda},\, \CC_{\lambda},\, \UCV,\,$ and $\, \es_p,\,  0 < \lambda  \leq 1$}
\bthm\label{ch2thm2eq001}
 Let, $a,\, b \in {\Bbb C} \backslash \{ 0 \} ,\, c > 0,\, \, |a|\neq1,\, |b| \neq1,\,2,\,3$ and $c > \max\{|a|+2, |a|+|b|-1\}$. For  $ 0 < \lambda \leq 1$ and $0 \leq \beta < 1$, assume that
 \beq\label{ch2thm2eq1}
  \left(\frac{\Gamma(c)\,\Gamma(c-|a|-|b|)}{\Gamma(c-|a|)\, \Gamma(c-|b|)}\right)\bigg[\left(\frac{(\lambda-1)\,(c-|a|-1)\,(c-|a|-|b|)}{(|a|-1)\,(|b|-3)_3}\right)\qquad\qquad\qquad\qquad\qquad&& \nonumber\\
  \quad \times\, \sum_{n=0}^{\infty}\, \frac{(|a|-1)_n\, (-1)^n\, (|b|-3)_n}{n!\,(c-|a|-2)_n}\,  _{2}F_1(-n,\, |b|-3+n;\, c-|a|-2+n;-1)&& \cr
  \qquad\qquad\,+\,\,\sum_{n=0}^{\infty}\, \frac{(|a|)_n\, (-1)^n\, (|b|)_n}{n!\,(c-|a|)_n} \, _{2}F_1(-n,|b|+n;c-|a|+n;-1)\bigg]&& \nonumber \\
  \qquad\qquad\qquad\leq \lambda\left(1+\frac{1}{2(1-\beta)}\right)+\left(\frac{(\lambda-1)\,(c-3)_3}{(|a|-1)(|b|-3)_3}\right).&&
\eeq
 Then, the integral  operator $\mathcal{I}^{a,\,\frac{b}{3},\, \frac{b+1}{3},\, \frac{b+2}{3}}_{\frac{c}{3},\, \frac{c+1}{3},\, \frac{c+2}{3}}(f)$ maps $ \mathcal{R}(\beta)$ into $\es^{*}_{\lambda}$.
\ethm
\bthm\label{ch2thm11eq0}  Let, $a,\,b \in {\Bbb C} \backslash \{ 0 \} ,\, c > 0,\, c > |a|+|b|+1$, and $ 0 < \lambda \leq 1$.  For $0 \leq \beta <1 $, assume that
 \beq\label{ch2thm11eq1}
 \left(\frac{\Gamma(c)\,\Gamma(c-|a|-|b|)}{\Gamma(c-|a|)\, \Gamma(c-|b|)}\right)\, \bigg[\sum_{n=0}^{\infty}\, \left(\frac{\, (|a|)_{n+1}\,(-1)^n\,(|b|)_{n+3}}{n!\,(c-|a|)_{n+2}\, (c-|a|-|b|-1)}\right)\qquad\qquad\qquad\qquad\qquad\,\nonumber\cr
\,\times\, _{2}F_1(-n,|b|+3+n;c-|a|+2+n;-1)\qquad\qquad\qquad\qquad\qquad\cr
 + \lambda \left(\sum_{n=0}^{\infty}\frac{(|a|)_n\,(-1)^n\, (|b|)_{n}}{n!\,(c-|a|)_{n}}\right){} _{2}F_1(-n,|b|+n;c-|a|+n;-1)\bigg]\leq \lambda\left( \frac{1}{2(1-\beta)}+1\right).
\eeq
Then, the operator $\mathcal{I}^{a,\,\frac{b}{3},\, \frac{b+1}{3},\, \frac{b+2}{3}}_{\frac{c}{3},\, \frac{c+1}{3},\, \frac{c+2}{3}}(f)$ maps $ \mathcal{R}(\beta)$ into $ \mathcal{C}_{\lambda}$.
\ethm
\bthm\label{ch2thm8eq0}  Let, $a, \, b \in {\Bbb C} \backslash \{ 0 \} ,\, c > 0,\,$ $c > |a|+|b|+1$ and $0 \leq \beta <1 $. Assume that
 \beq\label{ch2thm8eq1}
\displaystyle\frac{\Gamma(c)\,\Gamma(c-|a|-|b|)}{\Gamma(c-|a|)\, \Gamma(c-|b|)}\,\bigg[ 2\, \sum_{n=0}^{\infty}\,\bigg(\frac{(|a|)_{n+1}\,(-1)^n\,(|b|)_{n+3}}{n!\,(c-|a|)_{n+2}\,(c-|a|-|b|-1)}\bigg)\qquad\qquad\qquad\qquad&&\cr
\qquad \qquad \times \, _{2}F_1(-n,|b|+3+n;c-|a|+2+n;-1)\, &&\cr
+\,\sum_{n=0}^{\infty}\bigg( \frac{(|a|)_n\, (-1)^n\,(|b|)_{n}}{n!\,(c-|a|)_{n}}\bigg){}_{2}F_1(-n,|b|+n;c-|a|+n;-1)\bigg] \leq \frac{1}{2(1-\beta)}&+1.&
\eeq
Then, $\mathcal{I}^{a,\,\frac{b}{3},\, \frac{b+1}{3},\, \frac{b+2}{3}}_{\frac{c}{3},\, \frac{c+1}{3},\, \frac{c+2}{3}}(f)$  maps $\mathcal{R}(\beta)$ into $ {\UCV}$.
\ethm
\bthm\label{ch2thm5eq0}  Let, $a,\, b \in {\Bbb C} \backslash \{ 0 \} ,\,|a|\neq1,\, |b|\neq 1,2,3,\,   c > 0 ,\,  c > \max\{|a|+2, |a|+|b|-1\}$ and $0 \leq \beta <1$. Assume that
\beq\label{ch4thm5eq1}
\frac{\Gamma(c)\, \Gamma(c-|a|-|b|)}{\Gamma(c-|a|)\,\Gamma(c-|b|)}\bigg[2\, \sum_{n=0}^{\infty}\,\bigg(\frac{(|a|)_{n}\,(-1)^n\,(|b|)_{n}}{n!\,(c-|a|)_{n}}\bigg)\,\, _{2}F_1(-n,|b|+n;c-|a|+n;-1)\qquad\mathfrak{}&&\cr
-\,\frac{(c-|a|-1)\, (c-|a|-|b|)}{(|a|-1)\,(|b|-3)_3}\,\bigg(\,\sum_{n=0}^{\infty} \frac{(|a|-1)_n\, (-1)^n\,(|b|-3)_{n}\,}{n!\,(c-|a|-2)_{n}}\bigg)\qquad&&\cr
\qquad\times\,_{2}F_1(-n,|b|-3+n;c-|a|-2+n;-1)\,+\frac{(c-3)_3}{(|a|-1)\,(|b|-3)_3}\bigg]\leq \frac{1}{2(1-\beta)}+&1.&
\eeq
Then, $\mathcal{I}^{a,\,\frac{b}{3},\, \frac{b+1}{3},\, \frac{b+2}{3}}_{\frac{c}{3},\, \frac{c+1}{3},\, \frac{c+2}{3}}(f)$ maps $\mathcal{R}(\beta)$ into $ \es_p$ class.
\ethm
\subsection{Results of $\mathcal{I}^{a,\,\frac{b}{3},\, \frac{b+1}{3},\, \frac{b+2}{3}}_{\frac{c}{3},\, \frac{c+1}{3},\, \frac{c+2}{3}}(f)$ maps $ \es$ into classless $ \es^{*}_{\lambda},\, \CC_{\lambda},\, and\, \es_p,\,  0 < \lambda  \leq 1$}
\bthm\label{ch2thm3eq0}   Let, $a, b \in {\Bbb C} \backslash \{ 0 \} ,\, c > 0,\,$ $c > |a|+|b|+2$ and $0 < \lambda \leq 1$. If
 \beq\label{ch2thm3eq1}
\frac{\Gamma(c)\,\Gamma(c-|a|-|b|)}{\Gamma(c-|a|)\, \Gamma(c-|b|)}\bigg[ \left(\sum_{n=0}^{\infty}\frac{(|a|)_{n+2}\,(-1)^n\,(|b|)_{n+6}}{(c-|a|-|b|-2)_2\,n!\,(c-|a|)_{n+4}}\right)\quad\qquad\qquad\qquad\qquad&&\cr
\times\,_{2}F_1(-n,|b|+6+n;c-|a|+4+n;-1)\, \qquad\qquad\qquad\qquad&&\cr  \, + \left(\sum_{n=0}^{\infty}\frac{(\lambda+2)\,(|a|)_{n+1}\,(-1)^n\,(|b|)_{n+3}}{(c-|a|-|b|-1)\,n!\,(c-|a|)_{n+2}}\right) \,_{2}F_1(-n,|b|+3+n;c-|a|+2+n;-1)&& \cr
+\,\lambda\, \left(\sum_{n=0}^{\infty}\frac{(|a|)_{n}\,(-1)^n\,(|b|)_{n}}{n!\,(c-|a|)_{n}}\right)\,\,_{2}F_1(-n,|b|+n;c-|a|+n;-1)\bigg] &\leq& 2\lambda.
\eeq
Then, the integral  operator $\mathcal{I}^{a,\,\frac{b}{3},\, \frac{b+1}{3}\, \frac{b+2}{3}}_{\frac{c}{3},\, \frac{c+1}{3}\, \frac{c+2}{3}}(f)(z)$ maps $\mathcal{S}$ to $\es^{*}_{\lambda}$.
\ethm
\bthm\label{ch2thm12eq0}  Let, $a,\,b \in {\Bbb C} \backslash \{ 0 \},\,  c > 0$, $c > |a|+|b|+3$, and $ 0 < \lambda \leq 1$. If
 \beq\label{ch2thm12eq1}
 \left(\frac{\Gamma(c)\,\Gamma(c-|a|-|b|)}{\Gamma(c-|a|)\, \Gamma(c-|b|)}\right)\qquad\qquad\qquad\qquad\qquad\qquad\qquad\qquad\qquad\qquad\qquad\qquad\qquad\qquad\cr
 \times\bigg[ \left(\sum_{n=0}^{\infty}\frac{(|a|)_{n+3}\, (-1)^n\,(|b|)_{n+9}}{n!\,(c-|a|)_{n+6}\, (c-|a|-|b|-3)_{3}}\,\right)\,_{2}F_1(-n,|b|+9+n;c-|a|+6+n;-1)\,\,\,&&\cr
 \,+\,\sum_{n=0}^{\infty}\,\left(\frac{\,(\lambda+5)\,(|a|)_{n+2}\,(-1)^n\,(|b|)_{n+6}}{n!\,(c-|a|)_{n+4}\, (c-|a|-|b|-2)_{2}}\,\right) \,_{2}F_1(-n,|b|+6+n;c-|a|+4+n;-1)\,\,\,&& \cr
 \,+ \,\sum_{n=0}^{\infty}\, \left(\frac{(3\,\lambda+4)\, (|a|)_{n+1}\,(-1)^n\,(|b|)_{n+3}}{n!\,(c-|a|)_{n+2}\, (c-|a|-|b|-1)}\right)\,\, _{2}F_1(-n,|b|+3+n;c-|a|+2+n;-1)\,\,\,\,&&\cr \qquad\qquad\qquad
 +\, \lambda\,\left( \sum_{n=0}^{\infty}\frac{(|a|)_n\,(-1)^n\, (|b|)_{n}}{n!\,(c-|a|)_{n}}\right)\, _{2}F_1(-n,|b|+n;c-|a|+n;-1)\bigg] \leq 2\lambda,&&
\eeq
then, $\mathcal{I}^{a,\,\frac{b}{3},\, \frac{b+1}{3},\, \frac{b+2}{3}}_{\frac{c}{3},\, \frac{c+1}{3},\, \frac{c+2}{3}}(f)$ maps $\es$ into $ \mathcal{C}_{\lambda} $.
\ethm
\bthm\label{ch2thm6eq0}  Let, $a,\, b \in {\Bbb C} \backslash \{ 0 \} ,\, c > 0,\,$ and $c > |a|+|b|+2$. Suppose $a,\, b$ and $c$ satisfy the condition
 \beq\label{ch2thm6eq1}
\left( \frac{\Gamma(c)\,\Gamma(c-|a|-|b|)}{\Gamma(c-|a|)\, \Gamma(c-|b|)}\right)\, \bigg[ 2\,\sum_{n=0}^{\infty} \left(\frac{ \,(|a|)_{n+2}(-1)^n\,(|b|)_{n+6}}{n!\,(c-|a|)_{n+4}\, (c-|a|-|b|-2)_{2}}\,\right)\qquad\qquad\qquad\,\, &&\cr
 \times \,_{2}F_1(-n,|b|+6+n;c-|a|+4+n;-1)\,\qquad\qquad\qquad &&\cr
  \,+\, 5\, \sum_{n=0}^{\infty}\left(\frac{(|a|)_{n+1}\, (-1)^n\,(|b|)_{n+3}}{n!\,(c-|a|)_{n+2}\, (c-|a|-|b|-1)}\right)\, _{2}F_1(-n,|b|+3+n;c-|a|+2+n;-1)&&\cr
  +\sum_{n=0}^{\infty}\frac{(|a|)_n\, (-1)^n\,(|b|)_n}{n!\, (c-|a|)_n}\, _{2}F_1(-n,|b|+n;c-|a|+n;-1)\bigg] \leq&2.&
\eeq
Then, $\mathcal{I}^{a,\,\frac{b}{3},\, \frac{b+1}{3},\, \frac{b+2}{3}}_{\frac{c}{3},\, \frac{c+1}{3},\, \frac{c+2}{3}}(f)$ maps $\es$ into $ \es_{p} $ class.
\ethm
\subsection{Proofs of Main Theorems.}

\bpf {\it of Theorem \ref{ch2thm1eq0}:} \\
Let, $f(z)=z\, _4F_3\left(^{a,\,\frac{b}{3},\, \frac{b+1}{3},\, \frac{b+2}{3}}_{\frac{c}{3},\, \frac{c+1}{3},\, \frac{c+1}{3}};z\right) $, then, by the equation (\ref{ch2inteq2}), it is enough to show that
\begin{eqnarray*}
  T &=& \sum_{n=2}^{\infty}(n+\lambda-1)|A_n|\leq \lambda.
\end{eqnarray*}
Using the fact $|(a)_n|\leq  (|a|)_n$, one can get
\begin{eqnarray*}
  T &\leq&  \sum_{n=0}^{\infty} ((n+1)+(\lambda-1))\left(\frac{(|a|)_{n}\left(\frac{b}{3}\right)_{n}\, \left(\frac{b+1}{3}\right)_{n}\, \left(\frac{b+2}{3}\right)_{n}}{\left(\frac{c}{3}\right)_{n}\, \left(\frac{c+1}{3}\right)_{n}\, \left(\frac{c+2}{3}\right)_{n-1}\,(1)_{n}}\right)-1-(\lambda-1).
\end{eqnarray*}
Using (\ref{ch2inteq6}) and the result (1) of Lemma \ref{ch2lem1eq1} in the aforesaid equation, we get
\begin{eqnarray*}
T &\leq& \frac{\Gamma(c)\, \Gamma(c-|a|-|b|)}{\Gamma(c-|a|)\, \Gamma(c-|b|)}\bigg(\sum_{n=0}^{\infty}\left(\frac{(|a|)_{n+1}\,(-1)^n\, (|b|)_{n+3}}{n!\,(c-|a|)_{n+2}\, (c -|a|-|b|-1)}\right)\cr
   &&  \qquad \qquad\times\, _{2}F_1(-n,|b|+3+n;c-|a|+2+n;-1)\cr
  && +\, \,\sum_{n=0}^{\infty}\, \frac{(|a|)_n\, (-1)^n\, (|b|)_n}{n!\,(c-|a|)_n}\, _{2}F_1(-n,|b|+n;c-|a|+n;-1)\cr
  && +\, (\lambda-1)\,\sum_{n=0}^{\infty}\, \frac{(|a|)_n\, (-1)^n\, (|b|)_n}{n!\,(c-|a|)_n}\, _{2}F_1(-n,|b|+n;c-|a|+n;-1)\bigg)-\lambda.
\end{eqnarray*}
With regard to (\ref{ch2thm1eq1}), the above expression is bounded above by $\lambda$, and hence,
\begin{eqnarray*}
  T &\leq& \frac{\Gamma(c)\, \Gamma(c-|a|-|b|)}{\Gamma(c-|a|)\, \Gamma(c-|b|)}\bigg(\sum_{n=0}^{\infty}\left(\frac{(|a|)_{n+1}\,(-1)^n\, (|b|)_{n+3}}{n!\,(c-|a|)_{n+2}\, (c -|a|-|b|-1)}\right)\cr
   &&  \qquad \qquad\times\, _{2}F_1(-n,|b|+3+n;c-|a|+2+n;-1)\cr
  && +\,\lambda \,\sum_{n=0}^{\infty}\, \frac{(|a|)_n\, (-1)^n\, (|b|)_n}{n!\,(c-|a|)_n}\, _{2}F_1(-n,|b|+n;c-|a|+n;-1)\bigg)-\lambda \leq \lambda.
\end{eqnarray*}
Therefore, $z\, _4F_3\left(^{a,\,\frac{b}{3},\, \frac{b+1}{3},\, \frac{b+2}{3}}_{\frac{c}{3},\, \frac{c+1}{3},\, \frac{c+1}{3}};z\right)$ belongs to the class $\es^{*}_{\lambda}. $
\epf
\bpf {\it of Theorem \ref{ch2thm10eq1}:}
The proof is similar to Theorem \ref{ch2thm3eq0}. So, we omit the details.
\epf
\bpf {\it of Theorem \ref{ch2thm7eq1}:} Let, $a,\, b \in {\Bbb C} \backslash \{ 0 \} ,\,c > 0$  and $c > |a|+|b|+2 .$\\
Let,  $\displaystyle f(z)=\, _4F_3\left(^{a,\,\frac{b}{3},\, \frac{b+1}{3},\, \frac{b+2}{3}}_{\frac{c}{3},\, \frac{c+1}{3},\, \frac{c+2}{3}};z\right)$. Then, by (\ref{ch2lem4eq1}), it is enough to show that
\begin{eqnarray*}
  T &:=& \sum_{n=2}^{\infty}\, n\, (2n-1)\,|A_n|\leq 1,
\end{eqnarray*}
where, $A_n$ is given by (\ref{ch2inteq007}). Using the fact $|(a)_n|\leq  (|a|)_n$,
\begin{eqnarray*}
  T &\leq&  \sum_{n=2}^{\infty} n\, (2n-1)\,\left(\frac{(|a|)_{n-1}\left(\frac{|b|}{3}\right)_{n-1}\, \left(\frac{|b|+1}{3}\right)_{n-1}\, \left(\frac{|b|+2}{3}\right)_{n-1}}{\left(\frac{c}{3}\right)_{n-1}\, \left(\frac{c+1}{3}\right)_{n-1}\, \left(\frac{c+2}{3}\right)_{n-1}(1)_{n-1}}\right)\\
   &=&  2\sum_{n=0}^{\infty} (n+1)^2\, \left(\frac{(|a|)_{n}\left(\frac{|b|}{3}\right)_{n}\, \left(\frac{|b|+1}{3}\right)_{n}\, \left(\frac{|b|+2}{3}\right)_{n}}{\left(\frac{c}{3}\right)_{n}\, \left(\frac{c+1}{3}\right)_{n}\, \left(\frac{c+2}{3}\right)_{n}(1)_{n}}\right)\\
   && -\sum_{n=0}^{\infty} (n+1)\, \left(\frac{(|a|)_{n}\left(\frac{|b|}{3}\right)_{n}\, \left(\frac{|b|+1}{3}\right)_{n}\, \left(\frac{|b|+2}{3}\right)_{n}}{\left(\frac{c}{3}\right)_{n}\, \left(\frac{c+1}{3}\right)_{n}\, \left(\frac{c+2}{3}\right)_{n}(1)_{n}}\right)-1.
\end{eqnarray*}
Using (1) and (2) of Lemma \ref{ch2lem1eq1} in the aforementioned  equation becomes
\begin{eqnarray*}
T &\leq&  \frac{\Gamma(c)\, \Gamma(c-|a|-|b|)}{\Gamma(c-|a|)\,\Gamma(c-|b|)}\,\,\bigg[2\,\sum_{n=0}^{\infty}\,\bigg(\frac{(|a|)_{n+2}\,(-1)^n\,(|b|)_{n+6}}{n!\,(c-|a|)_{n+4}\,(c-|a|-|b|-2)_2}\bigg)\,\cr
&&\qquad\qquad\times\, _{2}F_1(-n,|b|+6+n;c-|a|+4+n;-1)\,\cr
&& + 5\,\sum_{n=0}^{\infty}\,\bigg(\frac{(|a|)_{n+1}\,(-1)^n\,(|b|)_{n+3}}{n!\,(c-|a|)_{n+2}\,(c-|a|-|b|-1)}\bigg)\\
&&\qquad\qquad\times\, _{2}F_1(-n,|b|+3+n;c-|a|+2+n;-1)\,\cr
&& +\,\sum_{n=0}^{\infty}\,\bigg( \frac{(|a|)_n\, (-1)^n\,(|b|)_{n}\,}{n!\,(c-|a|)_{n}\,}\bigg)\,_{2}F_1(-n,|b|+n;c-|a|+n;-1)\bigg]-1.
\end{eqnarray*}
Because of (\ref{ch2thm7eq10}), the above expression is bounded above by 1, and hence
\begin{eqnarray*}
\frac{\Gamma(c)\, \Gamma(c-|a|-|b|)}{\Gamma(c-|a|)\,\Gamma(c-|b|)}\,\bigg( \sum_{n=0}^{\infty}\,\bigg(2\,\frac{(|a|)_{n+2}\,(-1)^n\,(|b|)_{n+6}}{n!\,(c-|a|)_{n+4}\,(c-|a|-|b|-2)_2}\bigg)\,\,\, \qquad\qquad\qquad\qquad&&\cr
  \qquad\qquad\times\, _{2}F_1(-n,|b|+6+n;c-|a|+4+n;-1)\,\qquad\qquad\qquad\qquad\nonumber&&\cr
+ 5\,\sum_{n=0}^{\infty}\,\bigg(\frac{(|a|)_{n+1}\,(-1)^n\,(|b|)_{n+3}}{n!\,(c-|a|)_{n+2}\,(c-|a|-|b|-1)}\bigg)\,_{2}F_1(-n,|b|+3+n;c-|a|+2+n;-1)\,&&\nonumber\cr
+\,\sum_{n=0}^{\infty}\,\bigg( \frac{(|a|)_n\, (-1)^n\,(|b|)_{n}\,}{n!\,(c-|a|)_{n}\,}\bigg)\,_{2}F_1(-n,|b|+n;c-|a|+n;-1)\bigg)-1\, \, \, \leq&1.&
\end{eqnarray*}
Therefore, $z\, _4F_3\left(^{a,\,\frac{b}{3},\, \frac{b+1}{3},\, \frac{b+2}{3}}_{\frac{c}{3},\, \frac{c+1}{3},\, \frac{c+3}{3}};z\right) $ belongs to the class $ {\UCV}. $
\epf
\bpf {\it of Theorem \ref{ch2thm4eq0}:}
The proof is similar to Theorem \ref{ch2thm8eq0}. So, we omit the details.
\epf
\bpf {\it of Theorem \ref{ch2thm2eq001}:} Let, $a,\, b \in {\Bbb C} \backslash \{ 0 \} ,\, c > 0,\, |a|\neq1,\, |b| \neq1,\,2,\,3$, \,  $c > \max\{|a|+2, |a|+|b|-1\}$,\, $ 0 < \lambda \leq 1$ and $0 \leq \beta < 1$.
Consider the integral operator $\mathcal{I}^{a,\,\frac{b}{3},\, \frac{b+1}{3},\, \frac{b+2}{3}}_{\frac{c}{3},\, \frac{c+1}{3},\, \frac{c+2}{3}}(f)$ defined by (\ref{ch2inteq7}). According to (\ref{ch2inteq2}), it is enough to show that
\begin{eqnarray}\label{ch2thm2eq002}
  T &=& \sum_{n=2}^{\infty}(n+\lambda-1)|A_n|\leq \lambda,
\end{eqnarray}
where, $A_n$ is given by (\ref{ch2inteq007}). Then, it is derived that
\begin{eqnarray*}
  T &\leq&  \sum_{n=2}^{\infty} (n+(\lambda-1)) \left(\frac{(|a|)_{n-1}\left(\frac{|b|}{3}\right)_{n-1}\, \left(\frac{|b|+1}{3}\right)_{n-1}\, \left(\frac{|b|+2}{3}\right)_{n-1}}{\left(\frac{c}{3}\right)_{n-1}\, \left(\frac{c+1}{3}\right)_{n-1}\, \left(\frac{c+2}{3}\right)_{n-1}(1)_{n-1}}\right) |a_n|.
\end{eqnarray*}
Using (\ref{ch2inteq3}) in the aforementioned equation, it found that
\begin{eqnarray*}
  T&\leq&  2(1-\beta)\bigg(\sum_{n=0}^{\infty} \left(\frac{(|a|)_{n}\left(\frac{b}{3}\right)_{n}\, \left(\frac{b+1}{3}\right)_{n}\, \left(\frac{b+2}{3}\right)_{n}}{\left(\frac{c}{3}\right)_{n}\, \left(\frac{c+1}{3}\right)_{n}\, \left(\frac{c+2}{3}\right)_{n}(1)_{n}}\right) -1 \cr
  && \qquad\qquad  + (\lambda-1)\sum_{n=1}^{\infty} \left(\frac{(|a|)_{n}\left(\frac{b}{3}\right)_{n}\, \left(\frac{b+1}{3}\right)_{n}\, \left(\frac{b+2}{3}\right)_{n}}{\left(\frac{c}{3}\right)_{n}\, \left(\frac{c+1}{3}\right)_{n}\, \left(\frac{c+2}{3}\right)_{n}(1)_{n}}\right) \left(\frac{1}{n+1}\right)  \bigg):=T_1.
\end{eqnarray*}
Using the formula (\ref{ch2inteq6}) and the results (1) and (4) of Lemma \ref{ch2lem1eq1}, we find that
\begin{eqnarray*}
  T_1 &\leq& 2(1-\beta)\bigg[ \frac{\Gamma(c)\,\Gamma(c-|a|-|b|)}{\Gamma(c-|a|)\, \Gamma(c-|b|)}\bigg( \left(\frac{(\lambda-1)\,(c-|a|-1)\,(c-|a|-|b|)}{(|a|-1)\, (|b|-3)_3}\right)\cr
  &&\qquad\qquad\times\left(\sum_{n=0}^{\infty}\frac{(|a|-1)_n\,(-1)^n\,(|b|-3)_n}{n!\,(c-|a|-2)_n}\right)\cr
   && \qquad\qquad\qquad\times\,\, _{2}F_1(-n,|b|-3+n;c-|a|-2+n;-1)\cr
  && +\left(\sum_{n=0}^{\infty}\frac{(|a|)_n\,(-1)^n\,(|b|)_n}{n!\,(c-|a|)_n}\right)\,_{2}F_1(-n,|b|+n;c-|a|+n;-1)\bigg)\cr
   && \qquad\qquad\qquad\,- \left(\frac{(\lambda-1)(c-3)_3}{(|a|-1)(|b|-3)_3} \right) -\lambda \bigg].
\end{eqnarray*}
Under the condition (\ref{ch2thm2eq1}), one can get
\begin{eqnarray*}
2(1-\beta)\bigg[ \frac{\Gamma(c)\,\Gamma(c-|a|-|b|)}{\Gamma(c-|a|)\, \Gamma(c-|b|)}\bigg( \left(\frac{(\lambda-1)\,(c-|a|-1)\, (c-|a|-|b|)}{(|a|-1)\, (|b|-3)_3}\right)\qquad\qquad&&\, \, \cr
  \times\left(\sum_{n=0}^{\infty}\frac{(|a|-1)_n\,(-1)^n\,(|b|-3)_n}{n!\,(c-|a|-2)_n}\right)\, _{2}F_1(-n,|b|-3+n;c-|a|-2+n;-1) &&\cr
   +\left(\sum_{n=0}^{\infty}\frac{(|a|)_n\,(-1)^n\,(|b|)_n}{n!\,(c-|a|)_n}\right)\,_{2}F_1(-n,|b|+n;c-|a|+n;-1)\bigg) &&\cr \cr
    \,-\, \left(\frac{(\lambda-1)\,(c-3)_3}{(|a|-1)(|b|-3)_3}\right) -\lambda \bigg] &\leq& \lambda.
\end{eqnarray*}
Thus, it is derived that, the inequalities $T \leq T_1 \leq \lambda $,  and hence (\ref{ch2thm2eq002}) hold. Therefore, it is concluded that the operator $\mathcal{I}^{a,\,\frac{b}{3},\, \frac{b+1}{3},\, \frac{b+2}{3}}_{\frac{c}{3},\, \frac{c+1}{3},\, \frac{c+2}{3}}(f)$ maps $ \mathcal{R}(\beta)$ into $\es^{*}_{\lambda}$, which completes the proof of the theorem.
\epf
\bpf {\it of Theorem \ref{ch2thm11eq0}:}
The proof is similar to Theorem \ref{ch2thm2eq001}. So, we omit the details.
\epf
\bpf {\it of Theorem \ref{ch2thm8eq0}:} Let, $a, \, b \in {\Bbb C} \backslash \{ 0 \} ,\, c > 0,\,$ $c > |a|+|b|+1$ and $0 \leq \beta <1 $.\\
Consider, the integral operator $\mathcal{I}^{a,\,\frac{b}{3},\, \frac{b+1}{3},\, \frac{b+2}{3}}_{\frac{c}{3},\, \frac{c+1}{3},\, \frac{c+2}{3}}(f)$ given in (\ref{ch2inteq7}). According to sufficient condition given in (\ref{ch2lem4eq1}), it is enough to show that
\begin{eqnarray*}
  T &:=& \sum_{n=2}^{\infty}n\, (2n-1)\,|A_n|\leq 1,
\end{eqnarray*}
where, $A_n$ is given by (\ref{ch2inteq007}). Using the fact $|(a)_n|\leq  (|a|)_n$ and $(\ref{ch2inteq3})$ in the aforementioned equation,  it is derived that
\begin{eqnarray*}
  T &\leq&  2(1-\beta)\sum_{n=2}^{\infty} n\,(2n-1)\, \left(\frac{(|a|)_{n-1}\left(\frac{|b|}{3}\right)_{n-1}\, \left(\frac{|b|+1}{3}\right)_{n-1}\, \left(\frac{|b|+2}{3}\right)_{n-1}}{\left(\frac{c}{3}\right)_{n-1}\, \left(\frac{c+1}{3}\right)_{n-1}\, \left(\frac{c+2}{3}\right)_{n-1}(1)_{n-1}\, n}\right)\\
&=&  2(1-\beta)\bigg(\sum_{n=0}^{\infty} \,2(n+1)\, \left(\frac{(|a|)_{n}\left(\frac{|b|}{3}\right)_{n}\, \, \left(\frac{|b|+1}{3}\right)_{n}\,\left(\frac{|b|+2}{3}\right)_{n}}{\left(\frac{c}{3}\right)_{n}\, \left(\frac{c+1}{3}\right)_{n}\, \left(\frac{c+2}{3}\right)_{n}(1)_{n}}\right)\\
   &&\qquad\qquad-\sum_{n=0}^{\infty} \left(\frac{(|a|)_{n}\left(\frac{|b|}{3}\right)_{n}\, \, \left(\frac{|b|+1}{3}\right)_{n}\,\left(\frac{|b|+2}{3}\right)_{n}}{\left(\frac{c}{3}\right)_{n}\, \left(\frac{c+1}{3}\right)_{n}\, \left(\frac{c+2}{3}\right)_{n}(1)_{n}}\right)-1\bigg).
 \end{eqnarray*}
Using the formula (\ref{ch2inteq6}) and (1) of Lemma \ref{ch2lem1eq1}, it is found that\\
\begin{eqnarray*}
T &\leq &  2(1-\beta)\bigg( \frac{\Gamma(c)\,\Gamma(c-|a|-|b|)}{\Gamma(c-|a|)\, \Gamma(c-|b|)}\, \bigg( 2\, \sum_{n=0}^{\infty}\,\bigg(\frac{(|a|)_{n+1}\,(-1)^n\,(|b|)_{n+3}}{n!\,(c-|a|)_{n+2}\,(c-|a|-|b|-1)}\bigg)\cr
&& \qquad \qquad\qquad\times \, _{2}F_1(-n,|b|+3+n;c-|a|+2+n;-1)\,\cr
&&\qquad+\,\sum_{n=0}^{\infty}\,\bigg( \frac{(|a|)_n\, (-1)^n\,(|b|)_{n}\,}{n!\,(c-|a|)_{n}\,}\bigg)\,_{2}F_1(-n,|b|+n;c-|a|+n;-1)\bigg)-1 \bigg).
\end{eqnarray*}
By (\ref{ch2thm8eq1}), the aforementioned expression is bounded above by $1$, and hence
\begin{eqnarray*}
 2(1-\beta)\bigg( \frac{\Gamma(c)\,\Gamma(c-|a|-|b|)}{\Gamma(c-|a|)\, \Gamma(c-|b|)}\, \bigg( 2\, \sum_{n=0}^{\infty}\,\bigg(\frac{(|a|)_{n+1}\,(-1)^n\,(|b|)_{n+3}}{n!\,(c-|a|)_{n+2}\,(c-|a|-|b|-1)}\bigg)\qquad&& \cr
\times \, _{2}F_1(-n,|b|+3+n;c-|a|+2+n;-1)\,\qquad\qquad&& \cr
 +\,\sum_{n=0}^{\infty}\,\bigg( \frac{(|a|)_n\, (-1)^n\,(|b|)_{n}\,}{n!\,(c-|a|)_{n}\,}\bigg)\,_{2}F_1(-n,|b|+n;c-|a|+n;-1)\bigg) -1 \bigg) &\leq& 1.
\end{eqnarray*}
Therefore, the operator $\mathcal{I}^{a,\,\frac{b}{3},\, \frac{b+1}{3},\, \frac{b+2}{3}}_{\frac{c}{3},\, \frac{c+1}{3},\, \frac{c+2}{3}}(f)$ maps $\mathcal{R}(\beta)$ into ${\UCV}$, and the result follows.
\epf
\bpf {\it of Theorem \ref{ch2thm5eq0}:} Let, $a,\, b \in {\Bbb C} \backslash \{ 0 \} ,\,|a|\neq1,\, |b|\neq 1,2,3,\,   c > 0 ,\,  c > \max\{|a|+2, |a|+|b|-1\}$ and $0 \leq \beta <1$. \\ \\
Consider, the integral operator $\mathcal{I}^{a,\,\frac{b}{3},\, \frac{b+1}{3},\, \frac{b+2}{3}}_{\frac{c}{3},\, \frac{c+1}{3},\, \frac{c+2}{3}}(f)$ given by (\ref{ch2inteq7}). In the view of (\ref{ch2lem2eq1}), it is enough to show that
\begin{eqnarray*}
  T &:=& \sum_{n=2}^{\infty}(2n-1)|A_n|\leq 1,
\end{eqnarray*}
where, $A_n$ is given by (\ref{ch2inteq007}). Using the inequality $|(a)_n|\leq  (|a|)_n$, it is proven that
\begin{eqnarray*}
  T &\leq & \sum_{n=2}^{\infty}(2n-1)\bigg(\frac{(|a|)_{n-1}\left(\frac{|b|}{3}\right)_{n-1}\, \left(\frac{|b|+1}{3}\right)_{n-1}\, \left(\frac{|b|+2}{3}\right)_{n-1}}{\left(\frac{c}{3}\right)_{n-1}\, \left(\frac{c+1}{3}\right)_{n-1}\, \left(\frac{c+2}{3}\right)_{n-1}(1)_{n-1}}\bigg)\, |a_n|.
\end{eqnarray*}
Using $(\ref{ch2inteq3})$ in the aforesaid equation, one can have
\begin{eqnarray*}
 &= &  2(1-\beta)\bigg(2\sum_{n=0}^{\infty} \frac{(|a|)_{n}\left(\frac{|b|}{3}\right)_{n}\, \left(\frac{|b|+1}{3}\right)_{n}\, \left(\frac{|b|+2}{3}\right)_{n}}{\left(\frac{c}{3}\right)_{n}\, \left(\frac{c+1}{3}\right)_{n}\, \left(\frac{c+2}{3}\right)_{n}(1)_{n}}\cr
&&\qquad\qquad-\sum_{n=0}^{\infty}  \frac{(|a|)_{n}\left(\frac{|b|}{3}\right)_{n}\, \left(\frac{|b|+1}{3}\right)_{n}\, \left(\frac{|b|+2}{3}\right)_{n}}{\left(\frac{c}{3}\right)_{n}\, \left(\frac{c+1}{3}\right)_{n}\, \left(\frac{c+2}{3}\right)_{n}(1)_{n+1}} -1 \bigg).
\end{eqnarray*}
Using the equation (\ref{ch2inteq6}) and the result (4) of Lemma \ref{ch2lem1eq1}, it is found that
\begin{eqnarray*}
  T &\leq&  2(1-\beta)\bigg(\frac{\Gamma(c)\, \Gamma(c-|a|-|b|)}{\Gamma(c-|a|)\,\Gamma(c-|b|)} \,\bigg(2\, \sum_{n=0}^{\infty}\,\bigg(\frac{(|a|)_{n}\,(-1)^n\,(|b|)_{n}}{n!\,(c-|a|)_{n}}\bigg)\cr
&&\qquad\qquad\times\, _{2}F_1(-n,|b|+n;c-|a|+n;-1)\,\cr
&&-\,\frac{(c-|a|-1)\, (c-|a|-|b|)}{\,(|a|-1)\,(|b|-3)_3}\,\sum_{n=0}^{\infty}\,\bigg( \frac{(|a|-1)_n\, (-1)^n\,(|b|-3)_{n}\,}{n!(c-|a|-2)_{n}}\bigg)\qquad\cr
&&\qquad\times\,_{2}F_1(-n,|b|-3+n;c-|a|-2+n;-1) + \frac{(c-3)_3}{(|a|-1)\,(|b|-3)_3}\bigg)-1  \bigg).
\end{eqnarray*}
By the condition (\ref{ch4thm5eq1}), the aforementioned expression is bounded above by 1, and hence
\begin{eqnarray*}
  2(1-\beta)\bigg(\frac{\Gamma(c)\, \Gamma(c-|a|-|b|)}{\Gamma(c-|a|)\,\Gamma(c-|b|)} \,\bigg( \sum_{n=0}^{\infty}\,\bigg(\frac{(2\,|a|)_{n}\,(-1)^n\,(|b|)_{n}}{n!\,(c-|a|)_{n}}\bigg)\, _{2}F_1(-n,|b|+n;c-|a|+n;-1)&& \cr
-\,\frac{(c-|a|-1)\,(c-|a|-|b|)}{\,(|a|-1)\,(|b|-3)_3}\,\sum_{n=0}^{\infty}\,\bigg( \frac{(|a|-1)_n\, (-1)^n\,(|b|-3)_{n}\,}{n!(c-|a|-2)_{n}}\bigg)&&\cr
\times\,_{2}F_1(-n,|b|-3+n;c-|a|-2+n;-1) + \frac{(c-3)_3}{(|a|-1)\,(|b|-3)_3}\bigg) -1  \bigg) \leq 1.&&
\end{eqnarray*}
Under the stated condition, the operator $\mathcal{I}^{a,\,\frac{b}{3},\, \frac{b+1}{3},\, \frac{b+2}{3}}_{\frac{c}{3},\, \frac{c+1}{3},\, \frac{c+2}{3}}(f)$ maps $\mathcal{R}(\beta)$ into $\es_p$ and  the proof is completed.
\epf
\bpf {\it of Theorem \ref{ch2thm3eq0}:}  Let,  $a, b \in {\Bbb C} \backslash \{ 0 \} ,\, c > 0,\,$ $c > |a|+|b|+2$ and $0 < \lambda \leq 1$. \\
If the integral  operator $\mathcal{I}^{a,\,\frac{b}{3},\, \frac{b+1}{3},\, \frac{b+2}{3}}_{\frac{c}{3},\, \frac{c+1}{3},\, \frac{c+2}{3}}(f)$ is defined by (\ref{ch2inteq7}), in view of (\ref{ch2inteq2}), it is enough to show that
\begin{eqnarray}
  T &:=& \sum_{n=2}^{\infty}(n+\lambda-1)|A_n|\leq \lambda,
\end{eqnarray}
where, $A_n$ is given by (\ref{ch2inteq007}).
Using the fact $|(a)_n|\leq  (|a|)_n$ and the equation ($\ref{ch2inteq00}$) in the aforementioned equation, it is derived that
\begin{eqnarray*}
  T &\leq & \sum_{n=2}^{\infty} (n+(\lambda-1)) \left(\frac{(|a|)_{n-1}\left(\frac{|b|}{3}\right)_{n-1}\, \left(\frac{|b|+1}{3}\right)_{n-1}\, \left(\frac{|b|+2}{3}\right)_{n-1}}{\left(\frac{c}{3}\right)_{n-1}\, \left(\frac{c+1}{3}\right)_{n-1}\, \left(\frac{c+3}{3}\right)_{n-1}(1)_{n-1}}\right)|a_n|\\
  & = & \sum_{n=0}^{\infty} \left(\frac{(n+1)^2 \,(|a|)_{n}\left(\frac{|b|}{3}\right)_{n}\, \left(\frac{|b|+1}{3}\right)_{n}\, \left(\frac{|b|+2}{3}\right)_{n}}{\left(\frac{c}{3}\right)_{n}\, \left(\frac{c+1}{3}\right)_{n}\, \left(\frac{c+3}{3}\right)_{n}(1)_{n}}\right)-1\\
  &&\qquad+(\lambda-1) \,\sum_{n=0}^{\infty}\left(\frac{ (n+1)\,(|a|)_{n}\left(\frac{|b|}{3}\right)_{n}\, \left(\frac{|b|+1}{3}\right)_{n}\, \left(\frac{|b|+2}{3}\right)_{n}}{\left(\frac{c}{3}\right)_{n}\, \left(\frac{c+1}{3}\right)_{n}\, \left(\frac{c+3}{3}\right)_{n}(1)_{n}}\right)-(\lambda-1).
\end{eqnarray*}
Using (1) and (2) of Lemma \ref{ch2lem1eq1}, it is found that
\begin{eqnarray*}
T &\leq&\frac{\Gamma(c)\,\Gamma(c-|a|-|b|)}{\Gamma(c-|a|)\, \Gamma(c-|b|)}\bigg[ \left(\sum_{n=0}^{\infty}\frac{(|a|)_{n+2}\,(-1)^n\,(|b|)_{n+6}}{(c-|a|-|b|-2)_2\,n!\,(c-|a|)_{n+4}}\right) \cr
&&\qquad\qquad\times\,_{2}F_1(-n,|b|+6+n;c-|a|+4+n;-1)\,\cr &&+\,\left(\sum_{n=0}^{\infty}\frac{(\lambda+2)\,(|a|)_{n+1}\,(-1)^n\,(|b|)_{n+3}}{(c-|a|-|b|-1)\,n!\,(c-|a|)_{n+2}}\right) \,\cr
&&\qquad\qquad\times\,_{2}F_1(-n,|b|+3+n;c-|a|+2+n;-1)\, \cr
&&+\,\lambda\, \left(\sum_{n=0}^{\infty}\frac{(|a|)_{n}\,(-1)^n\,(|b|)_{n}}{n!\,(c-|a|)_{n}}\right)\,_{2}F_1(-n,|b|+n;c-|a|+n;-1)\bigg]-\lambda.
\end{eqnarray*}
By (\ref{ch2thm3eq1}), the above expression is bounded above by $\lambda$, and hence
\begin{eqnarray*}
&&\frac{\Gamma(c)\,\Gamma(c-|a|-|b|)}{\Gamma(c-|a|)\, \Gamma(c-|b|)}\times\bigg[ \left(\sum_{n=0}^{\infty}\frac{(|a|)_{n+2}\,(-1)^n\,(|b|)_{n+6}}{(c-|a|-|b|-2)_2\,n!\,(c-|a|)_{n+4}}\right)\,\cr
&&\qquad\qquad\qquad\qquad\times\,_{2}F_1(-n,|b|+6+n;c-|a|+4+n;-1)\,\cr &&\qquad+\,\left(\sum_{n=0}^{\infty}\frac{(\lambda+2)\,(|a|)_{n+1}\,(-1)^n\,(|b|)_{n+3}}{(c-|a|-|b|-1)\,n!\,(c-|a|)_{n+2}}\right)\cr
&& \qquad\qquad\qquad\qquad\times\, \,_{2}F_1(-n,|b|+3+n;c-|a|+2+n;-1)\, \cr
&& \qquad+\,\lambda\, \left(\sum_{n=0}^{\infty}\frac{(|a|)_{n}\,(-1)^n\,(|b|)_{n}}{n!\,(c-|a|)_{n}}\right)\,_{2}F_1(-n,|b|+n;c-|a|+n;-1)\bigg]-\lambda \leq \lambda.
\end{eqnarray*}
Under the stated condition, the integral operator $\mathcal{I}^{a,\,\frac{b}{3},\, \frac{b+1}{3},\, \frac{b+2}{3}}_{\frac{c}{3},\, \frac{c+1}{3},\, \frac{c+2}{3}}(f)$ maps $\es$ into $\es^{*}_{\lambda}$.
\epf
\bpf {\it of Theorem \ref{ch2thm12eq0}:}
Let, $a,\,b \in {\Bbb C} \backslash \{ 0 \},\,  c > 0$, $c > |a|+|b|+3$, and $ 0 < \lambda \leq 1$.\\
Suppose, the integral operator $\mathcal{I}^{a,\,\frac{b}{3},\, \frac{b+1}{3},\, \frac{b+2}{3}}_{\frac{c}{3},\, \frac{c+1}{3},\, \frac{c+2}{3}}(f)$ is defined by (\ref{ch2inteq7}). In view of the sufficient condition given in (\ref{ch2inteq}), it is enough to prove that
\begin{eqnarray*}
  T &:=& \sum_{n=2}^{\infty}\,n\, (n+\lambda-1)\, |A_n|\leq \lambda.
\end{eqnarray*}
Using the fact that $|(a)_n|\leq  (|a|)_n$ and $(\ref{ch2inteq00})$ in the aforementioned equation, it is derived
\begin{eqnarray*}
T &\leq& \sum_{n=2}^{\infty}\,n^2\, (n+\lambda-1)\, \left(\frac{(|a|)_{n-1}\left(\frac{|b|}{3}\right)_{n-1}\, \left(\frac{|b|+1}{3}\right)_{n-1}\, \left(\frac{|b|+2}{3}\right)_{n-1}}{\left(\frac{c}{3}\right)_{n-1}\, \left(\frac{c+1}{3}\right)_{n-1}\, \left(\frac{c+2}{3}\right)_{n-1}\,(1)_{n-1}}\right)\\
 &=&  \sum_{n=0}^{\infty} \left(\frac{((n+1)^3\,|a|)_{n}\left(\frac{|b|}{3}\right)_{n}\, \left(\frac{|b|+1}{3}\right)_{n}\, \left(\frac{|b|+2}{3}\right)_{n}}{\left(\frac{c}{3}\right)_{n}\, \left(\frac{c+1}{3}\right)_{n}\, \left(\frac{c+2}{3}\right)_{n}\,(1)_{n}}\right)\\
  &&+(\lambda-1)\sum_{n=0}^{\infty} \left((n+1)^2\,\frac{(|a|)_{n}\left(\frac{|b|}{3}\right)_{n}\, \left(\frac{|b|+1}{3}\right)_{n}\, \left(\frac{|b|+2}{3}\right)_{n}}{\left(\frac{c}{3}\right)_{n}\, \left(\frac{c+1}{3}\right)_{n}\, \left(\frac{c+2}{3}\right)_{n}\,(1)_{n}}\right)-\lambda.
\end{eqnarray*}
Using (2) and (3) of Lemma \ref{ch2lem1eq1}, it follows that
\begin{eqnarray*}
 &=& \frac{\Gamma(c)\,\Gamma(c-|a|-|b|)}{\Gamma(c-|a|)\, \Gamma(c-|b|)}\,\bigg[ \left(\sum_{n=0}^{\infty}\frac{(|a|)_{n+3}\, (-1)^n\,(|b|)_{n+9}}{n!\,(c-|a|)_{n+6}\, (c-|a|-|b|-3)_{3}}\,\right)  \\
 && \qquad\qquad \times \,_{2}F_1(-n,|b|+9+n;c-|a|+6+n;-1)\,\cr
 &&\,+\,\sum_{n=0}^{\infty}\,\left(\frac{(\lambda+5)\,(|a|)_{n+2}\,(-1)^n\,(|b|)_{n+6}}{n!\,(c-|a|)_{n+4}\, (c-|a|-|b|-2)_{2}}\,\right)  \\
 && \qquad\qquad \times \,_{2}F_1(-n,|b|+6+n;c-|a|+4+n;-1)\,\cr
 &&\, \,+ \,\sum_{n=0}^{\infty}\, \left(\frac{(3\lambda+4)\, (|a|)_{n+1}\,(-1)^n\,(|b|)_{n+3}}{n!\,(c-|a|)_{n+2}\, (c-|a|-|b|-1)}\right) \\
 && \qquad\qquad \times\, _{2}F_1(-n,|b|+3+n;c-|a|+2+n;-1)\,\cr
 &&\, +\, \lambda\, \sum_{n=0}^{\infty}\frac{(|a|)_n\,(-1)^n\, (|b|)_{n}}{n!\,(c-|a|)_{n}}\, _{2}F_1(-n,|b|+n;c-|a|+n;-1)\bigg] -\lambda.
\end{eqnarray*}
By equation (\ref{ch2thm12eq1}), the above expression is bounded above by $\lambda$, and hence
\begin{eqnarray*}
\frac{\Gamma(c)\,\Gamma(c-|a|-|b|)}{\Gamma(c-|a|)\, \Gamma(c-|b|)}\,\bigg[ \left(\sum_{n=0}^{\infty}\frac{(|a|)_{n+3}\, (-1)^n\,(|b|)_{n+9}}{n!\,(c-|a|)_{n+6}\, (c-|a|-|b|-3)_{3}}\,\right) \quad&& \cr
  \qquad\qquad \times \,_{2}F_1(-n,|b|+9+n;c-|a|+6+n;-1)\,
 &\cr  \,+\,\sum_{n=0}^{\infty}\,\left(\frac{\,(\lambda+5)\,(|a|)_{n+2}\,(-1)^n\,(|b|)_{n+6}}{n!\,(c-|a|)_{n+4}\, (c-|a|-|b|-2)_{2}}\,\right)\qquad\qquad\qquad\qquad&& \cr
  \qquad\qquad \times\,_{2}F_1(-n,|b|+6+n;c-|a|+4+n;-1)\,
 &&\cr  \,+ \,\sum_{n=0}^{\infty}\, \left(\frac{(3\lambda+4)\, (|a|)_{n+1}\,(-1)^n\,(|b|)_{n+3}}{n!\,(c-|a|)_{n+2}\, (c-|a|-|b|-1)}\right) \qquad\qquad\qquad\qquad&&\cr
  \qquad\qquad \times\, _{2}F_1(-n,|b|+3+n;c-|a|+2+n;-1)
 && \cr  \qquad +\, \lambda\, \sum_{n=0}^{\infty}\frac{(|a|)_n\,(-1)^n\, (|b|)_{n}}{n!\,(c-|a|)_{n}}\, _{2}F_1(-n,|b|+n;c-|a|+n;-1)\bigg] -\lambda &\leq& \lambda.
\end{eqnarray*}
Hence, the integral operator $\mathcal{I}^{a,\,\frac{b}{3},\, \frac{b+1}{3},\, \frac{b+2}{3}}_{\frac{c}{3},\, \frac{c+1}{3},\, \frac{c+2}{3}}(f)$ maps $\es$ into $\mathcal{C}_{\lambda}$ and the proof is completed.
\epf
\bpf  {\it of Theorem \ref{ch2thm6eq0}:}
The proof is similar to Theorem \ref{ch2thm7eq1}. So, we omit the details.
\epf
\end{sloppypar}
\newpage
\begin{sloppypar}
\section[Convolutions with Generalized Hypergeometric Functions]{Convolutions with Generalized Hypergeometric Functions}
\begin{abstract}
In this article, we consider an operator $\mathcal{I}^{  a,\frac{b}{4},\frac{b+1}{4},\frac{b+2}{4},\frac{b+3}{4} }_{ \frac{c}{4}, \frac{c+1}{4}, \frac{c+2}{4},\frac{c+3}{4} }(f)(z)$$= z\, _5F_4\left(^{a,\frac{b}{4},\frac{b+1}{4},\frac{b+2}{4},\frac{b+3}{4}}_{\frac{c}{4}, \frac{c+1}{4}, \frac{c+2}{4},\frac{c+3}{4}}; z\right)*f(z)$, where, $_5F_4(z)$ hypergeometric function and the $*$ is usual Hadamard  product. In the main results, conditions are determined on $ a,b,$ and $c$ such that the function $z\, _5F_4\left(^{a,\frac{b}{4},\frac{b+1}{4},\frac{b+2}{4},\frac{b+3}{4}}_{\frac{c}{4}, \frac{c+1}{4}, \frac{c+2}{4},\frac{c+3}{4}}; z\right)$ is in the each of the classes  $ \es^{*}_{\lambda} $, $ \mathcal{C}_{\lambda}$,  ${\UCV} $ and $\es_p$. Subsequently, conditions on $a,\,b,\,c,\, \lambda,$ and $\beta$ are determined using the integral operator such that functions belonging to $\mathcal{R}(\beta)$ and $\es$ are mapped onto each of the classes $ \es^{*}_{\lambda}  $, $ \mathcal{C}_{\lambda}$,  ${\UCV},$ and $\es_p$.
\end{abstract}
\maketitle

\subsection{Introduction and preliminaries}

In geometric function theory, the integral operator plays an important role in  characterizing, various subclasses of the univalent functions. The important integral operators, viz. Alexander, Libera, and Bernadi are particular case of the Hohlov convolution operator. The Gaussian Hypergeometric function $_2F_1(a,b;c;z)$ have been studied extensively by a large number of researchers in the context of deriving geometric properties such as convexity, starlikeness, close-to-convexity, and univalency. Whereas, the generalized hypergeometric function $_3F_2(a,b,c;d,e;z)$ have been studied only by a few researchers.\\
In \cite{Chandru-prabha-2019}, Chandrasekran and Prabhakaran introduced an integral operator involving the Clausen’s hypergeometric function and with the help of operator, they derived geometric properties of various subclasses of univalent function. Subsequently, the conditions on the parameters $ a,b,$ and $c$ are determined using the integral operator $\mathcal{I}^{a,\frac{b}{2},\frac{b+1}{2}}_{\frac{c}{2}, \frac{c+1}{2}}(f)(z)$ for various subclasses of univalent functions by the authors Chandrasekran and Prabhakaran \cite{Chandru-prabha-2020}.  Recently, using convolution for $_4F_3(z)$ hypergeometric function, an operator $\mathcal{I}^{a_1,a_2,a_3,a_4}_{b_1,b_2,b_3}(f)(z)$ is introduced in \cite{Chandru-prabha-2021}. Univalence, starlikenes, and convexity properties of ${_4F_3}\left(^{  a,\frac{b}{3},\frac{b+1}{3},\frac{b+2}{3}}_{ \frac{c}{3}, \frac{c+1}{3}, \frac{c+2}{3}};z\right)$  hypergeometric functions are discussed based on its Taylor's coefficient of various subclasses of univalent functions. In sequel we are interested in operators that involves more generalized hypergeometric function with more parameters.\\
Now, we will state the basic concepts of the geometric function theory, which will be helpful to prove our main results.\\
Let, $\A$ denote the family of functions $f(z)= z + a_2\,z^2 + a_3\,z^3 + \cdots$ regular in the open unit disc $\D =\{z:\, |z|<1\}$ of the complex plane and normalized by the conditions $f(0)=f^{\prime}(0)-1=0$. The subfamilies consisting of functions in $\A$ that are univalent, starlike of order $\alpha$ and convex of order $\alpha$, where, $0\leq \alpha <1$ denoted by $\es$, \, $\es^{*}(\alpha)$ and $\CC(\alpha)$, respectively. For $\alpha = 0$, the class $\es^{*}(0)=\es^{*}$ represents the class of starlike and univalent functions and $\CC(0)=\CC$, is the class of convex and univalent functions in the unit disc. For further details \cite{Peter-L-Duren-book-1983,A-W-Goodman-1983-book}) can be referred. If $\displaystyle f \in \A$ belongs to $\es$, then,
\beq\label{inteq00}
|a_n|\leq n,\text{ for}\, n\geq 2.
\eeq%
Our main focus is on the following subclasses of starlike and convex functions that are defined as follows. For $\lambda >0,\, \es^{*}_{\lambda},\,  $ is  defined as follows $$\es^{*}_{\lambda}\, =\, \left\{f(z)\in \A \, |\, \left|\frac{zf'(z)}{f(z)}-1\right| \, < \, \lambda,\, z\in \D \right\}$$ and a sufficient condition for which the function $f(z)$ to be in $\es^{*}_{\lambda}$ is
\beq\label{inteq2}
\displaystyle \sum_{n=2}^{\infty}(n+\lambda-1)|a_n| \leq \lambda.
\eeq%

For $\lambda >0,$ $\CC_{\lambda}$ is defined as follows $\CC_{\lambda}=\left\{f(z)\in \A\, |\, zf'(z)\in \es^{*}_{\lambda}\right\}$ and a sufficient condition for which the function $f(z)$ to be in $\CC_{\lambda}$, is as follows:
\beq\label{inteq}
\displaystyle \sum_{n=2}^{\infty}\, n\, (n+\lambda-1)|a_n| \leq \lambda.
\eeq%

Goodman \cite{Good-1991-Ann-PM} gave a two-variable analytic characterization for the class of uniformly convex, denoted by $\UCV$, that is,
$\UCV = \left\{ f\in \es: Re\left(1+\frac{(z-\zeta)f^{\prime\prime}(z)}{f^{\prime}(z)}\right) > 0, \zeta, z \in \D\right\}$
while R{\o}nning \cite{Ronn-1993-Proc-ams} and Ma and Minda \cite{Ma-Minda-1992-Ann-PM} independently gave the one variable analytic characterization of the class $\UCV$ using the minimum principle for harmonic functions: $$f\in \UCV\, \Leftrightarrow \, \displaystyle  \left|\frac{zf''(z)}{f'(z)}\right| < Re\left(1+\frac{zf^{''}(z)}{f'(z)}\right), \, z\in \D.$$

In \cite{Subram-Murugu-1995}, it was shown that if
\beq\label{lem4eq1}
\sum_{n=2}^{\infty}\,n\, (2n-1)|a_n|\leq 1,
\eeq then, the function $f \in \A$ is in $\UCV$.\\%\\

The  subclass $\es_p$ of starlike functions introduced by R{\o}nning \cite{Ronn-1993-Proc-ams} is defined as
$$f(z) \in \es, \text{is in}\, \es_p \, \Leftrightarrow \, \displaystyle  \left|\frac{zf'(z)}{f(z)}-1\right| < Re\left(\frac{zf^{'}(z)}{f(z)}\right), \, z\in \D.$$
Equivalently, $\displaystyle \es_p = \{ F \in \es^{*}|F(z)=zf'(z),\, f(z) \in \UCV\}.$  A sufficient condition for a function $f\in \A$ belong to $\es_p$ is given by
 \beq \label{lem2eq1}
\sum_{n=2}^{\infty}(2n-1)|a_n|\leq 1.
\eeq%

The above result was proved for more general case $\es_p(\alpha)$ in \cite{Subra-Sudharsan-1998}. For $\beta <1$, let, $${\mathcal R}(\beta)=\{f\in {\mathcal A}:  \exists \  \varphi \in \IR \, |\, {\rm Re}\, \ [ e^{i\varphi}(f'(z)-\beta)] > 0, \quad z\in\D \}.$$

It is evident that ${\mathcal R}(\beta) \subset {\mathcal S}$ when $\beta \ge 0$ and for each $\beta <0,\ \  {\mathcal R}(\beta) $ also contains non-univalent functions. This class has been widely used to study certain integral transforms. For instance, consider \cite{Anbu-Parva-2000,Parva-prabha-2001-Far-East} and the reference presented. Suppose that $\displaystyle f \in \A$ is in the class $ {\mathcal R}(\beta)$, then, applying \cite{MacGregor-1962-Trans-ams}, we get
\beq\label{inteq3}
|a_n|\leq \frac{2(1-\beta)}{n},\, n \geq 2.
\eeq%

For any non-zero complex variable $a$, the Pochhammer symbol (or shifted factorial) is defined as $(a)_0\,=\,1,\, {\rm and }\, (a)_n\,=\,a(a+1)\cdots (a+n-1),\, {\rm for}\, n\,=\,1,2,3,\cdots.$ The Hadamard product (i.e., convolution) of two power series $ \displaystyle f(z)= z+\sum_{n=2}^{\infty}\, a_n\,z^n $ and  $ \displaystyle g(z)= z+\sum_{n=2}^{\infty}\, b_n\,z^n $ be analytic in $\D$ is defined as the power series $$(f*g)(z)= z+\sum_{n=2}^{\infty} a_nb_n z^n.$$

The hypergeometric function $_5F_4(z)$ is defined as
\beq\label{inteq5}
_5F_4\left(^{a_1,\, a_2,\, a_3,\, a_4,\, a_5}_{b_1,\, b_2,\, b_3,\, b_4} ;z\right)=\sum_{n=0}^{\infty}\frac{(a_1)_n(a_2)_n(a_3)_n(a_4)_n(a_5)_n}{(b_1)_n(b_2)_n(b_3)_n(b_4)_n(1)_n}z^n;\, |z|<1
\eeq with $a_1,a_2,a_3,a_4,a_5,b_1,b_2,b_3,b_4 \in \IC$ provided $b_1,\, b_2,\, b_3,\, b_4 \neq 0,-1,-2,-3\cdots,$ which is an analytic functions in unit disc $\D$.\\
In this work, the authors introduced an integral operator $\mathcal{I}^{a_1,a_2,a_3,a_4,a_5}_{b_1,b_2,b_3,b_4}(f)(z)$ involving $_5F_4(z)$ Hypergeometric Function by means of convolution (i.e., Hadamard  product).\\
For $f\in \mathcal{A}$, we define the operator $\mathcal{I}^{a_1,a_2,a_3,a_4,a_5}_{b_1,b_2,b_3,b_4}(f)(z)$
\beq\label{inteq7}
\mathcal{I}^{a_1,a_2,a_3,a_4,a_5}_{b_1,b_2,b_3,b_4}(f)(z) &=& z\, _5F_4\left(^{a_1,\, a_2,\, a_3,\, a_4,\, a_5}_{b_1,\, b_2,\, b_3,\, b_4} ;z\right)*f(z)\\
&=& z+\sum_{n=2}^{\infty} A_n\, z^n,\nonumber
\eeq
with $A_1=1$ and for $n > 1,$
\beq\label{inteq007}
A_n&=&\frac{(a_1)_{n-1}\left(a_2\right)_{n-1}\left(a_3\right)_{n-1}\left(a_4\right)_{n-1}\left(a_5\right)_{n-1}}{\left(b_1\right)_{n-1}
\left(b_2\right)_{n-1}\left(b_3\right)_{n-1}\left(b_4\right)_{n-1}(1)_{n-1}}\, a_n.
\eeq%

In this article, we focus on the following particular values of the parameter for our investigation,
\beq\label{inteq7}
\mathcal{I}^{  a,\frac{b}{4},\frac{b+1}{4},\frac{b+2}{4},\frac{b+3}{4} }_{ \frac{c}{4}, \frac{c+1}{4}, \frac{c+2}{4},\frac{c+3}{4} }(f)(z)&=& z\,
_5F_4\left(^{a,\frac{b}{4},\frac{b+1}{4},\frac{b+2}{4},\frac{b+3}{4}}_{\frac{c}{4}, \frac{c+1}{4}, \frac{c+2}{4},\frac{c+3}{4}}; z\right)*f(z)\\
&=& z+\sum_{n=2}^{\infty} A_n\, z^n,\nonumber
\eeq
with $A_1=1$ and for $n > 1,$
\beq\label{inteq007}
A_n&=&\frac{(a)_{n-1}\left(\frac{b}{4}\right)_{n-1}\left(\frac{b+1}{4}\right)_{n-1}\left(\frac{b+2}{4}\right)_{n-1}\left(\frac{b+3}{4}\right)_{n-1}}
{\left(\frac{c}{4}\right)_{n-1}\left(\frac{c+1}{4}\right)_{n-1}\left(\frac{c+2}{4}\right)_{n-1} \left(\frac{c+3}{4}\right)_{n-1} (1)_{n-1}}\, a_n.
\eeq%

In 2009, Coffey and Johnston \cite{Coffey-johnston-2009} derived a summation formula for $_5F_4\left(1\right)$ hypergeometric function in terms of Gaussian hypergeometric function. We recall their summation formula as follows:
\begin{sloppypar}
\begin{flushleft}
$\displaystyle $
\end{flushleft}
\beq\label{inteq6}
 _5F_4\left(^{a,\frac{b}{4},\frac{b+1}{4},\frac{b+2}{4},\frac{b+3}{4}}_{\frac{c}{4}, \frac{c+1}{4}, \frac{c+2}{4},\frac{c+3}{4}}; 1\right)&=& \frac{\Gamma(c)\,\Gamma(c-a-b)}{\Gamma(b)\,\Gamma(c-b)}\,\sum_{n=0}^{\infty}  \binom{-a}{n} \left(\frac{\Gamma(b+2n)}{ \Gamma(c-a+2n)}\right)\\
 && \qquad\qquad\qquad\times\,_{2}F_1(a,b+2n;c-a+2n;-1)\nonumber
\eeq  provide $Re(c-a-b)\, >\, 0.$
\end{sloppypar}
\subsection{Main Theorems and its Proofs}
\begin{sloppypar}
The following Lemma is useful to prove our main results.
\end{sloppypar}
\begin{sloppypar}
\blem \label{lem1eq1}
Let, $a,b,c > 0$. Then we have the following:
\begin{enumerate}
\item For $ c > a+b+1$, we have\\
\begin{flushleft}
$\displaystyle\sum_{n=0}^{\infty} \frac{(n+1)\,(a)_n\, \left(\frac{b}{4}\right)_n\, \left(\frac{b+1}{4}\right)_n\, \left(\frac{b+2}{4}\right)_n\,\left(\frac{b+3}{4}\right)_n }
{\left(\frac{c}{4}\right)_n\, \left(\frac{c+1}{4}\right)_n\, \left(\frac{c+2}{4}\right)_n\,\left(\frac{c+3}{4}\right)_n\, (1)_n}$
\end{flushleft}
\begin{eqnarray*}
&=& \frac{\Gamma(c)\, \Gamma(c-a-b)}{\Gamma(b)\,\Gamma(c-b)} \, \bigg( \sum_{n=0}^{\infty}\, \binom{-(a+1)}{n}\, \bigg(\frac{a }{c-a-b-1}\,  \bigg)  \frac{\Gamma(b+4+2n)}{\Gamma(c-a+3+2n)}\cr
&& \qquad \times \, _{2}F_1(a+1,b+4+2n;c-a+3+2n;-1)\,\cr
&& +\,\sum_{n=0}^{\infty}\, \binom{-a}{n}\, \frac{\Gamma(b+2n)}{\Gamma(c-a+2n)}\,   _{2}F_1(a,b+2n;c-a+2n;-1)\bigg)
\end{eqnarray*}
\item For $c > a+b+2$, we have\\
\begin{flushleft}
$\displaystyle\sum_{n=0}^{\infty}\frac{(n+1)^2\,(a)_n\, \left(\frac{b}{4}\right)_n\, \left(\frac{b+1}{4}\right)_n\, \left(\frac{b+2}{4}\right)_n\,\left(\frac{b+3}{4}\right)_n }
{\left(\frac{c}{4}\right)_n\, \left(\frac{c+1}{4}\right)_n\, \left(\frac{c+2}{4}\right)_n\,\left(\frac{c+3}{4}\right)_n\, (1)_n}$
\end{flushleft}
\begin{eqnarray*}
&=& \frac{\Gamma(c)\, \Gamma(c-a-b)}{\Gamma(b)\,\Gamma(c-b)} \, \bigg( \sum_{n=0}^{\infty}\, \binom{-(a+2)}{n}\, \bigg(\frac{(a)_2 }{(c-a-b-2)_{2}}\,  \bigg)  \frac{\Gamma(b+8+2n)}{\Gamma(c-a+6+2n)}\cr
&& \qquad \times \, _{2}F_1(a+2,b+8+2n;c-a+6+2n;-1)\,\cr
&& +3\, \quad \sum_{n=0}^{\infty}\, \binom{-(a+1)}{n}\, \bigg(\frac{a }{c-a-b-1}\,  \bigg)  \frac{\Gamma(b+4+2n)}{\Gamma(c-a+3+2n)}\cr
&& \qquad \times \, _{2}F_1(a+1,b+4+2n;c-a+3+2n;-1)\,\cr
&& +\,\sum_{n=0}^{\infty}\, \binom{-a}{n}\, \frac{\Gamma(b+2n)}{\Gamma(c-a+2n)}\,   _{2}F_1(a,b+2n;c-a+2n;-1)\bigg)
\end{eqnarray*}
\item For $ c > a+b+3$, we have\\
\begin{flushleft}
$\displaystyle\sum_{n=0}^{\infty}\frac{(n+1)^3\,(a)_n\, \left(\frac{b}{4}\right)_n\, \left(\frac{b+1}{4}\right)_n\, \left(\frac{b+2}{4}\right)_n\,\left(\frac{b+3}{4}\right)_n }
{\left(\frac{c}{4}\right)_n\, \left(\frac{c+1}{4}\right)_n\, \left(\frac{c+2}{4}\right)_n\,\left(\frac{c+3}{4}\right)_n\, (1)_n}$
\end{flushleft}
\begin{eqnarray*}
&=& \frac{\Gamma(c)\, \Gamma(c-a-b)}{\Gamma(b)\,\Gamma(c-b)} \, \bigg( \sum_{n=0}^{\infty}\, \binom{-(a+3)}{n}\, \bigg(\frac{(a)_3 }{(c-a-b-3)_{3}}\,  \bigg)  \frac{\Gamma(b+12+2n)}{\Gamma(c-a+9+2n)}\cr
&& \qquad \times \, _{2}F_1(a+3,b+12+2n;c-a+9+2n;-1)\,\cr
&& +6\sum_{n=0}^{\infty}\, \binom{-(a+2)}{n}\, \bigg(\frac{(a)_2 }{(c-a-b-2)_{2}}\,  \bigg)  \frac{\Gamma(b+8+2n)}{\Gamma(c-a+6+2n)}\cr
&& \qquad \times \, _{2}F_1(a+2,b+8+2n;c-a+6+2n;-1)\,\cr
&& +7 \sum_{n=0}^{\infty}\, \binom{-(a+1)}{n}\, \bigg(\frac{a }{c-a-b-1}\,  \bigg)  \frac{\Gamma(b+4+2n)}{\Gamma(c-a+3+2n)}\cr
&& \qquad \times \, _{2}F_1(a+1,b+4+2n;c-a+3+2n;-1)\,\cr
&& +\,\sum_{n=0}^{\infty}\, \binom{-a}{n}\, \frac{\Gamma(b+2n)}{\Gamma(c-a+2n)}\,   _{2}F_1(a,b+2n;c-a+2n;-1)\bigg)
\end{eqnarray*}
\item For $a\neq 1,\, b\neq 1,\,2,\,3,\,4$ and $c >\max\{a+3,  a+b-1\}$, we have\\
\begin{flushleft}
$\displaystyle\sum_{n=0}^{\infty} \frac{(a)_n\, \left(\frac{b}{4}\right)_n\, \left(\frac{b+1}{4}\right)_n\, \left(\frac{b+2}{4}\right)_n\,\left(\frac{b+3}{4}\right)_n }
{\left(\frac{c}{4}\right)_n\, \left(\frac{c+1}{4}\right)_n\, \left(\frac{c+2}{4}\right)_n\,\left(\frac{c+3}{4}\right)_n\, (1)_{n+1}}$
\end{flushleft}
\begin{eqnarray*}
&=& \frac{\Gamma(c)\, \Gamma(c-a-b)}{\Gamma(b)\,\Gamma(c-b)} \, \left( \frac{c-a-b}{a-1} \right)
\times\,\sum_{n=0}^{\infty}\, \binom{-(a-1)}{n}\, \frac{\Gamma(b-4+2n)}{\Gamma(c-a-3+2n)} \\
 &&\times   _{2}F_1(a-1,b-4+2n;c-a-3+2n;-1) - \frac{  (c-4)_{4} }{ (a-1) (b-4)_{4}}
\end{eqnarray*}
\end{enumerate}
\elem
\end{sloppypar}
\bpf(1) Using Pochhammer symbol, we can formulate
\begin{flushleft}
$\displaystyle\sum_{n=0}^{\infty} \frac{(n+1)\,(a)_n\, \left(\frac{b}{4}\right)_n\, \left(\frac{b+1}{4}\right)_n\, \left(\frac{b+2}{4}\right)_n\,\left(\frac{b+3}{4}\right)_n }
{\left(\frac{c}{4}\right)_n\, \left(\frac{c+1}{4}\right)_n\, \left(\frac{c+2}{4}\right)_n\,\left(\frac{c+3}{4}\right)_n\, (1)_n}$
\end{flushleft}
\begin{eqnarray*}
&=&\displaystyle\sum_{n=0}^{\infty} \frac{ n\,(a)_n\, \left(\frac{b}{4}\right)_n\, \left(\frac{b+1}{4}\right)_n\, \left(\frac{b+2}{4}\right)_n\,\left(\frac{b+3}{4}\right)_n }
{\left(\frac{c}{4}\right)_n\, \left(\frac{c+1}{4}\right)_n\, \left(\frac{c+2}{4}\right)_n\,\left(\frac{c+3}{4}\right)_n\, (1)_n}
+\displaystyle\sum_{n=0}^{\infty} \frac{\,(a)_n\, \left(\frac{b}{4}\right)_n\, \left(\frac{b+1}{4}\right)_n\, \left(\frac{b+2}{4}\right)_n\,\left(\frac{b+3}{4}\right)_n }
{\left(\frac{c}{4}\right)_n\, \left(\frac{c+1}{4}\right)_n\, \left(\frac{c+2}{4}\right)_n\,\left(\frac{c+3}{4}\right)_n\, (1)_n} \\
&=&\displaystyle\sum_{n=0}^{\infty} \frac{ (a)_{n+1}\, \left(\frac{b}{4}\right)_{n+1}\, \left(\frac{b+1}{4}\right)_{n+1}\, \left(\frac{b+2}{4}\right)_{n+1}\,\left(\frac{b+3}{4}\right)_{n+1} }
{\left(\frac{c}{4}\right)_{n+1}\, \left(\frac{c+1}{4}\right)_{n+1}\, \left(\frac{c+2}{4}\right)_{n+1}\,\left(\frac{c+3}{4}\right)_{n+1}\, (1)_n}
+\displaystyle\sum_{n=0}^{\infty} \frac{\,(a)_n\, \left(\frac{b}{4}\right)_n\, \left(\frac{b+1}{4}\right)_n\, \left(\frac{b+2}{4}\right)_n\,\left(\frac{b+3}{4}\right)_n }
{\left(\frac{c}{4}\right)_n\, \left(\frac{c+1}{4}\right)_n\, \left(\frac{c+2}{4}\right)_n\,\left(\frac{c+3}{4}\right)_n\, (1)_n} \\
\end{eqnarray*}
Using the formula (\ref{inteq6}) and using the fact that $\Gamma(a+1)= a\Gamma(a)$, the aforementioned equation reduces to
\begin{flushleft}
$\displaystyle\sum_{n=0}^{\infty} \frac{(n+1)\,(a)_n\, \left(\frac{b}{4}\right)_n\, \left(\frac{b+1}{4}\right)_n\, \left(\frac{b+2}{4}\right)_n\,\left(\frac{b+3}{4}\right)_n }
{\left(\frac{c}{4}\right)_n\, \left(\frac{c+1}{4}\right)_n\, \left(\frac{c+2}{4}\right)_n\,\left(\frac{c+3}{4}\right)_n\, (1)_n}$
\end{flushleft}
\begin{eqnarray*}
&=& \frac{\Gamma(c)\, \Gamma(c-a-b)}{\Gamma(b)\,\Gamma(c-b)} \, \bigg( \sum_{n=0}^{\infty}\, \binom{-(a+1)}{n}\, \bigg(\frac{a }{c-a-b-1}\,  \bigg)  \frac{\Gamma(b+4+2n)}{\Gamma(c-a+3+2n)}\cr
&& \qquad \times \, _{2}F_1(a+1,b+4+2n;c-a+3+2n;-1)\,\cr
&& +\,\sum_{n=0}^{\infty}\, \binom{-a}{n}\, \frac{\Gamma(b+2n)}{\Gamma(c-a+2n)}\,   _{2}F_1(a,b+2n;c-a+2n;-1)\bigg)
\end{eqnarray*}
Hence, (1) is proved.

(2) Using ascending factorial notation and by adjusting coefficients suitably, we can write
\begin{flushleft}
  $\displaystyle \sum_{n=0}^{\infty} \frac{(n+1)^2\, (a)_n\, \left(\frac{b}{4}\right)_n\, \left(\frac{b+1}{4}\right)_n\, \left(\frac{b+2}{4}\right)_n\, \left(\frac{b+3}{4}\right)_n}{\left(\frac{c}{4}\right)_n\, \left(\frac{c+1}{4}\right)_n\, \left(\frac{c+2}{4}\right)_n\, \left(\frac{c+3}{4}\right)_n\, (1)_n}$
\end{flushleft}
\begin{eqnarray*}
&=& \displaystyle\sum_{n=0}^{\infty} \frac{(a)_{n+2}\, \left(\frac{b}{4}\right)_{n+2}\, \left(\frac{b+1}{4}\right)_{n+2}\, \left(\frac{b+2}{4}\right)_{n+2}\, \left(\frac{b+3}{4}\right)_{n+2}}{\left(\frac{c}{4}\right)_{n+2}\, \left(\frac{c+1}{4}\right)_{n+2}\, \left(\frac{c+2}{4}\right)_{n+2}\, \left(\frac{c+3}{4}\right)_{n+2}\, (1)_{n}} \\ \\
&&\qquad\qquad + 3 \displaystyle\sum_{n=0}^{\infty} \frac{(a)_{n+1}\, \left(\frac{b}{4}\right)_{n+1}\, \left(\frac{b+1}{4}\right)_{n+1}\, \left(\frac{b+2}{4}\right)_{n+1}\, \left(\frac{b+3}{4}\right)_{n+1}}{\left(\frac{c}{4}\right)_{n+1}\, \left(\frac{c+1}{4}\right)_{n+1}\, \left(\frac{c+2}{4}\right)_{n+1}\, \left(\frac{c+3}{4}\right)_{n+1}\, (1)_{n}} \\  \\
  && \qquad \qquad \qquad+ \displaystyle\sum_{n=0}^{\infty} \frac{\,(a)_n\, \left(\frac{b}{4}\right)_n\, \left(\frac{b+1}{4}\right)_n\, \left(\frac{b+2}{4}\right)_n\,\left(\frac{b+3}{4}\right)_n }
{\left(\frac{c}{4}\right)_n\, \left(\frac{c+1}{4}\right)_n\, \left(\frac{c+2}{4}\right)_n\,\left(\frac{c+3}{4}\right)_n\, (1)_n}
\end{eqnarray*}
Using the formula (\ref{inteq6}) and using the fact that $\Gamma(a+1)= a\Gamma(a)$, the aforementioned equation reduces to
\begin{flushleft}
$\displaystyle\sum_{n=0}^{\infty}\frac{(n+1)^2\,(a)_n\, \left(\frac{b}{4}\right)_n\, \left(\frac{b+1}{4}\right)_n\, \left(\frac{b+2}{4}\right)_n\,\left(\frac{b+3}{4}\right)_n }
{\left(\frac{c}{4}\right)_n\, \left(\frac{c+1}{4}\right)_n\, \left(\frac{c+2}{4}\right)_n\,\left(\frac{c+3}{4}\right)_n\, (1)_n}$
\end{flushleft}
\begin{eqnarray*}
&=& \frac{\Gamma(c)\, \Gamma(c-a-b)}{\Gamma(b)\,\Gamma(c-b)} \, \bigg( \sum_{n=0}^{\infty}\, \binom{-(a+2)}{n}\, \bigg(\frac{(a)_2 }{(c-a-b-2)_{2}}\,  \bigg)  \frac{\Gamma(b+8+2n)}{\Gamma(c-a+6+2n)}\cr \cr
&& \qquad \times \, _{2}F_1(a+2,b+8+2n;c-a+6+2n;-1)\,\cr  \cr
&& +\quad 3\sum_{n=0}^{\infty}\, \binom{-(a+1)}{n}\, \bigg(\frac{a }{c-a-b-1}\,  \bigg)  \frac{\Gamma(b+4+2n)}{\Gamma(c-a+3+2n)}\cr  \cr
&& \qquad \times \, _{2}F_1(a+1,b+4+2n;c-a+3+2n;-1)\,\cr \cr
&& +\,\sum_{n=0}^{\infty}\, \binom{-a}{n}\, \frac{\Gamma(b+2n)}{\Gamma(c-a+2n)}\,   _{2}F_1(a,b+2n;c-a+2n;-1)\bigg)
\end{eqnarray*}
which completes the proof of (2).

(3)  Using $(n+1)^3 = n(n-1)(n-2)+6(n)(n-1)+7(n)+1$, we can easily obtain that\\

\begin{flushleft}
  $\displaystyle \sum_{n=0}^{\infty} \frac{(n+1)^3\, (a)_n\, \left(\frac{b}{4}\right)_n\, \left(\frac{b+1}{4}\right)_n\, \left(\frac{b+2}{4}\right)_n\, \left(\frac{b+3}{4}\right)_n}{\left(\frac{c}{4}\right)_n\, \left(\frac{c+1}{4}\right)_n\, \left(\frac{c+2}{4}\right)_n\, \left(\frac{c+3}{4}\right)_n\, (1)_n}$
\end{flushleft}
\begin{eqnarray*}
&=& \displaystyle\sum_{n=0}^{\infty} \frac{(a)_{n+3}\, \left(\frac{b}{4}\right)_{n+3}\, \left(\frac{b+1}{4}\right)_{n+3}\, \left(\frac{b+2}{4}\right)_{n+3}\, \left(\frac{b+3}{4}\right)_{n+3}}{\left(\frac{c}{4}\right)_{n+3}\, \left(\frac{c+1}{4}\right)_{n+3}\, \left(\frac{c+2}{4}\right)_{n+3}\, \left(\frac{c+3}{4}\right)_{n+3}\, (1)_{n}} \\ \\
&& \qquad +6\displaystyle\sum_{n=0}^{\infty} \frac{(a)_{n+2}\, \left(\frac{b}{4}\right)_{n+2}\, \left(\frac{b+1}{4}\right)_{n+2}\, \left(\frac{b+2}{4}\right)_{n+2}\, \left(\frac{b+3}{4}\right)_{n+2}}{\left(\frac{c}{4}\right)_{n+2}\, \left(\frac{c+1}{4}\right)_{n+2}\, \left(\frac{c+2}{4}\right)_{n+2}\, \left(\frac{c+3}{4}\right)_{n+2}\, (1)_{n}} \\ \\
&&\qquad\qquad + 7 \displaystyle\sum_{n=0}^{\infty} \frac{(a)_{n+1}\, \left(\frac{b}{4}\right)_{n+1}\, \left(\frac{b+1}{4}\right)_{n+1}\, \left(\frac{b+2}{4}\right)_{n+1}\, \left(\frac{b+3}{4}\right)_{n+1}}{\left(\frac{c}{4}\right)_{n+1}\, \left(\frac{c+1}{4}\right)_{n+1}\, \left(\frac{c+2}{4}\right)_{n+1}\, \left(\frac{c+3}{4}\right)_{n+1}\, (1)_{n}} \\ \\
  && \qquad \qquad \qquad+ \displaystyle\sum_{n=0}^{\infty} \frac{\,(a)_n\, \left(\frac{b}{4}\right)_n\, \left(\frac{b+1}{4}\right)_n\, \left(\frac{b+2}{4}\right)_n\,\left(\frac{b+3}{4}\right)_n }
{\left(\frac{c}{4}\right)_n\, \left(\frac{c+1}{4}\right)_n\, \left(\frac{c+2}{4}\right)_n\,\left(\frac{c+3}{4}\right)_n\, (1)_n}
\end{eqnarray*}
Using the formula (\ref{inteq6}) and using the fact that $\Gamma(a+1)= a\Gamma(a)$, the aforementioned equation reduces to
\begin{flushleft}
$\displaystyle\sum_{n=0}^{\infty}\frac{(n+1)^3\,(a)_n\, \left(\frac{b}{4}\right)_n\, \left(\frac{b+1}{4}\right)_n\, \left(\frac{b+2}{4}\right)_n\,\left(\frac{b+3}{4}\right)_n }
{\left(\frac{c}{4}\right)_n\, \left(\frac{c+1}{4}\right)_n\, \left(\frac{c+2}{4}\right)_n\,\left(\frac{c+3}{4}\right)_n\, (1)_n}$
\end{flushleft}
\begin{eqnarray*}
&=& \frac{\Gamma(c)\, \Gamma(c-a-b)}{\Gamma(b)\,\Gamma(c-b)} \, \bigg( \sum_{n=0}^{\infty}\, \binom{-(a+3)}{n}\, \bigg(\frac{(a)_3 }{(c-a-b-3)_{3}}\,  \bigg)  \frac{\Gamma(b+12+2n)}{\Gamma(c-a+9+2n)}\cr \cr
&& \qquad \qquad\qquad\times \, _{2}F_1(a+3,b+12+2n;c-a+9+2n;-1)\,\cr \cr
&& \qquad+6 \sum_{n=0}^{\infty}\, \binom{-(a+2)}{n}\, \bigg(\frac{(a)_2 }{(c-a-b-2)_{2}}\,  \bigg)  \frac{\Gamma(b+8+2n)}{\Gamma(c-a+6+2n)}\cr \cr
&& \qquad\qquad \qquad\times \, _{2}F_1(a+2,b+8+2n;c-a+6+2n;-1)\,\cr \cr
&& \qquad+7\sum_{n=0}^{\infty}\, \binom{-(a+1)}{n}\, \bigg(\frac{a }{c-a-b-1}\,  \bigg)  \frac{\Gamma(b+4+2n)}{\Gamma(c-a+3+2n)}\cr \cr
&& \qquad \qquad\qquad\times \, _{2}F_1(a+1,b+4+2n;c-a+3+2n;-1)\,\cr \cr
&& \qquad+\sum_{n=0}^{\infty}\, \binom{-a}{n}\, \frac{\Gamma(b+2n)}{\Gamma(c-a+2n)}\,   _{2}F_1(a,b+2n;c-a+2n;-1)\bigg)
\end{eqnarray*}
which completes the proof.

(4)  Let, $a\neq 1$, $b\neq 1,\, 2,\, 3,\, 4$, $c > \max \{a+3, a+b-1\}$, it is found that
\begin{flushleft}
$\displaystyle\sum_{n=0}^{\infty} \frac{(a)_n\, \left(\frac{b}{4}\right)_n\, \left(\frac{b+1}{4}\right)_n\, \left(\frac{b+2}{4}\right)_n\,\left(\frac{b+3}{4}\right)_n }
{\left(\frac{c}{4}\right)_n\, \left(\frac{c+1}{4}\right)_n\, \left(\frac{c+2}{4}\right)_n\,\left(\frac{c+3}{4}\right)_n\, (1)_{n+1}}$
\end{flushleft}
\begin{eqnarray*}
&=& \left(\frac{(c-4)\, (c-3)\, (c-2)\, (c-1)}{(a-1)\, (b-4)\, (b-3)\, (b-2)\, (b-1)}\right)\\
&&\qquad\qquad\times\sum_{n=0}^{\infty} \frac{\,(a-1)_n\, \left(\frac{b}{4}-1 \right)_n\, \left(\frac{b+1}{4}-1 \right)_n\, \left(\frac{b+2}{4}-1 \right)_n\,\left(\frac{b+3}{4}-1\right)_n }
{\left(\frac{c}{4} - 1 \right)_n\, \left(\frac{c+1}{4} -1 \right)_n\, \left(\frac{c+2}{4} - 1\right)_n\,\left(\frac{c+3}{4} -1\right)_n\, (1)_n} \\
&=& \frac{\Gamma(c)\, \Gamma(c-a-b)}{\Gamma(b)\,\Gamma(c-b)} \, \left( \frac{c-a-b}{a-1} \right)
\times\,\sum_{n=0}^{\infty}\, \binom{-(a-1)}{n}\, \frac{\Gamma(b-4+2n)}{\Gamma(c-a-3+2n)} \\
 &&\qquad\qquad\times   _{2}F_1(a-1,b-4+2n;c-a-3+2n;-1) - \frac{  (c-4)_{4} }{ (a-1) (b-4)_{4}}
\end{eqnarray*}
Hence, the desired result follows.
\epf%
\subsubsection{Starlikeness Property}
\bthm\label{thm1eq0}
Let, $a,\, b \in {\Bbb C} \backslash \{ 0 \} $,\, $c > 0$\, and $c > |a|+|b|+1.$ A sufficient condition for the function $z\, _5F_4\left(^{a,\,\frac{b}{4},\, \frac{b+1}{4}\, \frac{b+2}{4}\, \frac{b+3}{4}}_{\frac{c}{4},\, \frac{c+1}{4}\, \frac{c+2}{4}\, \frac{c+3}{4}}; z\right) $ to belong to the class $ \es^{*}_{\lambda}, \,  0 < \lambda  \leq 1 $ is that
\beq\label{thm1eq1}
&&\frac{\Gamma(c)\, \Gamma(c-|a|-|b|)}{\Gamma(|b|)\, \Gamma(c-|b|)}\bigg[\sum_{n=0}^{\infty}\, \binom{-(|a|+1)}{n}\, \bigg(\frac{|a| }{c-|a|-|b|-1}\,  \bigg)  \frac{\Gamma(|b|+4+2n)}{\Gamma(c-|a|+3+2n)}\cr
&& \qquad \qquad\qquad\times \, _{2}F_1(|a|+1,|b|+4+2n;c-|a|+3+2n;-1)\,\cr
&& \qquad+\lambda\sum_{n=0}^{\infty}\, \binom{-|a|}{n}\, \frac{\Gamma(|b|+2n)}{\Gamma(c-|a|+2n)}\,   _{2}F_1(|a|,|b|+2n;c-|a|+2n;-1)\bigg] \leq 2\lambda.
\eeq
\ethm
\bpf  Let, $f(z)=z\, _5F_4\left(^{a,\,\frac{b}{4},\, \frac{b+1}{4}\, \frac{b+2}{4}\, \frac{b+3}{4}}_{\frac{c}{4},\, \frac{c+1}{4}\, \frac{c+2}{4}\, \frac{c+3}{4}}; z\right) $, then, by the equation (\ref{inteq2}), it is enough to show that
\begin{eqnarray*}
  T &=& \sum_{n=2}^{\infty}(n+\lambda-1)|A_n|\leq \lambda.
\end{eqnarray*}
Using the fact $|(a)_n|\leq  (|a|)_n$, one can get
\begin{eqnarray*}
  T &\leq& \sum_{n=0}^{\infty} ((n+1)+(\lambda-1))\left(\frac{(|a|)_{n}\left(\frac{|b|}{4}\right)_{n}\, \left(\frac{|b|+1}{4}\right)_{n}\, \left(\frac{|b|+2}{4}\right)_{n}\, \left(\frac{|b|+3}{4}\right)_{n}}{\left(\frac{c}{4}\right)_{n}\, \left(\frac{c+1}{4}\right)_{n}\, \left(\frac{c+2}{4}\right)_{n}\, \left(\frac{c+3}{4}\right)_{n}\,(1)_{n}}\right)-\lambda.
\end{eqnarray*}
Using (\ref{inteq6}) and the result (1) of Lemma \ref{lem1eq1} in the aforesaid equation, we get
\begin{eqnarray*}
T &\leq& \left(\frac{\Gamma(c)\, \Gamma(c-|a|-|b|)}{\Gamma(|b|)\, \Gamma(c-|b|)}\right)\bigg(\sum_{n=0}^{\infty}\, \binom{-(|a|+1)}{n}\, \bigg(\frac{|a| }{c-|a|-|b|-1}\,  \bigg)  \frac{\Gamma(|b|+4+2n)}{\Gamma(c-|a|+3+2n)}\cr
&& \qquad \times \, _{2}F_1(|a|+1,|b|+4+2n;c-|a|+3+2n;-1)\,\cr
&& +\,\sum_{n=0}^{\infty}\, \binom{-|a|}{n}\, \frac{\Gamma(|b|+2n)}{\Gamma(c-|a|+2n)}\,   _{2}F_1(|a|,|b|+2n;c-|a|+2n;-1)\cr
  && +\, (\lambda-1)\sum_{n=0}^{\infty}\, \binom{-|a|}{n}\, \frac{\Gamma(|b|+2n)}{\Gamma(c-|a|+2n)}\,   _{2}F_1(|a|,|b|+2n;c-|a|+2n;-1)\bigg)-\lambda.
\end{eqnarray*}
Because of (\ref{thm1eq1}), the aforementioned expression is bounded by $\lambda$, and hence,
\begin{eqnarray*}
  T &\leq& \left(\frac{\Gamma(c)\, \Gamma(c-|a|-|b|)}{\Gamma(|b|)\, \Gamma(c-|b|)}\right)\bigg(\sum_{n=0}^{\infty}\, \binom{-(|a|+1)}{n}\, \bigg(\frac{|a| }{c-|a|-|b|-1}\,  \bigg)  \frac{\Gamma(|b|+4+2n)}{\Gamma(c-|a|+3+2n)}\cr
&& \qquad \times \, _{2}F_1(|a|+1,|b|+4+2n;c-|a|+3+2n;-1)\,\cr
&& +\, \lambda\,\sum_{n=0}^{\infty}\, \binom{-|a|}{n}\, \frac{\Gamma(|b|+2n)}{\Gamma(c-|a|+2n)}\,   _{2}F_1(|a|,|b|+2n;c-|a|+2n;-1)\bigg)-\lambda \leq \lambda.
\end{eqnarray*}
Therefore, $z\, _5F_4\left(^{a,\,\frac{b}{4},\, \frac{b+1}{4}\, \frac{b+2}{4}\, \frac{b+3}{4}}_{\frac{c}{4},\, \frac{c+1}{4}\, \frac{c+2}{4}\, \frac{c+3}{4}}; z\right)$ belongs to the class $\es^{*}_{\lambda}. $
\epf
\bthm\label{thm2eq001}
 Let, $a,\, b \in {\Bbb C} \backslash \{ 0 \} ,\, c > 0,\, \, |a|\neq1,\, |b| \neq1,\,2,\,3,\,4$ and $c > \max\{|a|+3, |a|+|b|-1\}$. For  $ 0 < \lambda \leq 1$ and $0 \leq \beta < 1$, assume that
\beq\label{thm2eq1}
&&\left(\frac{\Gamma(c)\,\Gamma(c-|a|-|b|)}{\Gamma(|b|)\, \Gamma(c-|b|)}\right)\bigg[ \left( \frac{(\lambda-1)\,(c-|a|-|b|)}{|a|-1} \right)\sum_{n=0}^{\infty}\, \binom{-(|a|-1)}{n}\cr
&&\qquad\times\, \frac{\Gamma(|b|+4+2n)}{\Gamma(c-|a|-3+2n)}\,  _{2}F_1(|a|-1,|b|+4+2n;c-|a|-3+2n;-1)\cr   &&\qquad\qquad\qquad\,+\,\,\sum_{n=0}^{\infty}\, \binom{-|a|}{n}\, \frac{\Gamma(|b|+2n)}{\Gamma(c-|a|+2n)}\,   _{2}F_1(|a|,|b|+2n;c-|a|+2n;-1)\bigg]\nonumber \cr
&&\qquad\qquad\qquad\qquad\qquad\qquad\qquad\qquad\leq \lambda\left(1+\frac{1}{2(1-\beta)}\right)+\left(\frac{(\lambda-1)\,(c-4)_4}{(|a|-1)(|b|-4)_4}\right).
\eeq
 Then, the integral  operator $\mathcal{I}^{a,\,\frac{b}{4},\, \frac{b+1}{4}\, \frac{b+2}{4}\, \frac{b+3}{4}}_{\frac{c}{4},\, \frac{c+1}{4}\, \frac{c+2}{4}\, \frac{c+2}{4}}(f)(z)$ maps $ \mathcal{R}(\beta)$ into $\es^{*}_{\lambda}$.
\ethm
\bpf Let, $a,\, b \in {\Bbb C} \backslash \{ 0 \} ,\, c > 0,\, |a|\neq1,\, |b| \neq1,\,2,\,3,\, 4$, \,  $c > \max\{|a|+3, |a|+|b|-1\}$,\, $ 0 < \lambda \leq 1$ and $0 \leq \beta < 1$.

Consider, the integral operator $\mathcal{I}^{a,\,\frac{b}{4},\, \frac{b+1}{4}\, \frac{b+2}{4}\, \frac{b+3}{4}}_{\frac{c}{4},\, \frac{c+1}{4}\, \frac{c+2}{4}\, \frac{c+2}{4}}(f)(z)$ defined by (\ref{inteq7}). According to (\ref{inteq2}), it is enough to show that
\begin{eqnarray}\label{thm2eq002}
  T &=& \sum_{n=2}^{\infty}(n+\lambda-1)|A_n|\leq \lambda,
\end{eqnarray}
where, $A_n$ is given by (\ref{inteq007}). Then, it is derived that
\begin{eqnarray*}
  T &\leq&  \sum_{n=2}^{\infty} (n+(\lambda-1)) \left(\frac{(|a|)_{n-1}\left(\frac{|b|}{4}\right)_{n-1}\, \left(\frac{|b|+1}{4}\right)_{n-1}\, \left(\frac{|b|+2}{4}\right)_{n-1}\, \left(\frac{|b|+3}{4}\right)_{n-1}}{\left(\frac{c}{4}\right)_{n-1}\, \left(\frac{c+1}{4}\right)_{n-1}\, \left(\frac{c+2}{4}\right)_{n-1}\, \left(\frac{c+3}{4}\right)_{n}(1)_{n-1}}\right) |a_n|
\end{eqnarray*}
Using (\ref{inteq3})  in the aforementioned equation, we have
\begin{eqnarray*}
  T&\leq&  2(1-\beta)\bigg(\sum_{n=0}^{\infty} \left(\frac{(|a|)_{n}\left(\frac{|b|}{4}\right)_{n}\, \left(\frac{|b|+1}{4}\right)_{n}\, \left(\frac{|b|+2}{4}\right)_{n}\, \left(\frac{|b|+3}{4}\right)_{n}}{\left(\frac{c}{4}\right)_{n}\, \left(\frac{c+1}{4}\right)_{n}\, \left(\frac{c+2}{4}\right)_{n}\, \left(\frac{c+3}{4}\right)_{n}(1)_{n}}\right) -1 \cr
  && \qquad + (\lambda-1)\sum_{n=1}^{\infty} \left(\frac{(|a|)_{n}\left(\frac{|b|}{4}\right)_{n}\, \left(\frac{|b|+1}{4}\right)_{n}\, \left(\frac{|b|+2}{4}\right)_{n}\, \left(\frac{|b|+3}{4}\right)_{n}}{\left(\frac{c}{4}\right)_{n}\, \left(\frac{c+1}{4}\right)_{n}\, \left(\frac{c+2}{4}\right)_{n}\, \left(\frac{c+3}{4}\right)_{n}(1)_{n}}\right) \left(\frac{1}{n+1}\right)  \bigg):=T_1.
\end{eqnarray*}
Using the formula (\ref{inteq6}) and the result (4) of Lemma \ref{lem1eq1}, we find that
\begin{eqnarray*}
  T_1 &\leq& 2(1-\beta)\bigg[\left(\frac{\Gamma(c)\,\Gamma(c-|a|-|b|)}{\Gamma(|b|)\, \Gamma(c-|b|)}\right)\bigg( \left( \frac{(\lambda-1)\,(c-|a|-|b|)}{|a|-1} \right)\,\sum_{n=0}^{\infty}\, \binom{-(|a|-1)}{n}\cr
&&\times\, \frac{\Gamma(|b|-4+2n)}{\Gamma(c-|a|-3+2n)}\, _{2}F_1(|a|-1,|b|-4+2n;c-|a|-3+2n;-1)\cr   &&\qquad\qquad\,+\,\,\sum_{n=0}^{\infty}\, \binom{-|a|}{n}\, \frac{\Gamma(|b|+2n)}{\Gamma(c-|a|+2n)}\,   _{2}F_1(|a|,|b|+2n;c-|a|+2n;-1)\bigg)\cr   &&\qquad\qquad\qquad\, - \frac{(\lambda-1)(c-4)_4}{(|a|-1)(|b|-4)_4} -\lambda \bigg].
\end{eqnarray*}
Under the condition (\ref{thm2eq1})
\begin{eqnarray*}
&&2(1-\beta)\bigg[\left(\frac{\Gamma(c)\,\Gamma(c-|a|-|b|)}{\Gamma(|b|)\, \Gamma(c-|b|)}\right)\bigg( \left( \frac{(\lambda-1)\,(c-|a|-|b|)}{|a|-1} \right)\,\sum_{n=0}^{\infty}\, \binom{-(|a|-1)}{n}\cr
&&\qquad\times\, \frac{\Gamma(|b|-4+2n)}{\Gamma(c-|a|-3+2n)}\, _{2}F_1(|a|-1,|b|-4+2n;c-|a|-3+2n;-1)\cr
&&\qquad\qquad\qquad\,+\,\,\sum_{n=0}^{\infty}\, \binom{-|a|}{n}\, \frac{\Gamma(|b|+2n)}{\Gamma(c-|a|+2n)}\,   _{2}F_1(|a|,|b|+2n;c-|a|+2n;-1)\bigg)\cr   &&\qquad\qquad\qquad\qquad\, - \frac{(\lambda-1)(c-4)_4}{(|a|-1)(|b|-4)_4} -\lambda \bigg] \leq \lambda
\end{eqnarray*}
Thus, we have the inequalities $T \leq T_1 \leq \lambda $,  and hence, (\ref{thm2eq002}) is held. Therefore, we conclude that the operator $\mathcal{I}^{a,\,\frac{b}{4},\, \frac{b+1}{4}\, \frac{b+2}{4}\, \frac{b+3}{4}}_{\frac{c}{4},\, \frac{c+1}{4}\, \frac{c+2}{4}\, \frac{c+2}{4}}(f)(z)$ maps $ \mathcal{R}(\beta)$ into $\es^{*}_{\lambda}$, which completes the proof of the theorem.
\epf

When, $\lambda =1 $, we get the following result from Theorem \ref{thm2eq001}.
\bcor
 Let, $a,\, b \in {\Bbb C} \backslash \{ 0 \} ,\, c > 0$,\, $c > |a|+|b|$ and $0 \leq \beta < 1$. Assume that
\beq\label{cor2eq1}
&&\left(\frac{\Gamma(c)\,\Gamma(c-|a|-|b|)}{\Gamma(|b|)\, \Gamma(c-|b|)}\right)\,\sum_{n=0}^{\infty}\, \binom{-|a|}{n}\, \frac{\Gamma(|b|+2n)}{\Gamma(c-|a|+2n)}\,\nonumber \cr
&&\qquad\qquad\qquad\qquad\qquad\qquad\times   _{2}F_1(|a|,|b|+2n;c-|a|+2n;-1) \leq 1+\frac{1}{2(1-\beta)}.
\eeq
 Then, the integral  operator $\mathcal{I}^{a,\,\frac{b}{4},\, \frac{b+1}{4}\, \frac{b+2}{4}\, \frac{b+3}{4}}_{\frac{c}{4},\, \frac{c+1}{4}\, \frac{c+2}{4}\, \frac{c+2}{4}}(f)(z)$ maps $ \mathcal{R}(\beta)$ into $\es^{*}$.
\ecor
\bthm\label{thm3eq0}  Let, $a, b \in {\Bbb C} \backslash \{ 0 \} ,\, c > 0,\,$ $c > |a|+|b|+2$ and $0 < \lambda \leq 1$. If
 \beq\label{thm3eq1}
&&\left(\frac{\Gamma(c)\,\Gamma(c-|a|-b)}{\Gamma(c-a)\, \Gamma(c-b)}\right) \cr
&&\qquad\times\bigg[\bigg(\frac{(|a|)_2 }{(c-|a|-|b|-2)_2}\,  \bigg)\sum_{n=0}^{\infty}\, \binom{-(|a|+2)}{n}\, \frac{\Gamma(|b|+8+2n)}{\Gamma(c-|a|+6+2n)}\cr
&& \qquad\qquad \qquad\qquad\qquad\times \, _{2}F_1(|a|+2,|b|+8+2n;c-|a|+6+2n;-1)\,\cr
&& \qquad+ \bigg(\frac{(\lambda+2)|a| }{c-|a|-|b|-1}\,  \bigg)\sum_{n=0}^{\infty}\, \binom{-(|a|+1)}{n}\, \frac{\Gamma(|b|+4+2n)}{\Gamma(c-|a|+3+2n)}\cr
&& \qquad\qquad \qquad\qquad\qquad\times \, _{2}F_1(|a|+1,|b|+4+2n;c-|a|+3+2n;-1)\,\cr
&& \qquad+\, \lambda\,\sum_{n=0}^{\infty}\, \binom{-|a|}{n}\, \frac{\Gamma(|b|+2n)}{\Gamma(c-|a|+2n)}\,   _{2}F_1(|a|,|b|+2n;c-|a|+2n;-1)\bigg] \leq 2\lambda.
\eeq
then, the integral  operator $\mathcal{I}^{a,\,\frac{b}{4},\, \frac{b+1}{4}\, \frac{b+2}{4}\, \frac{b+3}{4}}_{\frac{c}{4},\, \frac{c+1}{4}\, \frac{c+2}{4}\, \frac{c+2}{4}}(f)(z)$ maps $\mathcal{S}$ to $\es^{*}_{\lambda}$.
\ethm
\bpf  Let,  $a, b \in {\Bbb C} \backslash \{ 0 \} ,\, c > 0,\,$ $c > |a|+|b|+2$ and $0 < \lambda \leq 1$.

If the integral  operator $\mathcal{I}^{a,\,\frac{b}{4},\, \frac{b+1}{4}\, \frac{b+2}{4}\, \frac{b+3}{4}}_{\frac{c}{4},\, \frac{c+1}{4}\, \frac{c+2}{4}\, \frac{c+2}{4}}(f)(z)$ is defined by (\ref{inteq7}), in view of (\ref{inteq2}), it is enough to show that
\begin{eqnarray}
  T &:=& \sum_{n=2}^{\infty}(n+\lambda-1)|A_n|\leq \lambda,
\end{eqnarray}
where, $A_n$ is given by (\ref{inteq007}).
Using the fact $|(a)_n|\leq  (|a|)_n$ and the equation ($\ref{inteq2}$) in the aforementioned equation, it is derived that
\begin{eqnarray*}
  T &\leq & \sum_{n=2}^{\infty} (n+(\lambda-1)) \left(\frac{(|a|)_{n-1}\left(\frac{|b|}{4}\right)_{n-1}\, \left(\frac{|b|+1}{4}\right)_{n-1}\, \left(\frac{|b|+2}{4}\right)_{n-1}\, \left(\frac{|b|+3}{4}\right)_{n-1}}{\left(\frac{c}{4}\right)_{n-1}\, \left(\frac{c+1}{4}\right)_{n-1}\, \left(\frac{c+2}{4}\right)_{n-1}\, \left(\frac{c+3}{4}\right)_{n-1}(1)_{n-1}}\right)|a_n|\\
  & = & \sum_{n=0}^{\infty} \left(\frac{(n+1)^2 \,(|a|)_{n}\left(\frac{|b|}{4}\right)_{n}\, \left(\frac{|b|+1}{4}\right)_{n}\, \left(\frac{|b|+2}{4}\right)_{n}\, \left(\frac{|b|+3}{4}\right)_{n}}{\left(\frac{c}{4}\right)_{n}\, \left(\frac{c+1}{4}\right)_{n}\, \left(\frac{c+2}{4}\right)_{n}\, \left(\frac{c+3}{4}\right)_{n}(1)_{n}}\right)-1\\
  &&\qquad+(\lambda-1) \,\sum_{n=0}^{\infty}\left(\frac{ (n+1)\,(|a|)_{n}\left(\frac{|b|}{4}\right)_{n}\, \left(\frac{|b|+1}{4}\right)_{n}\, \left(\frac{|b|+2}{4}\right)_{n}\, \left(\frac{|b|+3}{4}\right)_{n}}{\left(\frac{c}{4}\right)_{n}\, \left(\frac{c+1}{4}\right)_{n}\, \left(\frac{c+2}{4}\right)_{n}\, \left(\frac{c+3}{4}\right)_{n}(1)_{n}}\right)-(\lambda-1).
\end{eqnarray*}
Using (1) and (2) of Lemma \ref{lem1eq1}, we find that
\begin{eqnarray*}
T &\leq&\left(\frac{\Gamma(c)\,\Gamma(c-|a|-b)}{\Gamma(c-a)\, \Gamma(c-b)}\right) \cr
&&\qquad\times\bigg[\bigg(\frac{(|a|)_2 }{(c-|a|-|b|-2)_2}\,  \bigg)\sum_{n=0}^{\infty}\, \binom{-(|a|+2)}{n}\, \frac{\Gamma(|b|+8+2n)}{\Gamma(c-|a|+6+2n)}\cr
&& \qquad\qquad \qquad\qquad\qquad\times \, _{2}F_1(|a|+2,|b|+8+2n;c-|a|+6+2n;-1)\,\cr
&& \qquad+ \bigg(\frac{(\lambda+2)|a| }{c-|a|-|b|-1}\,  \bigg)\sum_{n=0}^{\infty}\, \binom{-(|a|+1)}{n}\, \frac{\Gamma(|b|+4+2n)}{\Gamma(c-|a|+3+2n)}\cr
&& \qquad\qquad \qquad\qquad\qquad\times \, _{2}F_1(|a|+1,|b|+4+2n;c-|a|+3+2n;-1)\,\cr
&& \qquad+\, \lambda\,\sum_{n=0}^{\infty}\, \binom{-|a|}{n}\, \frac{\Gamma(|b|+2n)}{\Gamma(c-|a|+2n)}\,   _{2}F_1(|a|,|b|+2n;c-|a|+2n;-1)\bigg]-\lambda
\end{eqnarray*}
By (\ref{thm3eq1}), the above expression is bounded above by $\lambda$, and hence,
\begin{eqnarray*}
&&\left(\frac{\Gamma(c)\,\Gamma(c-|a|-b)}{\Gamma(c-a)\, \Gamma(c-b)}\right) \cr
&&\qquad\times\bigg[\bigg(\frac{(|a|)_2 }{(c-|a|-|b|-2)_2}\,  \bigg)\sum_{n=0}^{\infty}\, \binom{-(|a|+2)}{n}\, \frac{\Gamma(|b|+8+2n)}{\Gamma(c-|a|+6+2n)}\cr
&& \qquad\qquad \qquad\qquad\qquad\times \, _{2}F_1(|a|+2,|b|+8+2n;c-|a|+6+2n;-1)\,\cr
&& \qquad+ \bigg(\frac{(\lambda+2)|a| }{c-|a|-|b|-1}\,  \bigg)\sum_{n=0}^{\infty}\, \binom{-(|a|+1)}{n}\, \frac{\Gamma(|b|+4+2n)}{\Gamma(c-|a|+3+2n)}\cr
&& \qquad\qquad \qquad\qquad\qquad\times \, _{2}F_1(|a|+1,|b|+4+2n;c-|a|+3+2n;-1)\,\cr
&& \qquad+\, \lambda\,\sum_{n=0}^{\infty}\, \binom{-|a|}{n}\, \frac{\Gamma(|b|+2n)}{\Gamma(c-|a|+2n)}\,   _{2}F_1(|a|,|b|+2n;c-|a|+2n;-1)\bigg]-\lambda \leq \lambda
\end{eqnarray*}
Under the stated condition, the integral operator $\mathcal{I}^{a,\,\frac{b}{4},\, \frac{b+1}{4}\, \frac{b+2}{4}\, \frac{b+3}{4}}_{\frac{c}{4},\, \frac{c+1}{4}\, \frac{c+2}{4}\, \frac{c+2}{4}}(f)(z)$ maps $\es$ into $\es^{*}_{\lambda}$.
\epf
\subsubsection{Convexity}
\bthm\label{thm10eq1}
 Let,  $a, b \in {\Bbb C} \backslash \{ 0 \} ,\, c > 0,\,$ $c > |a|+|b|+2$ and $0 < \lambda \leq 1$. A sufficient condition for the function $z\, _5F_4\left(^{a,\,\frac{b}{4},\, \frac{b+1}{4}\, \frac{b+2}{4}\, \frac{b+3}{4}}_{\frac{c}{4},\, \frac{c+1}{4}\, \frac{c+2}{4}\, \frac{c+3}{4}}; z\right)$ to belong to the class $ \mathcal{C}_{\lambda}$ is that
 \beq\label{thm10eq10}
&&\left(\frac{\Gamma(c)\,\Gamma(c-|a|-b)}{\Gamma(c-a)\, \Gamma(c-b)}\right) \cr
&&\qquad\times\bigg[\bigg(\frac{(|a|)_2 }{(c-|a|-|b|-2)_2}\,  \bigg)\sum_{n=0}^{\infty}\, \binom{-(|a|+2)}{n}\, \frac{\Gamma(|b|+8+2n)}{\Gamma(c-|a|+6+2n)}\cr
&& \qquad\qquad \qquad\qquad\qquad\times \, _{2}F_1(|a|+2,|b|+8+2n;c-|a|+6+2n;-1)\,\cr
&& \qquad+ \bigg(\frac{(\lambda+2)|a| }{c-|a|-|b|-1}\,  \bigg)\sum_{n=0}^{\infty}\, \binom{-(|a|+1)}{n}\, \frac{\Gamma(|b|+4+2n)}{\Gamma(c-|a|+3+2n)}\cr
&& \qquad\qquad \qquad\qquad\qquad\times \, _{2}F_1(|a|+1,|b|+4+2n;c-|a|+3+2n;-1)\,\cr
&& \qquad+\, \lambda\,\sum_{n=0}^{\infty}\, \binom{-|a|}{n}\, \frac{\Gamma(|b|+2n)}{\Gamma(c-|a|+2n)}\,   _{2}F_1(|a|,|b|+2n;c-|a|+2n;-1)\bigg] \leq 2\lambda.\nonumber
   \eeq
\ethm
\bpf
The proof is similar to Theorem \ref{thm3eq0}. So, we omit the details.
\epf
\bthm\label{thm11eq0}  Let, $a,\,b \in {\Bbb C} \backslash \{ 0 \} ,\, c > 0,\, c > |a|+|b|+1$, and $ 0 < \lambda \leq 1$.  For $0 \leq \beta <1 $, assume that
 \beq\label{thm11eq1}
\left(\frac{\Gamma(c)\,\Gamma(c-|a|-|b|)}{\Gamma(|b|)\, \Gamma(c-|b|)}\right)\bigg[ \left( \frac{|a|}{(c-|a|-|b|)} \right)\,\sum_{n=0}^{\infty}\, \binom{-(|a|+1)}{n}\,\qquad \qquad\qquad\qquad&&\cr
\times\left(\frac{\Gamma(|b|+4+2n)}{\Gamma(c-|a|-3+2n)}\right)\, _{2}F_1(|a|,|b|+4+2n;c-|a|+3+2n;-1)&&\cr
\,+\lambda\,\sum_{n=0}^{\infty}\, \binom{-|a|}{n}\,\left( \frac{\Gamma(|b|+2n)}{\Gamma(c-|a|+2n)}\right)\,   _{2}F_1(|a|,|b|+2n;c-|a|+2n;-1)\bigg]\nonumber \qquad\qquad&&\cr
\qquad\qquad\qquad\qquad\qquad\qquad\qquad\qquad\qquad \leq \lambda\left(1+\frac{1}{2(1-\beta)}\right).&&
\eeq
Then, the operator $\mathcal{I}^{a,\,\frac{b}{4},\, \frac{b+1}{4},\, \frac{b+2}{4},\, \frac{b+3}{4}}_{\frac{c}{4},\, \frac{c+1}{4},\, \frac{c+2}{4},\, \frac{c+2}{4}}(f)(z)$ maps $ \mathcal{R}(\beta)$ into $ \mathcal{C}_{\lambda}$.
\ethm
\bpf
The proof is similar to Theorem \ref{thm2eq001}. So, we omit the details.
\epf
\bthm\label{thm12eq0} Let, $a,\,b \in {\Bbb C} \backslash \{ 0 \},\,  c > 0$, $c > |a|+|b|+3$, and $ 0 < \lambda \leq 1$. If
 \beq\label{thm12eq1}
&&\left(\frac{\Gamma(c)\,\Gamma(c-|a|-|b|)}{\Gamma(|b|)\, \Gamma(c-|b|)}\right)\bigg[ \left( \frac{(|a|)_3}{(c-|a|-|b|-3)_3} \right)\,\sum_{n=0}^{\infty}\, \binom{-(|a|+3)}{n}\,\cr
&&\qquad\qquad\times\left(\frac{\Gamma(|b|+12+2n)}{\Gamma(c-|a|+9+2n)}\right)\, _{2}F_1(|a|+3,|b|+12+2n;c-|a|+9+2n;-1)\cr
&&\,+\left( \frac{(\lambda+5)(|a|)_2}{(c-|a|-|b|-2)_2} \right)\,\sum_{n=0}^{\infty}\, \binom{-(|a|+2)}{n}\,\qquad \qquad\qquad\qquad\cr
&&\qquad\qquad\times\left(\frac{\Gamma(|b|+8+2n)}{\Gamma(c-|a|+6+2n)}\right)\, _{2}F_1(|a|+2,|b|+8+2n;c-|a|+6+2n;-1)\cr
&&\,+\left( \frac{(3\lambda+4)(|a|)}{(c-|a|-|b|-1)} \right)\,\sum_{n=0}^{\infty}\, \binom{-(|a|+1)}{n}\,\qquad \qquad\qquad\qquad\cr
&&\qquad\qquad\times\left(\frac{\Gamma(|b|+4+2n)}{\Gamma(c-|a|+3+2n)}\right)\, _{2}F_1(|a|+1,|b|+4+2n;c-|a|+3+2n;-1)\cr
&&\qquad\qquad\,+\lambda\,\sum_{n=0}^{\infty}\, \binom{-|a|}{n}\,\left( \frac{\Gamma(|b|+2n)}{\Gamma(c-|a|+2n)}\right)\,   _{2}F_1(|a|,|b|+2n;c-|a|+2n;-1)\bigg] \leq 2\lambda.
\eeq
then, $\mathcal{I}^{a,\,\frac{b}{4},\, \frac{b+1}{4},\, \frac{b+2}{4},\, \frac{b+3}{4}}_{\frac{c}{4},\, \frac{c+1}{4},\, \frac{c+2}{4},\, \frac{c+2}{4}}(f)(z)$ maps $\es$ into $ \mathcal{C}_{\lambda} $.
\ethm
\bpf  Let, $a,\,b \in {\Bbb C} \backslash \{ 0 \},\,  c > 0$, $c > |a|+|b|+3$, and $ 0 < \lambda \leq 1$.

Suppose, the integral operator $\mathcal{I}^{a,\,\frac{b}{4},\, \frac{b+1}{4}\, \frac{b+2}{4}\, \frac{b+3}{4}}_{\frac{c}{4},\, \frac{c+1}{4}\, \frac{c+2}{4}\, \frac{c+2}{4}}(f)(z)$ is defined by (\ref{inteq7}). In view of the sufficient condition given in (\ref{inteq}), it is enough to prove that
\begin{eqnarray*}
  T &:=& \sum_{n=2}^{\infty}\,n\, (n+\lambda-1)\, |A_n|\leq \lambda.
\end{eqnarray*}
Using the fact that $|(a)_n|\leq  (|a|)_n$ and $(\ref{inteq2})$ in the aforementioned equation, it is found that
\begin{eqnarray*}
 T &\leq& \sum_{n=2}^{\infty}n\, (n+\lambda-1)\, \left(\frac{(|a|)_{n-1}\left(\frac{|b|}{4}\right)_{n-1}\, \left(\frac{|b|+1}{4}\right)_{n-1}\, \left(\frac{|b|+2}{4}\right)_{n-1}\, \left(\frac{|b|+3}{4}\right)_{n-1}}{\left(\frac{c}{4}\right)_{n-1}\, \left(\frac{c+1}{4}\right)_{n-1}\, \left(\frac{c+2}{4}\right)_{n-1}\, \left(\frac{c+3}{4}\right)_{n-1}\,(1)_{n-1}}\right)\, |a_n|\\
 &=& \sum_{n=0}^{\infty}\,((n+1)^3+(n+1)^2(\lambda-1))\, \left(\frac{(|a|)_{n}\left(\frac{|b|}{4}\right)_{n}\, \left(\frac{|b|+1}{4}\right)_{n}\, \left(\frac{|b|+2}{4}\right)_{n}\, \left(\frac{|b|+3}{4}\right)_{n}}{\left(\frac{c}{4}\right)_{n}\, \left(\frac{c+1}{4}\right)_{n}\, \left(\frac{c+2}{4}\right)_{n}\, \left(\frac{c+3}{4}\right)_{n}\,(1)_{n}}\right)-\lambda.\\
 &=& \sum_{n=0}^{\infty}\,((n+1)^3\, \left(\frac{(|a|)_{n}\left(\frac{|b|}{4}\right)_{n}\, \left(\frac{|b|+1}{4}\right)_{n}\, \left(\frac{|b|+2}{4}\right)_{n}\, \left(\frac{|b|+3}{4}\right)_{n}}{\left(\frac{c}{4}\right)_{n}\, \left(\frac{c+1}{4}\right)_{n}\, \left(\frac{c+2}{4}\right)_{n}\, \left(\frac{c+3}{4}\right)_{n}\,(1)_{n}}\right)\\
 &&\qquad\qquad+(\lambda-1)\sum_{n=0}^{\infty}\,(n+1)^2\, \left(\frac{(|a|)_{n}\left(\frac{|b|}{4}\right)_{n}\, \left(\frac{|b|+1}{4}\right)_{n}\, \left(\frac{|b|+2}{4}\right)_{n}\, \left(\frac{|b|+3}{4}\right)_{n}}{\left(\frac{c}{4}\right)_{n}\, \left(\frac{c+1}{4}\right)_{n}\, \left(\frac{c+2}{4}\right)_{n}\, \left(\frac{c+3}{4}\right)_{n}\,(1)_{n}}\right)-\lambda.
\end{eqnarray*}
Using (2) and (3) of Lemma \ref{lem1eq1}, it follows that
\begin{eqnarray*}
T &\leq& \left(\frac{\Gamma(c)\,\Gamma(c-|a|-|b|)}{\Gamma(|b|)\, \Gamma(c-|b|)}\right)\bigg[ \left( \frac{(|a|)_3}{(c-|a|-|b|-3)_3} \right)\,\sum_{n=0}^{\infty}\, \binom{-(|a|+3)}{n}\,\cr
&&\qquad\qquad\times\left(\frac{\Gamma(|b|+12+2n)}{\Gamma(c-|a|+9+2n)}\right)\, _{2}F_1(|a|+3,|b|+12+2n;c-|a|+9+2n;-1)\cr
&&\,+\left( \frac{(\lambda+5)(|a|)_2}{(c-|a|-|b|-2)_2} \right)\,\sum_{n=0}^{\infty}\, \binom{-(|a|+2)}{n}\,\qquad \qquad\qquad\qquad\cr
&&\qquad\qquad\times\left(\frac{\Gamma(|b|+8+2n)}{\Gamma(c-|a|+6+2n)}\right)\, _{2}F_1(|a|+2,|b|+8+2n;c-|a|+6+2n;-1)\cr
&&\,+\left( \frac{(3\lambda+4)(|a|)}{(c-|a|-|b|-1)} \right)\,\sum_{n=0}^{\infty}\, \binom{-(|a|+1)}{n}\,\qquad \qquad\qquad\qquad\cr
&&\qquad\qquad\times\left(\frac{\Gamma(|b|+4+2n)}{\Gamma(c-|a|+3+2n)}\right)\, _{2}F_1(|a|+1,|b|+4+2n;c-|a|+3+2n;-1)\cr
&&\qquad\,+\lambda\,\sum_{n=0}^{\infty}\, \binom{-|a|}{n}\,\left( \frac{\Gamma(|b|+2n)}{\Gamma(c-|a|+2n)}\right)\,   _{2}F_1(|a|,|b|+2n;c-|a|+2n;-1)\bigg] -\lambda.
\end{eqnarray*}
By the equation (\ref{thm12eq1}), the above expression is bounded above by $\lambda$, and hence,
\begin{eqnarray*}
&&\left(\frac{\Gamma(c)\,\Gamma(c-|a|-|b|)}{\Gamma(|b|)\, \Gamma(c-|b|)}\right)\bigg[ \left( \frac{(|a|)_3}{(c-|a|-|b|-3)_3} \right)\,\sum_{n=0}^{\infty}\, \binom{-(|a|+3)}{n}\,\cr
&&\qquad\qquad\times\left(\frac{\Gamma(|b|+12+2n)}{\Gamma(c-|a|+9+2n)}\right)\, _{2}F_1(|a|+3,|b|+12+2n;c-|a|+9+2n;-1)\cr \
&&\,+\left( \frac{(\lambda+5)(|a|)_2}{(c-|a|-|b|-2)_2} \right)\,\sum_{n=0}^{\infty}\, \binom{-(|a|+2)}{n}\,\qquad \qquad\qquad\qquad\cr
&&\qquad\qquad\times\left(\frac{\Gamma(|b|+8+2n)}{\Gamma(c-|a|+6+2n)}\right)\, _{2}F_1(|a|+2,|b|+8+2n;c-|a|+6+2n;-1)\cr
&&\,+\left( \frac{(3\lambda+4)(|a|)}{(c-|a|-|b|-1)} \right)\,\sum_{n=0}^{\infty}\, \binom{-(|a|+1)}{n}\,\qquad \qquad\qquad\qquad\cr
&&\qquad\qquad\times\left(\frac{\Gamma(|b|+4+2n)}{\Gamma(c-|a|+3+2n)}\right)\, _{2}F_1(|a|+1,|b|+4+2n;c-|a|+3+2n;-1)\cr
&&\qquad\,+\lambda\,\sum_{n=0}^{\infty}\, \binom{-|a|}{n}\,\left( \frac{\Gamma(|b|+2n)}{\Gamma(c-|a|+2n)}\right)\,   _{2}F_1(|a|,|b|+2n;c-|a|+2n;-1)\bigg] -\lambda \leq \lambda.
\end{eqnarray*}
Hence, the integral operator $\mathcal{I}^{a,\,\frac{b}{4},\, \frac{b+1}{4},\, \frac{b+2}{4},\, \frac{b+3}{4}}_{\frac{c}{4},\, \frac{c+1}{4},\, \frac{c+2}{4},\, \frac{c+2}{4}}(f)(z)$ maps $\es$ into $\mathcal{C}_{\lambda}$ and the proof is complete.
\epf
\subsubsection{Uniformly Convexity.}
\bthm\label{thm7eq1}
  Let, $a,\, b \in {\Bbb C} \backslash \{ 0 \} ,\,c > 0$  and $c > |a|+|b|+2 .$  A sufficient condition for the function $z\, _5F_4\left(^{a,\,\frac{b}{4},\, \frac{b+1}{4}\, \frac{b+2}{4}\, \frac{b+3}{4}}_{\frac{c}{4},\, \frac{c+1}{4}\, \frac{c+2}{4}\, \frac{c+3}{4}};  z\right)$ to belong to the class $\UCV$ is that
 \beq\label{thm7eq10}
&&\left(\frac{\Gamma(c)\,\Gamma(c-|a|-|b|)}{\Gamma(|b|)\, \Gamma(c-|b|)}\right)\bigg[\left( \frac{2\,(|a|)_2}{(c-|a|-|b|-2)_2} \right)\,\sum_{n=0}^{\infty}\, \binom{-(|a|+2)}{n}\,\qquad \qquad\qquad\qquad\cr
&&\qquad\qquad\times\left(\frac{\Gamma(|b|+8+2n)}{\Gamma(c-|a|+6+2n)}\right)\, _{2}F_1(|a|+2,|b|+8+2n;c-|a|+6+2n;-1)\cr
&&\qquad\,+\left( \frac{5(|a|)}{(c-|a|-|b|-1)} \right)\,\sum_{n=0}^{\infty}\, \binom{-(|a|+1)}{n}\,\qquad \qquad\qquad\qquad\cr
&&\qquad\qquad\times\left(\frac{\Gamma(|b|+4+2n)}{\Gamma(c-|a|+3+2n)}\right)\, _{2}F_1(|a|+1,|b|+4+2n;c-|a|+3+2n;-1)\cr
&&\qquad\, \, \qquad\,+\,\sum_{n=0}^{\infty}\, \binom{-|a|}{n}\,\left( \frac{\Gamma(|b|+2n)}{\Gamma(c-|a|+2n)}\right)\,   _{2}F_1(|a|,|b|+2n;c-|a|+2n;-1)\bigg]\leq 2.
\eeq
\ethm
\bpf
Let, $a,\, b \in {\Bbb C} \backslash \{ 0 \} ,\,c > 0$  and $c > |a|+|b|+2 .$

Let, $f(z)=z\, _5F_4\left(^{a,\,\frac{b}{4},\, \frac{b+1}{4}\, \frac{b+2}{4}\, \frac{b+3}{4}}_{\frac{c}{4},\, \frac{c+1}{4}\, \frac{c+2}{4}\, \frac{c+3}{4}};  z\right)$ . Then, by (\ref{lem4eq1}), it is enough to show that
\begin{eqnarray*}
  T &=& \sum_{n=2}^{\infty}\, n\, (2n-1)\,|A_n|\leq 1.
\end{eqnarray*}
where, $A_n$ is given by (\ref{inteq007}). Using the fact $|(a)_n|\leq  (|a|)_n$,
\begin{eqnarray*}
  T &\leq&  \sum_{n=2}^{\infty} n\, (2n-1)\,\left(\frac{(|a|)_{n-1}\left(\frac{|b|}{4}\right)_{n-1}\, \left(\frac{|b|+1}{4}\right)_{n-1}\, \left(\frac{|b|+2}{4}\right)_{n-1}\, \left(\frac{|b|+3}{4}\right)_{n-1}}{\left(\frac{c}{4}\right)_{n-1}\, \left(\frac{c+1}{4}\right)_{n-1}\, \left(\frac{c+2}{4}\right)_{n-1}\, \left(\frac{c+3}{4}\right)_{n-1}(1)_{n-1}}\right)\\
   &=&  2\sum_{n=0}^{\infty} (n+1)^2\, \left(\frac{(|a|)_{n}\left(\frac{|b|}{4}\right)_{n}\, \left(\frac{|b|+1}{4}\right)_{n}\, \left(\frac{|b|+2}{4}\right)_{n}\, \left(\frac{|b|+3}{4}\right)_{n}}{\left(\frac{c}{4}\right)_{n}\, \left(\frac{c+1}{4}\right)_{n}\, \left(\frac{c+2}{4}\right)_{n}\, \left(\frac{c+3}{4}\right)_{n}(1)_{n}}\right)\\
   &&\qquad -\sum_{n=0}^{\infty} (n+1)\, \left(\frac{(|a|)_{n}\left(\frac{|b|}{4}\right)_{n}\, \left(\frac{|b|+1}{4}\right)_{n}\, \left(\frac{|b|+2}{4}\right)_{n}\, \left(\frac{|b|+3}{4}\right)_{n}}{\left(\frac{c}{4}\right)_{n}\, \left(\frac{c+1}{4}\right)_{n}\, \left(\frac{c+2}{4}\right)_{n}\, \left(\frac{c+3}{4}\right)_{n}(1)_{n}}\right)-1.
\end{eqnarray*}
Using (1) and (2) of Lemma \ref{lem1eq1} in the aforementioned  equation, we find that
\begin{eqnarray*}
T &\leq& \left(\frac{\Gamma(c)\,\Gamma(c-|a|-|b|)}{\Gamma(|b|)\, \Gamma(c-|b|)}\right)\bigg[\left( \frac{2\,(|a|)_2}{(c-|a|-|b|-2)_2} \right)\,\sum_{n=0}^{\infty}\, \binom{-(|a|+2)}{n}\,\qquad \qquad\qquad\qquad\cr
&&\qquad\qquad\times\left(\frac{\Gamma(|b|+8+2n)}{\Gamma(c-|a|+6+2n)}\right)\, _{2}F_1(|a|+2,|b|+8+2n;c-|a|+6+2n;-1)\cr
&&\qquad\,+\left( \frac{5(|a|)}{(c-|a|-|b|-1)} \right)\,\sum_{n=0}^{\infty}\, \binom{-(|a|+1)}{n}\,\qquad \qquad\qquad\qquad\cr
&&\qquad\qquad\times\left(\frac{\Gamma(|b|+4+2n)}{\Gamma(c-|a|+3+2n)}\right)\, _{2}F_1(|a|+1,|b|+4+2n;c-|a|+3+2n;-1)\cr
&&\qquad\,+\,\sum_{n=0}^{\infty}\, \binom{-|a|}{n}\,\left( \frac{\Gamma(|b|+2n)}{\Gamma(c-|a|+2n)}\right)\,   _{2}F_1(|a|,|b|+2n;c-|a|+2n;-1)\bigg]-1.
\end{eqnarray*}
Because of (\ref{thm7eq10}), the aforementioned expression is bounded above by 1, and hence,
\begin{eqnarray*}
&&\left(\frac{\Gamma(c)\,\Gamma(c-|a|-|b|)}{\Gamma(|b|)\, \Gamma(c-|b|)}\right)\bigg[\left( \frac{2\,(|a|)_2}{(c-|a|-|b|-2)_2} \right)\,\sum_{n=0}^{\infty}\, \binom{-(|a|+2)}{n}\,\qquad \qquad\qquad\qquad\cr
&&\qquad\qquad\times\left(\frac{\Gamma(|b|+8+2n)}{\Gamma(c-|a|+6+2n)}\right)\, _{2}F_1(|a|+2,|b|+8+2n;c-|a|+6+2n;-1)\cr
&&\qquad\,+\left( \frac{5(|a|)}{(c-|a|-|b|-1)} \right)\,\sum_{n=0}^{\infty}\, \binom{-(|a|+1)}{n}\,\qquad \qquad\qquad\qquad\cr
&&\qquad\qquad\times\left(\frac{\Gamma(|b|+4+2n)}{\Gamma(c-|a|+3+2n)}\right)\, _{2}F_1(|a|+1,|b|+4+2n;c-|a|+3+2n;-1)\cr
&&\qquad\,+\,\sum_{n=0}^{\infty}\, \binom{-|a|}{n}\,\left( \frac{\Gamma(|b|+2n)}{\Gamma(c-|a|+2n)}\right)\,   _{2}F_1(|a|,|b|+2n;c-|a|+2n;-1)\bigg]-1\, \, \, \leq 1.
\end{eqnarray*}
Therefore, $z\, _5F_4\left(^{a,\,\frac{b}{4},\, \frac{b+1}{4}\, \frac{b+2}{4}\, \frac{b+3}{4}}_{\frac{c}{4},\, \frac{c+1}{4}\, \frac{c+2}{4}\, \frac{c+3}{4}}; z\right)$ belongs to the class $ \UCV. $
\epf
\bthm\label{thm8eq0}  Let, $a, \, b \in {\Bbb C} \backslash \{ 0 \} ,\, c > 0,\,$ $c > |a|+|b|+1$ and $0 \leq \beta <1 $. Assume that
 \beq\label{thm8eq1}
&&\left(\frac{\Gamma(c)\,\Gamma(c-|a|-|b|)}{\Gamma(|b|)\, \Gamma(c-|b|)}\right)\nonumber\\
&&\qquad\qquad\times\bigg[\left( \frac{2\,(|a|)}{(c-|a|-|b|-1)} \right)\,\sum_{n=0}^{\infty}\, \binom{-(|a|+1)}{n}\,\left(\frac{\Gamma(|b|+4+2n)}{\Gamma(c-|a|+3+2n)}\right)\nonumber\\
&&\qquad\qquad\qquad\qquad\qquad\qquad\times\, _{2}F_1(|a|,|b|+4+2n;c-|a|+3+2n;-1)\cr
&&\qquad\qquad\qquad\qquad\,+\,\sum_{n=0}^{\infty}\, \binom{-|a|}{n}\,\left( \frac{\Gamma(|b|+2n)}{\Gamma(c-|a|+2n)}\right)\nonumber\\
&&\qquad\qquad\qquad\qquad\qquad\qquad\times\,   _{2}F_1(|a|,|b|+2n;c-|a|+2n;-1)\bigg] \leq \frac{1}{2(1-\beta)}+1.
\eeq
Then, $\mathcal{I}^{a,\,\frac{b}{4},\, \frac{b+1}{4},\, \frac{b+2}{4},\, \frac{b+3}{4}}_{\frac{c}{4},\, \frac{c+1}{4},\, \frac{c+2}{4},\, \frac{c+2}{4}}(f)(z)$  maps $\mathcal{R}(\beta)$ into $ {\UCV}$.
\ethm
\bpf
Let, $a, \, b \in {\Bbb C} \backslash \{ 0 \} ,\, c > 0,\,$ $c > |a|+|b|+1$ and $0 \leq \beta <1 $.

Consider, the integral operator $\mathcal{I}^{a,\,\frac{b}{4},\, \frac{b+1}{4},\, \frac{b+2}{4},\, \frac{b+3}{4}}_{\frac{c}{4},\, \frac{c+1}{4},\, \frac{c+2}{4},\, \frac{c+2}{4}}(f)(z)$  given in (\ref{inteq7}). According to sufficient condition given in (\ref{lem4eq1}), it is enough to show that
\begin{eqnarray*}
  T &:=& \sum_{n=2}^{\infty}n\, (2n-1)\,|A_n|\leq 1,
\end{eqnarray*}
where, $A_n$ is given by (\ref{inteq007}). Using the fact $|(a)_n|\leq  (|a|)_n$ and (\ref{inteq3}) in the aforementioned equation,  it is derived that
\begin{eqnarray*}
  T &\leq&  2(1-\beta)\bigg[\sum_{n=0}^{\infty} \,(2(n+1)-1)\, \left(\frac{(|a|)_{n}\left(\frac{|b|}{4}\right)_{n}\, \, \left(\frac{|b|+1}{4}\right)_{n}\,\left(\frac{|b|+2}{4}\right)_{n}\,\left(\frac{|b|+3}{4}\right)_{n}}{\left(\frac{c}{4}\right)_{n}\, \left(\frac{c+1}{4}\right)_{n}\, \left(\frac{c+2}{4}\right)_{n}\, \left(\frac{c+3}{4}\right)_{n}(1)_{n}}\right)-1\bigg]\\
   &=&  2(1-\beta)\bigg[\sum_{n=0}^{\infty} \,2(n+1)\, \left(\frac{(|a|)_{n}\left(\frac{|b|}{4}\right)_{n}\, \, \left(\frac{|b|+1}{4}\right)_{n}\,\left(\frac{|b|+2}{4}\right)_{n}\,\left(\frac{|b|+3}{4}\right)_{n}}{\left(\frac{c}{4}\right)_{n}\, \left(\frac{c+1}{4}\right)_{n}\, \left(\frac{c+2}{4}\right)_{n}\,\left(\frac{c+3}{4}\right)_{n}(1)_{n}}\right)\cr
   &&\qquad\qquad\qquad-\sum_{n=0}^{\infty} \, \left(\frac{(|a|)_{n}\left(\frac{|b|}{4}\right)_{n}\, \, \left(\frac{|b|+1}{4}\right)_{n}\,\left(\frac{|b|+2}{4}\right)_{n}\,\left(\frac{|b|+3}{4}\right)_{n}}{\left(\frac{c}{4}\right)_{n}\, \left(\frac{c+1}{4}\right)_{n}\, \left(\frac{c+2}{4}\right)_{n}\,\left(\frac{c+3}{4}\right)_{n}(1)_{n}}\right)-1\bigg]
\end{eqnarray*}
Using the formula (\ref{inteq6}) and (1) of Lemma \ref{lem1eq1}, we find that
\begin{eqnarray*}
T&\leq &2(1-\beta)\left(\frac{\Gamma(c)\,\Gamma(c-|a|-|b|)}{\Gamma(|b|)\, \Gamma(c-|b|)}\right)\nonumber\\
&&\qquad\times\bigg[\left( \frac{2\,(|a|)}{(c-|a|-|b|-1)} \right)\,\sum_{n=0}^{\infty}\, \binom{-(|a|+1)}{n}\,\left(\frac{\Gamma(|b|+4+2n)}{\Gamma(c-|a|+3+2n)}\right)\nonumber\\
&&\qquad\qquad\qquad\times\, _{2}F_1(|a|+1,|b|+4+2n;c-|a|+3+2n;-1)\cr
&&\qquad\qquad\,+\,\sum_{n=0}^{\infty}\, \binom{-|a|}{n}\,\left( \frac{\Gamma(|b|+2n)}{\Gamma(c-|a|+2n)}\right)\nonumber\\
&&\qquad\qquad\qquad\qquad\qquad\qquad\times\,   _{2}F_1(|a|,|b|+2n;c-|a|+2n;-1)-1\bigg].
\end{eqnarray*}
By (\ref{thm8eq1}), the aforementioned expression is bounded above by $1$, and hence,
\begin{eqnarray*}
&&2(1-\beta)\left(\frac{\Gamma(c)\,\Gamma(c-|a|-|b|)}{\Gamma(|b|)\, \Gamma(c-|b|)}\right)\nonumber\\
&&\qquad\times\bigg[\left( \frac{2\,(|a|)}{(c-|a|-|b|-1)} \right)\,\sum_{n=0}^{\infty}\, \binom{-(|a|+1)}{n}\,\left(\frac{\Gamma(|b|+4+2n)}{\Gamma(c-|a|+3+2n)}\right)\nonumber\\
&&\qquad\qquad\qquad\qquad\times\, _{2}F_1(|a|+1,|b|+4+2n;c-|a|+3+2n;-1)\cr
&&\qquad\qquad\,+\,\sum_{n=0}^{\infty}\, \binom{-|a|}{n}\,\left( \frac{\Gamma(|b|+2n)}{\Gamma(c-|a|+2n)}\right)\nonumber\\
&&\qquad\qquad\qquad\qquad\qquad\qquad\qquad\times\,   _{2}F_1(|a|,|b|+2n;c-|a|+2n;-1)-1\bigg] \leq 1.
\end{eqnarray*}
Therefore, the operator $\mathcal{I}^{a,\,\frac{b}{4},\, \frac{b+1}{4},\, \frac{b+2}{4},\, \frac{b+3}{4}}_{\frac{c}{4},\, \frac{c+1}{4},\, \frac{c+2}{4},\, \frac{c+2}{4}}(f)(z)$ maps $\mathcal{R}(\beta)$ into $\UCV$, and the result follows.
\epf
\subsubsection{Inclusion Properties in $\es_p$-CLASS}
\bthm\label{thm4eq0}
  Let, $a,\, b \in {\Bbb C} \backslash \{ 0 \} $,\, $c > 0$   and $c > |a|+|b|+1.$  A sufficient condition for the function $z\, _5F_4\left(^{a,\,\frac{b}{4},\, \frac{b+1}{4}\, \frac{b+2}{4}\, \frac{b+3}{4}}_{\frac{c}{4},\, \frac{c+1}{4}\, \frac{c+2}{4}\, \frac{c+3}{4}}; z\right)$ to belong to the class $ \es_p$ is that
  \beq\label{thm4eq1}
&&\left(\frac{\Gamma(c)\,\Gamma(c-|a|-|b|)}{\Gamma(|b|)\, \Gamma(c-|b|)}\right)\nonumber\\
&&\qquad\times\bigg[\left( \frac{2\,(|a|)}{(c-|a|-|b|-1)} \right)\,\sum_{n=0}^{\infty}\, \binom{-(|a|+1)}{n}\,\left(\frac{\Gamma(|b|+4+2n)}{\Gamma(c-|a|+3+2n)}\right)\nonumber\\
&&\qquad\qquad\qquad\qquad\qquad\times\, _{2}F_1(|a|+1,|b|+4+2n;c-|a|+3+2n;-1)\cr
&&\qquad\qquad\,+\,\sum_{n=0}^{\infty}\, \binom{-|a|}{n}\,\left( \frac{\Gamma(|b|+2n)}{\Gamma(c-|a|+2n)}\right)\nonumber\\
&&\qquad\qquad\qquad\qquad\qquad\qquad\qquad\qquad\times\,   _{2}F_1(|a|,|b|+2n;c-|a|+2n;-1)\bigg] \leq 2.
\eeq
\ethm
\bpf
The proof is similar to Theorem \ref{thm8eq0}. So, we omit the details.
\epf
\bthm\label{thm5eq0}  Let, $a,\, b \in {\Bbb C} \backslash \{ 0 \} ,\,|a|\neq1,\, |b|\neq 1,2,3,4,\,   c > 0 ,\,  c > \max\{|a|+3, |a|+|b|-1\}$ and $0 \leq \beta <1$. Assume that
\beq\label{thm5eq1}
&&\left(\frac{\Gamma(c)\,\Gamma(c-|a|-|b|)}{\Gamma(|b|)\, \Gamma(c-|b|)}\right)\bigg(2\,\sum_{n=0}^{\infty}\, \binom{-|a|}{n}\,\left( \frac{\Gamma(|b|+2n)}{\Gamma(c-|a|+2n)}\right)\nonumber\cr
&&\qquad\qquad\qquad\qquad\qquad\qquad\times\,   _{2}F_1(|a|,|b|+2n;c-|a|+2n;-1)\cr
&&\qquad\qquad-\left( \frac{c-|a|-|b|}{(|a|-1)} \right)\,\sum_{n=0}^{\infty}\, \binom{-(|a|-1)}{n}\,\left(\frac{\Gamma(|b|-4+2n)}{\Gamma(c-|a|-3+2n)}\right)\nonumber\cr
&&\qquad\qquad\qquad\qquad\,\,\times\, _{2}F_1(|a|-1,|b|-4+2n;c-|a|-3+2n;-1)\cr
&&\qquad\qquad\qquad\qquad\qquad\qquad\qquad\qquad\qquad\qquad\,+\frac{(c-4)_4}{(|a|-1)(|b|-4)_4}\bigg) \leq \frac{1}{2(1-\beta)}+1.
\eeq
then, $\mathcal{I}^{a,\,\frac{b}{4},\, \frac{b+1}{4},\, \frac{b+2}{4},\, \frac{b+3}{4}}_{\frac{c}{4},\, \frac{c+1}{4},\, \frac{c+2}{4},\, \frac{c+2}{4}}(f)(z)$ maps $\mathcal{R}(\beta)$ into $ S_{p}$ class.
\ethm
\bpf Let, $a,\, b \in {\Bbb C} \backslash \{ 0 \} ,\,|a|\neq1,\, |b|\neq 1,2,3,4,\,   c > 0 ,\,  c > \max\{|a|+3, |a|+|b|-1\}$ and $0 \leq \beta <1$.

Consider, the integral operator $\mathcal{I}^{a,\,\frac{b}{4},\, \frac{b+1}{4},\, \frac{b+2}{4},\, \frac{b+3}{4}}_{\frac{c}{4},\, \frac{c+1}{4},\, \frac{c+2}{4},\, \frac{c+2}{4}}(f)(z)$ given by (\ref{inteq7}). In the view of (\ref{lem2eq1}), it is enough to show that
\begin{eqnarray*}
  T &=& \sum_{n=2}^{\infty}(2n-1)|A_n|\leq 1,
\end{eqnarray*}
where, $A_n$ is given by (\ref{inteq007}). Using the inequality $|(a)_n|\leq  (|a|)_n$, it is proven that
\begin{eqnarray*}
  T &\leq & \sum_{n=2}^{\infty}(2n-1)\bigg(\frac{(|a|)_{n-1}\left(\frac{|b|}{4}\right)_{n-1}\, \left(\frac{|b|+1}{4}\right)_{n-1}\, \left(\frac{|b|+2}{4}\right)_{n-1}\, \left(\frac{|b|+3}{4}\right)_{n-1}}{\left(\frac{c}{4}\right)_{n-1}\, \left(\frac{c+1}{4}\right)_{n-1}\, \left(\frac{c+2}{4}\right)_{n-1}\, \left(\frac{c+3}{4}\right)_{n-1}(1)_{n-1}}\bigg)\, |a_n|\leq 1
\end{eqnarray*}
Using $(\ref{inteq3})$ in the aforesaid equation, one can have
\begin{eqnarray*}
  T  &\leq &  2(1-\beta)\bigg(2\sum_{n=0}^{\infty} \frac{(|a|)_{n}\left(\frac{|b|}{4}\right)_{n}\, \left(\frac{|b|+1}{4}\right)_{n}\, \left(\frac{|b|+2}{4}\right)_{n}\, \left(\frac{|b|+3}{4}\right)_{n}}{\left(\frac{c}{4}\right)_{n}\, \left(\frac{c+1}{4}\right)_{n}\, \left(\frac{c+2}{4}\right)_{n}\, \left(\frac{c+3}{4}\right)_{n}(1)_{n}}\cr
&&\qquad-\sum_{n=0}^{\infty}  \frac{(|a|)_{n}\left(\frac{|b|}{4}\right)_{n}\, \left(\frac{|b|+1}{4}\right)_{n}\, \left(\frac{|b|+2}{4}\right)_{n}\, \left(\frac{|b|+2}{4}\right)_{n}\, \left(\frac{|b|+3}{4}\right)_{n}}{\left(\frac{c}{4}\right)_{n}\, \left(\frac{c+1}{4}\right)_{n}\, \left(\frac{c+2}{4}\right)_{n}\, \left(\frac{c+3}{4}\right)_{n}(1)_{n+1}} -1  \bigg)
\end{eqnarray*}
Using the equation (\ref{inteq6}), and the result (4) of Lemma \ref{lem1eq1}, it is found that
\begin{eqnarray*}
T &\leq&  2(1-\beta)\left(\frac{\Gamma(c)\,\Gamma(c-|a|-|b|)}{\Gamma(|b|)\, \Gamma(c-|b|)}\right)\bigg(2\,\sum_{n=0}^{\infty}\, \binom{-|a|}{n}\,\left( \frac{\Gamma(|b|+2n)}{\Gamma(c-|a|+2n)}\right)\nonumber\cr
&&\qquad\qquad\qquad\qquad\qquad\qquad\times\,   _{2}F_1(|a|,|b|+2n;c-|a|+2n;-1)\cr
&&\qquad\qquad-\left( \frac{c-|a|-|b|}{(|a|-1)} \right)\,\sum_{n=0}^{\infty}\, \binom{-(|a|-1)}{n}\,\left(\frac{\Gamma(|b|-4+2n)}{\Gamma(c-|a|-3+2n)}\right)\nonumber\cr
&&\qquad\qquad\qquad\qquad\,\,\times\, _{2}F_1(|a|-1,|b|-4+2n;c-|a|-3+2n;-1)\cr
&&\qquad\qquad\qquad\qquad\qquad\qquad\qquad\qquad\qquad\qquad\,+\frac{(c-4)_4}{(|a|-1)(|b|-4)_4}-1\bigg).
\end{eqnarray*}
By the condition (\ref{thm5eq1}), the aforementioned expression is bounded above by 1, and hence
\begin{eqnarray*}
&& 2(1-\beta)\left(\frac{\Gamma(c)\,\Gamma(c-|a|-|b|)}{\Gamma(|b|)\, \Gamma(c-|b|)}\right)\bigg(2\,\sum_{n=0}^{\infty}\, \binom{-|a|}{n}\,\left( \frac{\Gamma(|b|+2n)}{\Gamma(c-|a|+2n)}\right)\nonumber\cr
&&\qquad\qquad\qquad\qquad\qquad\qquad\times\,   _{2}F_1(|a|,|b|+2n;c-|a|+2n;-1)\cr
&&\qquad\qquad-\left( \frac{c-|a|-|b|}{(|a|-1)} \right)\,\sum_{n=0}^{\infty}\, \binom{-(|a|-1)}{n}\,\left(\frac{\Gamma(|b|-4+2n)}{\Gamma(c-|a|-3+2n)}\right)\nonumber\cr
&&\qquad\qquad\qquad\qquad\,\,\times\, _{2}F_1(|a|-1,|b|-4+2n;c-|a|-3+2n;-1)\cr
&&\qquad\qquad\qquad\qquad\qquad\qquad\qquad\qquad\qquad\qquad\,+\frac{(c-4)_4}{(|a|-1)(|b|-4)_4}-1\bigg) \leq 1.
\end{eqnarray*}
Under the stated condition, the operator $\mathcal{I}^{a,\,\frac{b}{4},\, \frac{b+1}{4},\, \frac{b+2}{4},\, \frac{b+3}{4}}_{\frac{c}{4},\, \frac{c+1}{4},\, \frac{c+2}{4},\, \frac{c+2}{4}}(f)(z)$ maps $\mathcal{R}(\beta)$ into $\es_p$ and  the proof is complete.
\epf
\bthm\label{thm6eq0}
 Let, $a,\, b \in {\Bbb C} \backslash \{ 0 \} ,\, c > 0,\,$ and $c > |a|+|b|+2$. Suppose, $a,\, b$, and $c$ satisfy the condition
 \beq\label{thm6eq1}
&&\left(\frac{\Gamma(c)\,\Gamma(c-|a|-|b|)}{\Gamma(|b|)\, \Gamma(c-|b|)}\right)\bigg[\left( \frac{2\,(|a|)_2}{(c-|a|-|b|-2)_2} \right)\,\sum_{n=0}^{\infty}\, \binom{-(|a|+2)}{n}\,\qquad \qquad\qquad\qquad\cr
&&\qquad\qquad\times\left(\frac{\Gamma(|b|+8+2n)}{\Gamma(c-|a|+6+2n)}\right)\, _{2}F_1(|a|+2,|b|+8+2n;c-|a|+6+2n;-1)\cr
&&\qquad\,+\left( \frac{5(|a|)}{(c-|a|-|b|-1)} \right)\,\sum_{n=0}^{\infty}\, \binom{-(|a|+1)}{n}\,\qquad \qquad\qquad\qquad\cr
&&\qquad\qquad\times\left(\frac{\Gamma(|b|+4+2n)}{\Gamma(c-|a|+3+2n)}\right)\, _{2}F_1(|a|+1,|b|+4+2n;c-|a|+3+2n;-1)\cr
&&\qquad\, \, \qquad\,-\,\sum_{n=0}^{\infty}\, \binom{-|a|}{n}\,\left( \frac{\Gamma(|b|+2n)}{\Gamma(c-|a|+2n)}\right)\,   _{2}F_1(|a|,|b|+2n;c-|a|+2n;-1)\bigg]\leq 2.\nonumber
\eeq
Then, the operator $\mathcal{I}^{a,\,\frac{b}{4},\, \frac{b+1}{4},\, \frac{b+2}{4},\, \frac{b+3}{4}}_{\frac{c}{4},\, \frac{c+1}{4},\, \frac{c+2}{4},\, \frac{c+2}{4}}(f)(z)$ maps $\es$ into $ \es_{p} $ class.
\ethm
\bpf
The proof is similar to Theorem \ref{thm7eq1}. So, we omit the details.
\epf

\end{sloppypar}
\end{document}